\documentclass{mscreport} 

\title{Solving Schr\"{o}dinger's equation \\by B-spline collocation} 
\author{Christian P. H. Salas} 
\pin{} 
\begin{document}

\maketitle

\begin{abstract}
B-splines and collocation techniques have been applied to the solution of Schr\"{o}dinger's equation in quantum mechanics since the early 1970s, but one aspect that is noticeably missing from this literature is the use of \emph{Gaussian points} (i.e., the zeros of Legendre polynomials) as the collocation points, which can significantly reduce approximation errors. Authors in the past have used equally spaced or nonlinearly distributed collocation points (noticing that the latter can increase approximation accuracy) but, strangely, have continued to avoid Gaussian collocation points so there are no published papers employing this approach. Using the methodology and computer routines provided by Carl de Boor's book \emph{A Practical Guide to Splines} as a `numerical laboratory', the present dissertation examines how the use of Gaussian collocation points can interact with other features such as box size, mesh size and the order of polynomial approximants to affect the accuracy of approximations to Schr\"{o}dinger's bound state wave functions for the electron in the hydrogen atom. In particular, we explore whether or not, and under what circumstances, B-spline collocation at Gaussian points can produce more accurate approximations to Schr\"{o}dinger's wave functions than equally spaced and nonlinearly distributed collocation points. We also apply B-spline collocation at Gaussian points to a Schr\"{o}dinger equation with cubic nonlinearity which has been used extensively in the past to study nonlinear phenomena. Our computer experiments show that in the case of the hydrogen atom, collocation at Gaussian points can be a highly successful approach,  consistently superior to equally spaced collocation points and often superior to nonlinearly distributed collocation points. However, we do encounter some situations, typically when the mesh is quite coarse relative to the box size for the hydrogen atom, and also in the cubic Schr\"{o}dinger equation case, in which nonlinearly distributed collocation points perform significantly better than Gaussian collocation points.         
\end{abstract}

\tableofcontents


\chapter{Introduction}

Numerous studies have demonstrated the usefulness of B-splines and collocation techniques for solving approximation problems in quantum mechanics (see, e.g., \cite{shore}, \cite{bottcher}, \cite{sapirstein}, \cite{morrison}, \cite{sanchez}, \cite{martin}, \cite{odero}, \cite{bachau}, \cite{ting}, \cite{chang}). B-spline collocation techniques involve the use of spline basis functions to construct piecewise polynomial approximations to the solutions of differential equations in such a way that the approximations are guaranteed to satisfy the differential equations at certain collocation points within subintervals of the domain of interest. The literature applying these techniques to atomic theory began with a seminal paper by Bruce Shore published in 1973 \cite{shore}. He showed how cubic spline collocation could be used to solve the radial Schr\"{o}dinger equation as an eigenvalue problem in a variety of spherically symmetric quantum systems with zero angular momentum. In particular, he solved Schr\"{o}dinger's radial equation for the hydrogen atom with a Coulomb potential, comparing cubic spline collocation with Galerkin methods and finding the latter somewhat superior. 

Decades later, authors were still revisiting and extending Shore's results (see \cite{sapirstein} and \cite{odero} in particular) and there is now a large literature encompassing a wide range of non-relativisic and relativistic quantum mechanical applications of B-splines and collocation. However, one aspect that is noticeably missing from this literature is the use of \emph{Gaussian points} (i.e., the zeros of Legendre polynomials) as the collocation points. Shore's paper was published a short time before another influential paper appeared in a numerical analysis journal in 1973, written by Carl de Boor and Blair Swartz \cite{deboor}, showing how collocation at Gaussian points can significantly reduce approximation errors. This approach exploits the orthogonality of Legendre polynomials to make some of the polynomial products making up the relevant Green's functions vanish, thus reducing the norm of the error terms particularly at the boundaries of the subintervals (a phenomenon called `superconvergence'). A later book by Carl de Boor published in 1978, \emph{A Practical Guide to Splines} \cite{deboor2}, made these ideas much more widely accessible by providing practical advice and relevant computer routines. Shore made no mention of collocation at Gaussian points in his 1973 paper, though he did emphasise that changing from equally spaced collocation points to nonlinearly distributed collocation points improved the accuracy of his approximations by several orders of magnitude. Strangely, authors who have revisited Shore's work, even quite recently, have continued to choose not to employ Gaussian points in their collocation approaches (see, e.g., \cite{odero}) preferring to use equally spaced or nonlinearly distributed collocation points instead.

It has become customary to use known solutions of Schr\"{o}dinger's equation, particularly for the hydrogen atom, as the prototypical problems with which to explore approximation methods in quantum mechanics. As there currently seems to be no published work in which collocation at Gaussian points is explored in this context, the aim of the present dissertation is to address this gap in the literature by revisiting and extending the work on the hydrogen atom in Shore's paper, and thoroughly studying how the use of Gaussian collocation points can interact with other features such as box size, mesh size and the order of polynomial approximants to affect the accuracy of approximations to Schr\"{o}dinger's wave functions. In particular, the dissertation will seek to determine whether or not, and under what circumstances, B-spline collocation at Gaussian points can produce more accurate approximations to Schr\"{o}dinger's wave functions than equally spaced and nonlinearly distributed collocation points. As in Shore's paper, bound state wave functions (negative energy) for the electron in the hydrogen atom will be studied using a Coulomb potential, but we will also extend Shore's framework by studying radial Schr\"{o}dinger equations for the hydrogen atom incorporating nonzero angular momentum.

The dissertation will also explore the applicability of B-spline collocation at Gaussian points to a particular \emph{nonlinear} extension of the radial equation in Shore's paper, which actually arises from a Schr\"{o}dinger equation with cubic nonlinearity and a potential. This form of nonlinear Schr\"{o}dinger equation was first proved to have standing wave solutions in a 1986 paper by Floer and Weinstein \cite{floer} and has been used extensively since then to study solitions and other nonlinear phenomena in areas such as optics, plasma physics, superconductivity and quantum field theory. The standing wave solutions arise when a certain perturbation parameter in the equation is close enough to zero and this setup seems somewhat similar to the nonlinear perturbation problem discussed in Chapter XV of \emph{A Practical Guide to Splines}. This nonlinear extension of the radial equation in Shore's paper therefore seems well worth exploring here, not only being well-suited to the machinery of de Boor's book, but also due to the fact that no previous use appears to have been made of collocation at Gaussian points in this literature. For a discussion of exact solutions of the Schr\"{o}dinger equation with cubic nonlinearity, see \cite{ebaid}, and for additional references and a discussion of the asymptotic behaviour of the solutions in the presence of a potential, see \cite{naumkin}.         

The strategy in this dissertation will be to use the methodology and computer routines provided by Carl de Boor's book \emph{A Practical Guide to Splines}, particularly the setup in Chapter XV, as a kind of `numerical laboratory' to explore the extent to which collocation at Gaussian points is feasible and can accurately approximate Schr\"{o}dinger wave functions. We will therefore treat the problem as a two-point BVP with all the parameters and exact solutions known, and experiments will then be carried out using different patterns of collocation points to investigate the effects on approximating the \emph{eigenfunctions} of Schr\"{o}dinger's equation accurately. Note that this approach is different from Shore's in that he focused primarily on finding individual \emph{eigenvalues} for Schr\"{o}dinger's equation, assuming these are unknown a priori. Looking at entire eigenfunction approximations, rather than single eigenvalues, will provide richer visual and numerical information for studying the detailed effects of varying the pattern of collocation points in conjunction with different box sizes, mesh sizes, and orders of the polynomial approximants. 

The dissertation is organised as follows. Chapter 2 derives the radial Schr\"{o}dinger differential equation used in Shore's paper, supported by detailed mathematical notes in Appendix~\ref{first-appendix} and Appendix~\ref{second-appendix}. It also explains how Shore's framework  for the hydrogen atom, and its extension to cases with nonzero angular momentum and to the nonlinear Schr\"{o}dinger equation, will be implemented in the computer experiments. It concludes by clarifying how our eigenfunction approach differs from Shore's eigenvalue approach. Chapter 3 provides the necessary background on B-splines and the concepts relating to collocation at Gaussian points in the de Boor and Swartz paper. The material here is tailored to Shore's radial equation, in particular to clarify how the choice of Gaussian collocation points can improve approximations in this particular context. Chapter 4 reports the results for bound state electronic wave functions (negative energy) in the hydrogen atom, while Chapter 5 reports the results for the nonlinear extension of Shore's radial equation relating to the Schr\"{o}dinger equation with cubic nonlinearity. Finally, Chapter 6 summarises and evaluates the findings of the dissertation, and suggests possible directions for future investigations. The key components of the computer routines used in the dissertation are provided in Appendices C to G. 


\chapter{Derivation and implementation of the equations in this study}

\section{Shore's radial equation for the hydrogen atom}
The differential equation used in Shore's paper is ultimately based on Schr\"{o}dinger's time-dependent equation
\begin{equation*}
i \hbar \frac{\partial \Psi}{\partial t} = H \Psi
\end{equation*}
where $H$ is the standard Hamiltonian
\begin{equation*}
H = -\frac{\hbar^2}{2m}\nabla^2 + U
\end{equation*}
and $m$ and $U$ are the mass of the quantum particle in question and the potential energy of the system respectively (see, e.g., \cite{griffiths}). In the case of the hydrogen atom with a Coulomb potential, for example, $m = m_e$ is the mass of the electron and the potential is
\begin{equation*}
U = -\frac{e^2}{4 \pi \epsilon_0 r}
\end{equation*}
where $e$ is the electronic charge, $\epsilon_0$ is the permittivity of free space and $r$ is the radial distance of the electron from the nucleus. By separation of variables, Schr\"{o}dinger's time-dependent equation is decomposed into a time-independent equation
\begin{equation}
H \psi = E \psi
\end{equation}
and an essentially trivial differential equation involving time whose solution is an exponential function of time and the parameter $E$. In (2.1), $\psi$ is a time-independent wave function representing a stationary quantum state, or eigenstate, of the system and $E$ in the case of a bound electron is the energy eigenvalue corresponding to this particular eigenstate. The solution $\Psi$ for Schr\"{o}dinger's time-dependent equation is then written as a superposition of products of the form
\begin{equation*}
\psi e^{-iEt/\hbar}
\end{equation*}       
such that this superposition contains all possible eigenstate-eigenvalue pairs. Observation of the system causes this superposition to collapse to one particular eigenstate, with the probability of observing that state being proportional to the modulus squared of its expansion coefficient in the superposition. 

Solving a bound-state quantum mechanics problem essentially involves finding the eigenvalues E and corresponding eigenstates $\psi$ of the time-independent equation (2.1) above, given the functional form of the potential energy $U$. In Appendix~\ref{first-appendix}, I provide a full derivation of the time-independent wave function for the electron in a hydrogen atom, which takes the form 
\begin{equation*}
\psi_{\tilde{n} l m_l}(r, \theta, \phi) \propto e^{-\rho/2}p^l L_{\tilde{n} - l - 1}^{(2l + 1)} P_l^{m_l}(\cos \theta) e^{i m_l \phi}
\end{equation*}
where
\begin{equation*}
\rho = \bigg(-\frac{8 m_e E_{\tilde{n}}}{\hbar^2}\bigg)^{1/2} r
\end{equation*}
and
\begin{equation*}
E_{\tilde{n}} = \bigg(-\frac{m_e}{2 \hbar^2}\bigg) \bigg(\frac{e^2}{4 \pi \epsilon_0}\bigg)^2 \frac{1}{\tilde{n}^2}
\end{equation*}
and where $\tilde{n} = 1, 2, 3, \ldots$ is the principal quantum number determining the electron's energy, $l = 0, 1, 2, \ldots, (\tilde{n} - 1)$ is the orbital quantum number determining its orbital angular-momentum magnitude, $m_l = 0, \pm 1, \pm 2, \ldots, \pm l$ is the magnetic quantum number determining its orbital angular-momentum direction, $L_{\tilde{n} - l - 1}^{(2l + 1)}$ are the associated Laguerre polynomials and $P_l^{m_l}(\cos \theta)$ are the associated Legendre functions, all of which are discussed in Appendix~\ref{first-appendix}.

\subsection{The radial part of the wave function}

Typically in numerical studies involving the bound states of the electron in the hydrogen atom, we are concerned only with the discrete bound states produced by Coulomb attraction in the radial direction, so we restrict our attention to the radial differential equation in Appendix~\ref{first-appendix}, namely
\begin{equation}
\frac{1}{r^2} \frac{d}{d r}\bigg( r^2 \frac{d R}{d r}\bigg) + \bigg[ \frac{2m_e}{\hbar^2}\bigg(\frac{e^2}{4 \pi \epsilon_0 r} + E\bigg) - \frac{l(l + 1)}{r^2} \bigg] R = 0
\end{equation}
whose solutions are 
\begin{equation*}
R_{\tilde{n}l}(r) \propto e^{-\rho/2}\rho^l L_{\tilde{n} - l - 1}^{(2l + 1)}
\end{equation*}
In Appendix~\ref{second-appendix}, I provide detailed derivations of the first few exact solutions of the radial Schr\"{o}dinger differential equation based on this formula, for use in assessing the accuracy of the approximations in this study. 

In Shore's paper the situation is restricted still further in that he only considers the spherically symmetric case in which wave functions have no dependence on angle whatsoever. These wave functions therefore have angular momentum quantum numbers $l = m_l = 0$ and under these circumstances the full wave function above reduces to
\begin{equation}
R_{\tilde{n}0}(r) \propto e^{-\rho/2} L_{\tilde{n} - 1}^{(1)}
\end{equation}
These are the solutions to the radial Schr\"{o}dinger equation
\begin{equation}
\frac{1}{r^2} \frac{d}{d r}\bigg( r^2 \frac{d R}{d r}\bigg) + \bigg[ \frac{2m_e}{\hbar^2}\bigg(\frac{e^2}{4 \pi \epsilon_0 r} + E\bigg)\bigg] R = 0
\end{equation}
obtained by setting $l = 0$ in (2.2). The differential equation used in Shore's paper is just a rescaled version of (2.4), resulting from expressing radial distances from the nucleus in terms of the \emph{Bohr radius}
\begin{equation}
a = \frac{4\pi\epsilon_0 \hbar^2}{e^2 m_e}
\end{equation}
(This is the radius of the innermost Bohr orbit, equal to $5.292 \times 10^{-11}$m). To see this, we can derive Shore's equation (equation (I.1) in his paper) directly from (2.4) as follows. Let 
\begin{equation*}
F(r) = rR
\end{equation*}
Then the first term in (2.4) becomes
\begin{equation*}
\frac{1}{r^2} \frac{d}{d r}\bigg( r^2 \frac{d }{d r}\bigg(\frac{F}{r}\bigg)\bigg) = \frac{1}{r} \frac{d^2F}{dr^2} 
\end{equation*}
so we can rewrite equation (2.4) as
\begin{equation*}
\frac{1}{r} \frac{d^2F}{dr^2} + \bigg[ \frac{2m_e}{\hbar^2}\bigg(\frac{e^2}{4 \pi \epsilon_0 r} + E\bigg)\bigg] \frac{F}{r} = 0
\end{equation*}
Multiplying through by $-\frac{\hbar^2 r}{2m_e}$ and rearranging we get
\begin{equation}
-\frac{\hbar^2}{2m_e} \frac{d^2F}{dr^2} - \frac{e^2}{4 \pi \epsilon_0 r}F = EF
\end{equation}
We can now make the change of variable $r = ax$ where $a$ is the Bohr radius defined in (2.5) above. We then have $dr^2 = a^2 dx^2$ and putting this in (2.6) we get
\begin{equation*}
-\frac{\hbar^2}{2m_e}\frac{1}{a^2} \frac{d^2F}{dx^2} - \frac{e^2}{4 \pi \epsilon_0 a}\frac{F}{x} = EF 
\end{equation*}
or
\begin{equation*}
-\frac{\hbar^2}{2m_e}\bigg(\frac{e^4 m_e^2}{(4 \pi \epsilon_0)^2 \hbar^4}\bigg) \frac{d^2F}{dx^2} - \bigg(\frac{e^2}{4 \pi \epsilon_0}\bigg)\bigg(\frac{e^2 m_e}{4 \pi \epsilon_0 \hbar^2}\bigg)\frac{F}{x} = EF 
\end{equation*}
which simplifies to
\begin{equation}
\frac{1}{2} \frac{d^2F}{dx^2} + \bigg[\frac{1}{x} + E^{\prime}\bigg]F = 0
\end{equation}
where 
\begin{equation}
E^{\prime} = \frac{(4\pi\epsilon_0)^2 \hbar^2}{e^4 m_e} E
\end{equation}
Equation (2.7) is equation (I.1) in Shore's paper, with the rescaled Coulomb potential $V(x) = -\frac{1}{x}$ and the rescaled energy $E^{\prime}$. To obtain $E^{\prime}$ explicitly, note that at the end of Appendix~\ref{first-appendix} we found that the unscaled energy for the hydrogen atom problem is given by
\begin{equation*}
E = \bigg(-\frac{m_e}{2 \hbar^2}\bigg) \bigg(\frac{e^2}{4 \pi \epsilon_0}\bigg)^2 \frac{1}{\tilde{n}^2}
\end{equation*}
Putting this in (2.8) yields the rescaled energy as
\begin{equation}
E^{\prime} = \bigg(\frac{(4\pi\epsilon_0)^2 \hbar^2}{e^4 m_e}\bigg) \bigg(-\frac{m_e}{2 \hbar^2}\bigg) \bigg(\frac{e^2}{4 \pi \epsilon_0}\bigg)^2 \frac{1}{\tilde{n}^2} = -\frac{1}{2\tilde{n}^2}
\end{equation}
Therefore, for example, the ground state energy for Shore's rescaled equation (corresponding to $\tilde{n} = 1$) is $-\frac{1}{2}$. 

\subsection{Implementation in computer experiments}
\begin{figure}[ht]
\begin{center}
\includegraphics[width=4in]{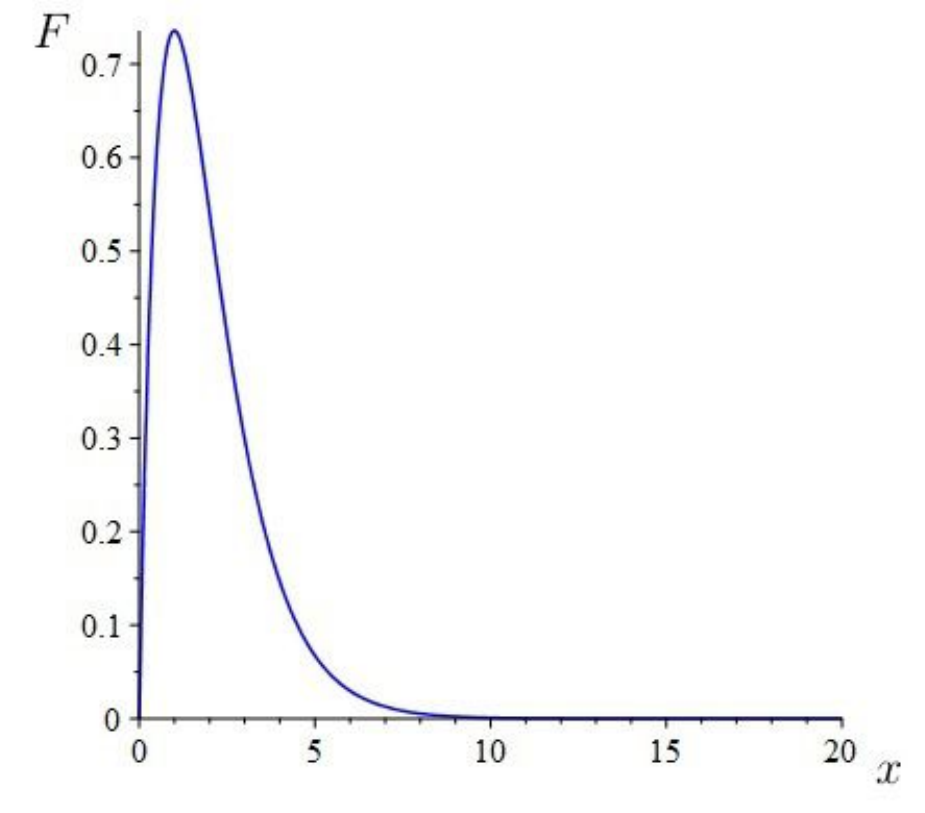}
\caption{Exact solution to Shore's radial equation for the electronic ground state in the hydrogen atom}\label{fig:R10Solution}
\end{center}
\end{figure}

The first equation to be implemented in our study is Shore's radial equation for the ground state of the electron in the hydrogen atom, using the rescaled Coulomb potential $V(x) = -\frac{1}{x}$ and the rescaled energy $-\frac{1}{2}$, giving a radial equation of the form
\begin{equation}
\frac{1}{2} \frac{d^2F}{dx^2} + \bigg[\frac{1}{x} - \frac{1}{2}\bigg]F = 0
\end{equation}
Since we obtained Shore's equation by making the change of variable $F(r) = rR$ in the unscaled radial equation, and by rescaling distances so that they are all expressed in terms of the Bohr radius $a$, the solutions to Shore's equation will be of the form $r R_{\tilde{n}0}$ where $R_{\tilde{n}0}$ is as given in (2.3) above, but with $a=1$ whenever $a$ arises in these solutions. Using equation (B.1) in Appendix~\ref{second-appendix}, the exact solution to (2.10) is then given by applying these changes to $rR_{10}$ to get
\begin{equation}
F(x) = 2xe^{-x}
\end{equation}
This is the exact solution we can use to gauge the accuracy of our computer approximations for the ground state of the electron in the hydrogen atom. Figure~\ref{fig:R10Solution} shows a plot of (2.11). 

Since we are assuming all parameters are known, we will use the boundary conditions $F^{\prime}(0) = 2$ and $F(\infty) = 0$, implementing the latter by ensuring that the box size is large enough to approximate this condition adequately at the right-hand endpoint of the interval. The relevance of box size to the accuracy of approximations is a feature that will be explored in the dissertation. The boundary condition at $0$ comes from the known solution in (2.11). Note that in his paper Shore used the boundary condition $F(0) = 0$ but, as explained in section 3.2.1 below, when applying the approach in Chapter XV of de Boor's book \emph{A Practical Guide to Splines} this causes the collocation procedure to find only the trivial solution $F(x) = 0$. For our numerical work in this dissertation in which we are focusing only on the relative performance of different patterns of collocation points assuming everything else is known, setting the first boundary condition as $F^{\prime}(0) = 2$ ensures that the exact solution in (2.11) is found. 

\begin{figure}[ht]
\begin{center}
\includegraphics[width=4in]{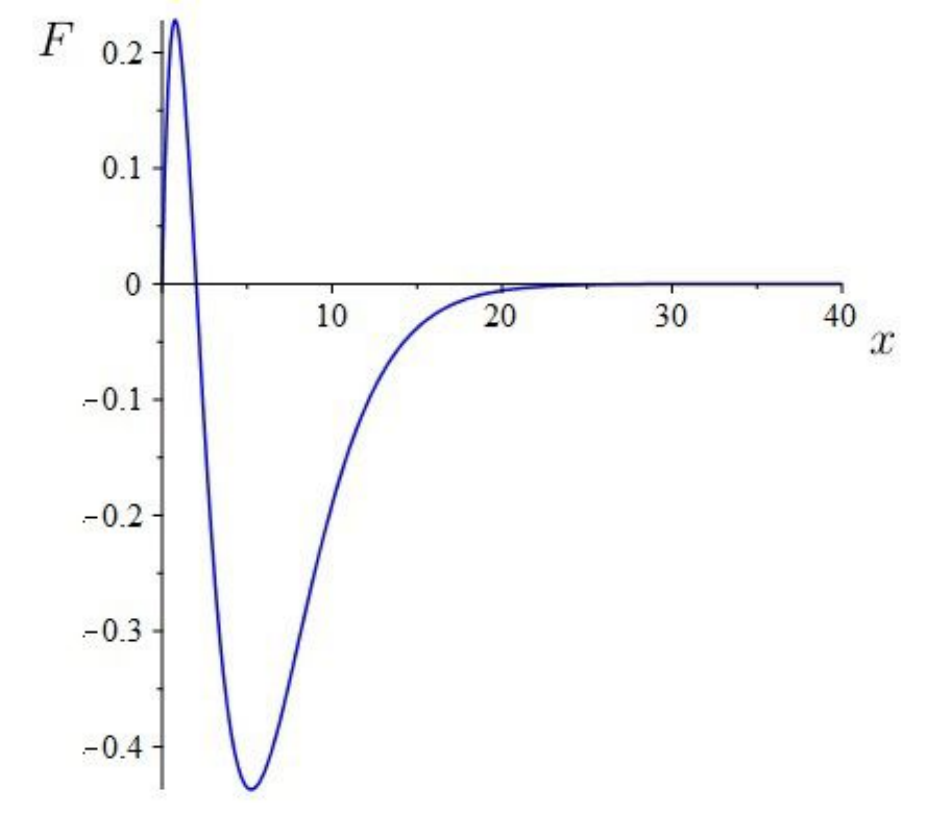}
\caption{Exact solution to Shore's radial equation for $\tilde{n} = 2$ in the hydrogen atom}\label{fig:R20Solution}
\end{center}
\end{figure}

We next implemented Shore's radial equation for the first excited state of the electron in the hydrogen atom. From (2.9), the rescaled energy for the first excited state corresponding to  $\tilde{n} = 2$ is $-\frac{1}{8}$, giving a radial equation of the form
\begin{equation}
\frac{1}{2} \frac{d^2F}{dx^2} + \bigg[\frac{1}{x} - \frac{1}{8}\bigg]F = 0
\end{equation} 
Using equation (B.2) in Appendix~\ref{second-appendix}, the exact solution to (2.12) is then given by setting $a=1$ in $rR_{20}$ to get
\begin{equation}
F(x) = \frac{1}{2\sqrt{2}}x(2 - x)e^{-x/2}
\end{equation}
Figure~\ref{fig:R20Solution} shows a plot of (2.13). 

In this case we use the boundary conditions $F^{\prime}(0) = \frac{1}{\sqrt{2}}$ and $F(\infty) = 0$ for the purposes of our experiments with different patterns of collocation points, again implementing the latter by ensuring that the box size is large enough to approximate this condition. As before, the boundary condition at $0$ comes from the known solution in (2.13).

\begin{figure}[ht]
\begin{center}
\includegraphics[width=4in]{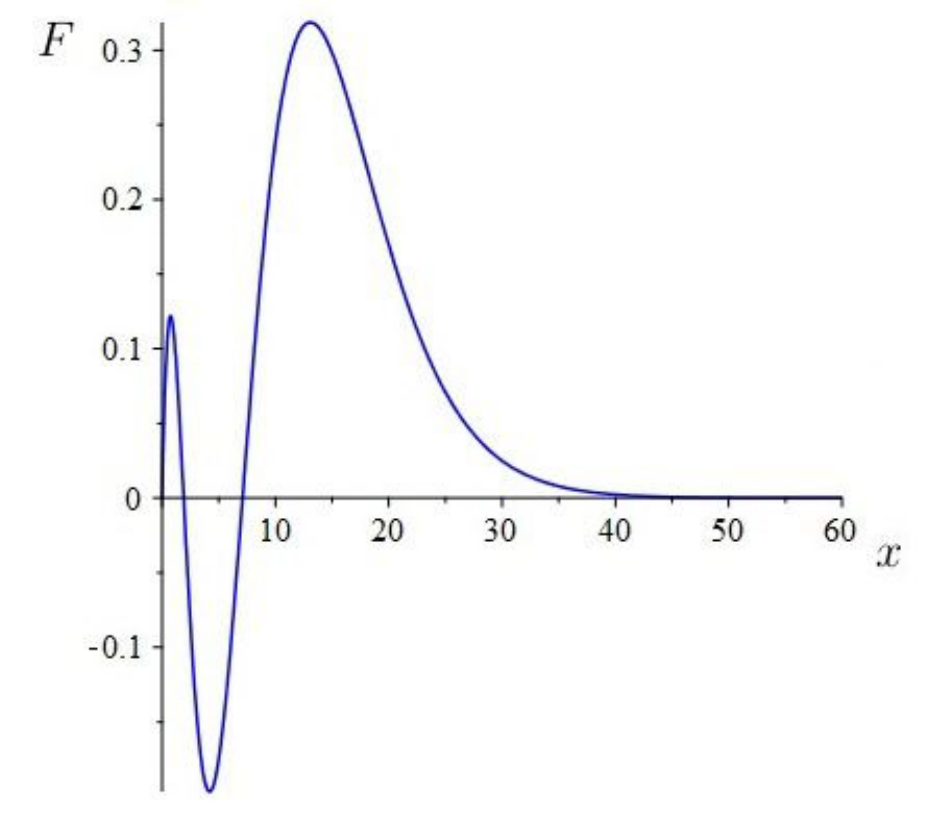}
\caption{Exact solution to Shore's radial equation for $\tilde{n} = 3$ in the hydrogen atom}\label{fig:R30Solution}
\end{center}
\end{figure}

Finally for the zero angular momentum case, we implemented Shore's radial equation for the second excited state of the electron in the hydrogen atom. From (2.9), the rescaled energy for the second excited state corresponding to  $\tilde{n} = 3$ is $-\frac{1}{18}$, giving a radial equation of the form
\begin{equation}
\frac{1}{2} \frac{d^2F}{dx^2} + \bigg[\frac{1}{x} - \frac{1}{18}\bigg]F = 0
\end{equation} 
Using equation (B.3) in Appendix~\ref{second-appendix}, the exact solution to (2.14) is then given by setting $a=1$ in $rR_{30}$ to get
\begin{equation}
F(x) = \frac{2}{81\sqrt{3}}x(27 - 18x + 2x^2)e^{-x/3}
\end{equation}
Figure~\ref{fig:R30Solution} shows a plot of (2.15). 

In this case we use the boundary conditions $F^{\prime}(0) = \frac{2}{3\sqrt{3}}$ and $F(\infty) = 0$ for our experiments with different patterns of collocation points.

\subsection{Incorporating angular momentum}

To extend the $l = 0$ case in Shore's paper, we can also implement radial equations with nonzero orbital angular momentum obtained by rescaling (2.2) in exactly the same way that we rescaled (2.4) earlier, to give
\begin{equation}
\frac{1}{2} \frac{d^2F}{dx^2} + \bigg[\frac{1}{x} + E^{\prime} -\frac{l(l + 1)}{2x^2}\bigg]F = 0
\end{equation}
where $E^{\prime}$ is as in (2.9) above and $l \leq \tilde{n} - 1$ (cf. \cite{odero}, p. 1098). 
\begin{figure}[ht]
\begin{center}
\includegraphics[width=4in]{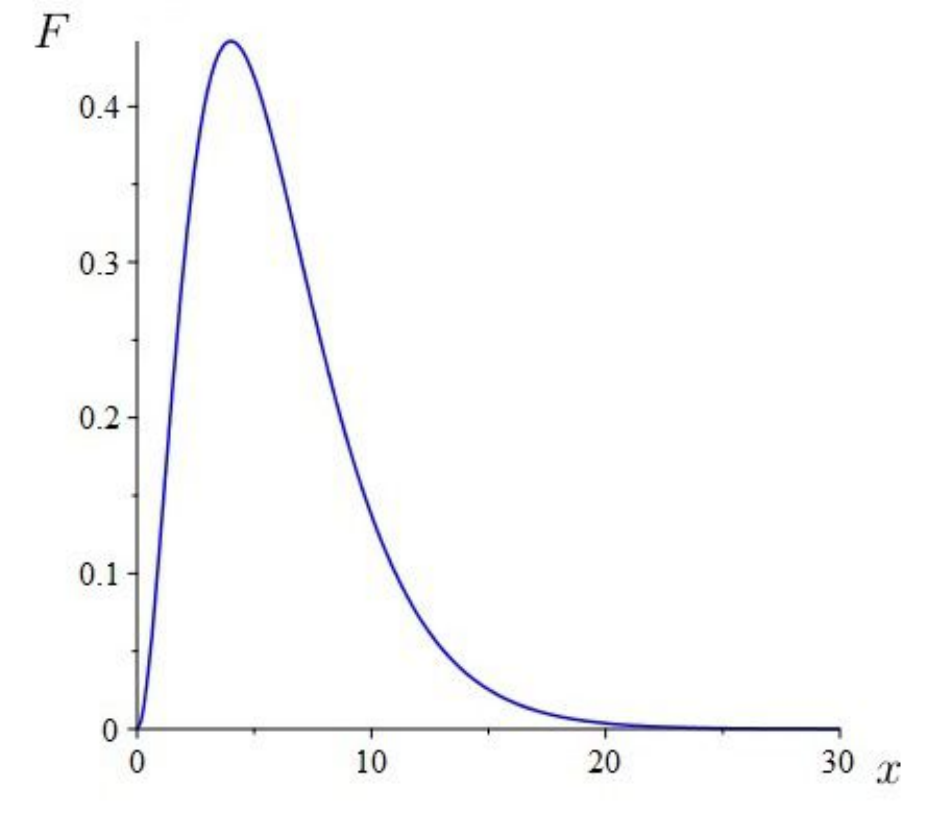}
\caption{Exact solution to the radial equation for $\tilde{n} = 2$, $l = 1$ in hydrogen}\label{fig:R21Solution}
\end{center}
\end{figure}

\begin{figure}[ht]
\begin{center}
\includegraphics[width=4in]{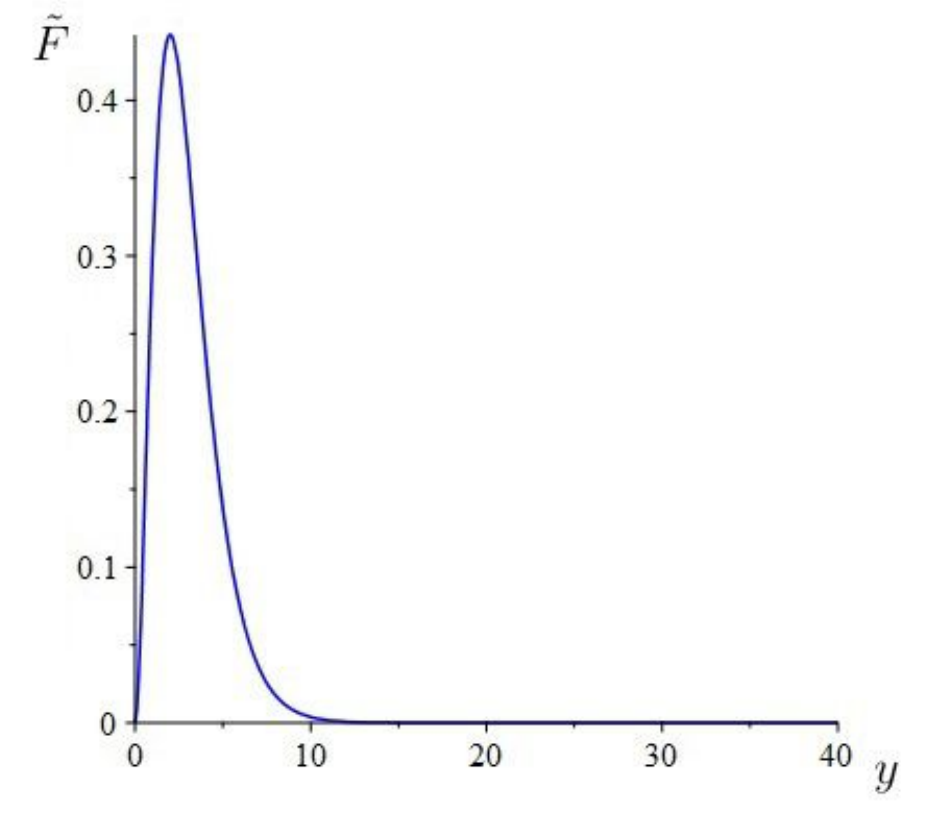}
\caption{Plot of $\tilde{F}(y) = yG(y)$}\label{fig:FTildeSolution}
\end{center}
\end{figure}

For the case $\tilde{n} = 2$, $l = 1$, the radial equation is 
\begin{equation}
\frac{1}{2} \frac{d^2F}{dx^2} + \bigg[\frac{1}{x} -\frac{1}{8} -\frac{1}{x^2}\bigg]F = 0
\end{equation}
Using equation (B.4) in Appendix~\ref{second-appendix}, the exact solution to (2.17) is then as shown in Figure~\ref{fig:R21Solution}, given by setting $a=1$ in $rR_{21}$ to get
\begin{equation}
F(x) = \frac{1}{2\sqrt{6}}x^2e^{-x/2} = \frac{2}{\sqrt{6}}\bigg(\frac{x}{2}\bigg)^2e^{-x/2}
\end{equation}
Our attempts to directly implement (2.17) again brought to light an interesting problem with Schr\"{o}dinger's equation, explained in section 3.2.1, that only the trivial solution $F(x) = 0$ can be found by our collocation procedure when, as is the case here with (2.18), the solution is such that both the function value and the first derivative are zero at the left and right boundaries. As our focus here is on numerically exploring the relative performance of different patterns of collocation points in this quantum system while treating everything else as known, \emph{not} on looking for unknown solutions, we overcame this problem to enable us to continue with our numerical experiments by applying simple transformations to (2.17) and (2.18) as follows. First, we make the change of variable $y = \frac{x}{2}$ in (2.18) to get
\begin{equation}
\tilde{F}(y) \equiv F(2y) =  \frac{2}{\sqrt{6}}y^2e^{-y}
\end{equation}
Putting $x = 2y$ into (2.17) we find that the differential equation satisfied by $\tilde{F}(y)$ is
\begin{equation}
\frac{1}{8} \frac{d^2\tilde{F}}{dy^2} + \bigg[\frac{1}{2y} -\frac{1}{8} -\frac{1}{4y^2}\bigg]\tilde{F} = 0
\end{equation}
Next, we define
\begin{equation}
G(y) \equiv \frac{\tilde{F}(y)}{y} = \frac{2}{\sqrt{6}}ye^{-y}
\end{equation}
Then putting $\tilde{F}(y) = yG(y)$ into (2.20), we find that the differential equation satisfied by $G(y)$ is
\begin{equation}
\frac{y}{8} \frac{d^2G}{dy^2} + \frac{1}{4}\frac{dG}{dy} + \bigg[\frac{1}{2} -\frac{y}{8} -\frac{1}{4y}\bigg]G = 0
\end{equation}
The exact solution to (2.22) is (2.21) and we find that 
\begin{equation}
G^{\prime}(0) = \frac{2}{\sqrt{6}} 
\end{equation} 
Therefore our numerical experiments with different patterns of collocation points in this quantum system will proceed by first implementing (2.22), with boundary conditions $G^{\prime}(0) = \frac{2}{\sqrt{6}}$ and $G(\infty) = 0$. The desired approximation of (2.18) can then be obtained simply by multiplying the output by $y$ and using $F(x) = \tilde{F}\big(\frac{x}{2}\big)$. A plot of $\tilde{F}(y) = yG(y)$ is shown in Figure~\ref{fig:FTildeSolution}.

\begin{figure}[ht]
\begin{center}
\includegraphics[width=4in]{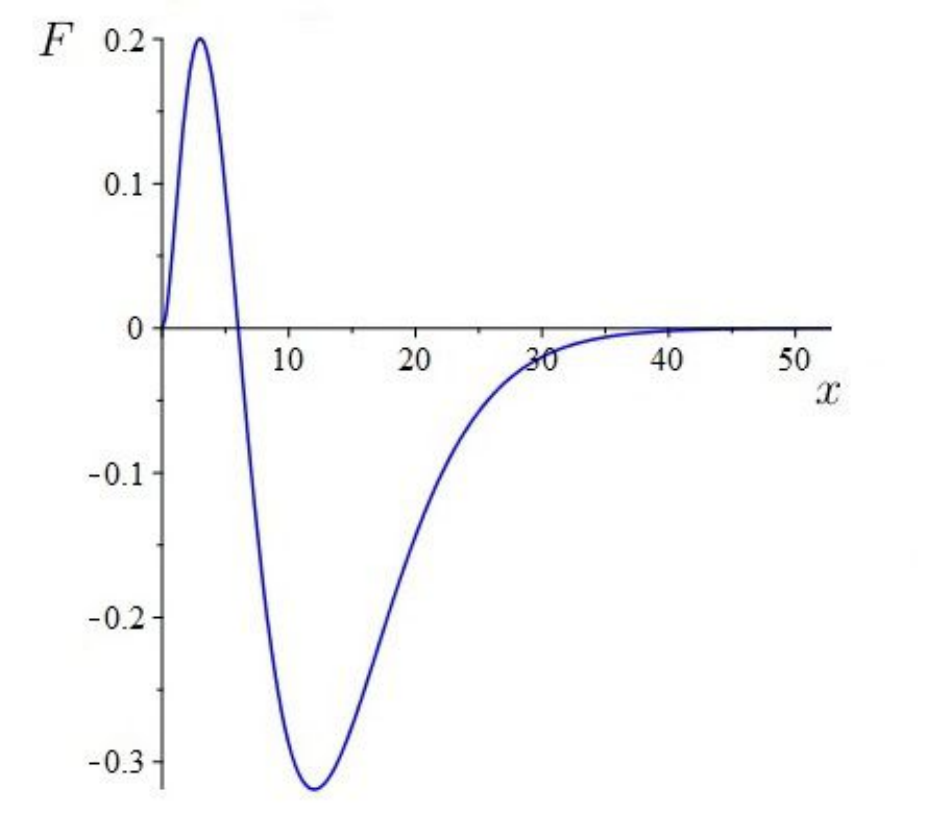}
\caption{Exact solution to the radial equation for $\tilde{n} = 3$, $l = 1$ in hydrogen}\label{fig:R31Solution}
\end{center}
\end{figure}

Similar issues arise in the cases  $\tilde{n} = 3$, $l = 1$ and  $\tilde{n} = 3$, $l = 2$, and they can be overcome in a very similar way. For the case $\tilde{n} = 3$, $l = 1$, the radial equation is 
\begin{equation}
\frac{1}{2} \frac{d^2F}{dx^2} + \bigg[\frac{1}{x} -\frac{1}{18} -\frac{1}{x^2}\bigg]F = 0
\end{equation}
Using equation (B.5) in Appendix~\ref{second-appendix}, the exact solution to (2.24) is then given by setting $a=1$ in $rR_{31}$ to get
\begin{equation}
F(x) = \frac{4}{81\sqrt{6}}(6 - x)x^2e^{-x/3} =  \frac{4}{3\sqrt{6}}\bigg(\frac{x}{3}\bigg)^2\bigg(2 - \frac{x}{3}\bigg)e^{-x/3}
\end{equation}
Figure~\ref{fig:R31Solution} shows a plot of (2.25). Here we can employ the same kind of transformation as before, beginning with the change of variable $y = \frac{x}{3}$ to get $\tilde{F}(y) \equiv F(3y)$ and then using $G(y) = \frac{\tilde{F}}{y^2}$ = $\frac{4}{3\sqrt{6}}(2 - y) e^{-y}$. Following the same procedure as before, we find that the differential equation satisfied by $G$ in this case is
\begin{equation*}
\frac{y}{18} \frac{d^2G}{dy^2} + \frac{2}{9}\frac{dG}{dy} + \bigg[\frac{1}{3} -\frac{y}{18}\bigg]G = 0
\end{equation*} 
so to carry out our experiments with different patterns of collocation points in this quantum system our strategy will be to implement this alternative differential equation first, with boundary conditions $G^{\prime}(0) = -\frac{4}{\sqrt{6}}$ and $G(\infty) = 0$, and then obtain the desired approximation of (2.25) simply by multiplying the output by $y^2$ and using $F(x) = \tilde{F}\big(\frac{x}{3}\big)$.

\begin{figure}[ht]
\begin{center}
\includegraphics[width=4in]{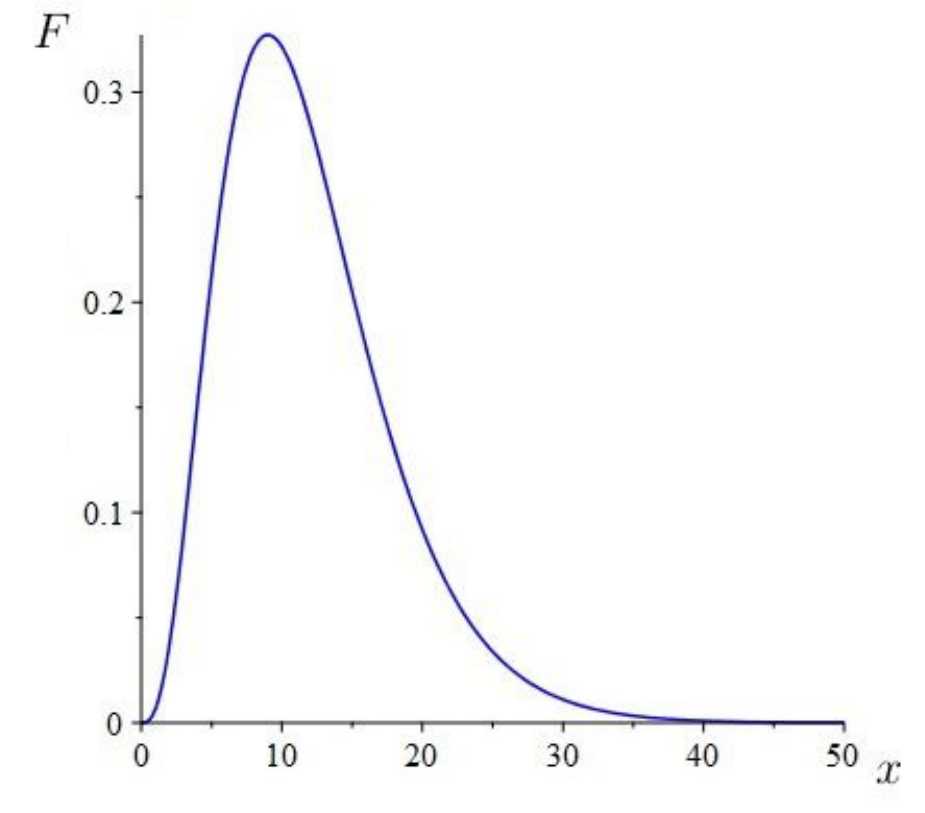}
\caption{Exact solution to the radial equation for $\tilde{n} = 3$, $l = 2$ in hydrogen}\label{fig:R32Solution}
\end{center}
\end{figure}

Finally, for the case $\tilde{n} = 3$, $l = 2$, the radial equation is 
\begin{equation}
\frac{1}{2} \frac{d^2F}{dx^2} + \bigg[\frac{1}{x} -\frac{1}{18} -\frac{3}{x^2}\bigg]F = 0
\end{equation}
Using equation (B.6) in Appendix~\ref{second-appendix}, the exact solution to (2.26) is then given by setting $a=1$ in $rR_{32}$ to get
\begin{equation}
F(x) = \frac{4}{81\sqrt{30}}x^3e^{-x/3} =  \frac{4}{3\sqrt{30}}\bigg(\frac{x}{3}\bigg)^3e^{-x/3}
\end{equation}
Figure~\ref{fig:R32Solution} shows a plot of (2.27). 

We can again employ the same kind of transformation as before, beginning with the change of variable $y = \frac{x}{3}$ and then using $G(y) = \frac{\tilde{F}}{y^2}$ = $\frac{4}{3\sqrt{30}} y e^{-y}$. In this case, we find that the differential equation satisfied by $G$ is 
\begin{equation*}
\frac{y}{18} \frac{d^2G}{dy^2} + \frac{2}{9}\frac{dG}{dy} + \bigg[\frac{1}{3} -\frac{y}{18} - \frac{2}{9y}\bigg]G = 0
\end{equation*} 
We can implement this with boundary conditions $G^{\prime}(0) = \frac{4}{3\sqrt{30}}$ and $G(\infty) = 0$, obtaining the desired approximation of (2.27) by multiplying the output by $y^2$ and using $F(x) = \tilde{F}\big(\frac{x}{3}\big)$. 

\newpage

\section{A nonlinear extension of Shore's framework}
To explore the performance of B-spline collocation at Gaussian points in a nonlinear Schr\"{o}dinger equation setting, we would like to implement a nonlinear version of Shore's basic equation
\begin{equation}
\frac{1}{2} \frac{d^2\psi}{dx^2} + [E - V]\psi = 0
\end{equation}   
incorporating a perturbation parameter and a nonlinear term analogous to the setup in the nonlinear perturbation problem discussed in Chapter XV of de Boor's book \emph{A Practical Guide to Splines}. That is to say, we would like to extend Shore's basic framework to a nonlinear equation of the form 
\begin{equation}
\frac{1}{2} \epsilon^2 \frac{d^2\psi}{dx^2} + [E - V]\psi + [\psi]^n = 0
\end{equation}   
where $\epsilon$ is a perturbation parameter (typically we want to explore solutions to this equation as $\epsilon \rightarrow 0$), and $n$ is an integer with $n > 1$. An equation exactly of the type (2.29) arises in the nonlinear Schr\"{o}dinger equation literature in relation to a Schr\"{o}dinger equation with cubic nonlinearity and a bounded potential of the form
\begin{equation}
i\epsilon\frac{\partial \Psi}{\partial t} = -\frac{\epsilon^2}{2m}\frac{\partial^2\Psi}{\partial x^2} + V\Psi - \gamma|\Psi|^2 \Psi
\end{equation}  
(which is shown in \cite{floer} to have standing wave solutions if $\gamma > 0$, $V$ is bounded, and $\epsilon$ is sufficiently small). To see this, by analogy with the usual linear Schr\"{o}dinger equation, we use separation of variables to seek solutions to (2.30) of the form 
\begin{equation}
\Psi(x, t) = \psi(x) e^{-iEt/\epsilon}
\end{equation} 
Putting (2.31) into (2.30), rearranging, and setting $m = 1$ and $\gamma = 1$ we get
\begin{equation}
\frac{1}{2} \epsilon^2 \frac{d^2\psi}{dx^2} + [E - V]\psi + [\psi]^3 = 0
\end{equation}   
which is exactly of the form (2.29) with $n = 3$. 

To implement this equation in our study using de Boor's methodology, we need to linearize it and also to find an exact solution for it in order to assess our approximations. We can linearize (2.32) by writing it as 
\begin{equation*}
\frac{1}{2} \epsilon^2 \frac{d^2\psi}{dx^2} = F(x, \psi(x), \psi^{\prime}(x)) \equiv -[E - V]\psi - [\psi]^3
\end{equation*}   
We then note that by Taylor's Theorem, expanding about the point $(v(x), v^{\prime}(x))$ (which in the iterative approximation process later we will treat as being derived from the result of the previous iteration), we have 
\begin{equation*}
F(x, \psi(x), \psi^{\prime}(x)) \approx F(x, v(x), v^{\prime}(x)) + (\psi(x) - v(x))\frac{\partial F}{\partial v(x)} +  (\psi^{\prime}(x) - v^{\prime}(x))\frac{\partial F}{\partial v^{\prime}(x)}
\end{equation*}
But
\begin{equation*}
F(x, v(x), v^{\prime}(x)) = -[E - V] v(x) - [v(x)]^3
\end{equation*}
\begin{equation*}
\frac{\partial F}{\partial v(x)} = -[E - V] - 3[v(x)]^2
\end{equation*}
\begin{equation*}
\frac{\partial F}{\partial v^{\prime}(x)} = 0 
\end{equation*}
Therefore 
\begin{equation*}
F(x, \psi(x), \psi^{\prime}(x)) \approx (-[E - V] - 3[v(x)]^2)\psi(x) + 2[v(x)]^3
\end{equation*}
so we can write the linearized form of the differential equation as 
\begin{equation}
\frac{1}{2} \epsilon^2 \frac{d^2\psi}{dx^2} + (3[v(x)]^2 + [E - V])\psi(x) = 2[v(x)]^3
\end{equation}
Given suitable choices of $E$ and $V$ and boundary conditions on $\psi(x)$, this can now be implemented using de Boor's methodology. 

To find an exact solution for (2.32) we need to specify $[E - V]$. For the purposes of our study, in which the focus is on exploring the numerical performance of B-spline collocation at Gaussian points rather than on physical applications of (2.32), we will assume an invariant potential (i.e., a quasi-free space) and set $[E - V] = -\frac{1}{2}$. This gives an equation of the form 
\begin{equation}
\frac{1}{2} \epsilon^2 \frac{d^2\psi}{dx^2} - \frac{1}{2}\psi + [\psi]^3 = 0
\end{equation} 
with linearized form 
\begin{equation}
\frac{1}{2} \epsilon^2 \frac{d^2\psi}{dx^2} + \bigg(3[v(x)]^2 - \frac{1}{2}\bigg)\psi(x) = 2[v(x)]^3
\end{equation}
and simple trial and error with functions of the form $\text{sech}(x)$ (mentioned in \cite{floer}, equation (1.3), p. 399) shows that an exact solution for (2.34) is
\begin{equation}
\psi(x) = \frac{1}{\cosh\big(\frac{x}{\epsilon}\big)}
\end{equation} 
We will therefore implement (2.34) using the linearized form (2.35), comparing our approximations for different values of $\epsilon$ with exact solutions of the form (2.36). 

\begin{figure}[ht]
\begin{center}
\includegraphics[width=4in]{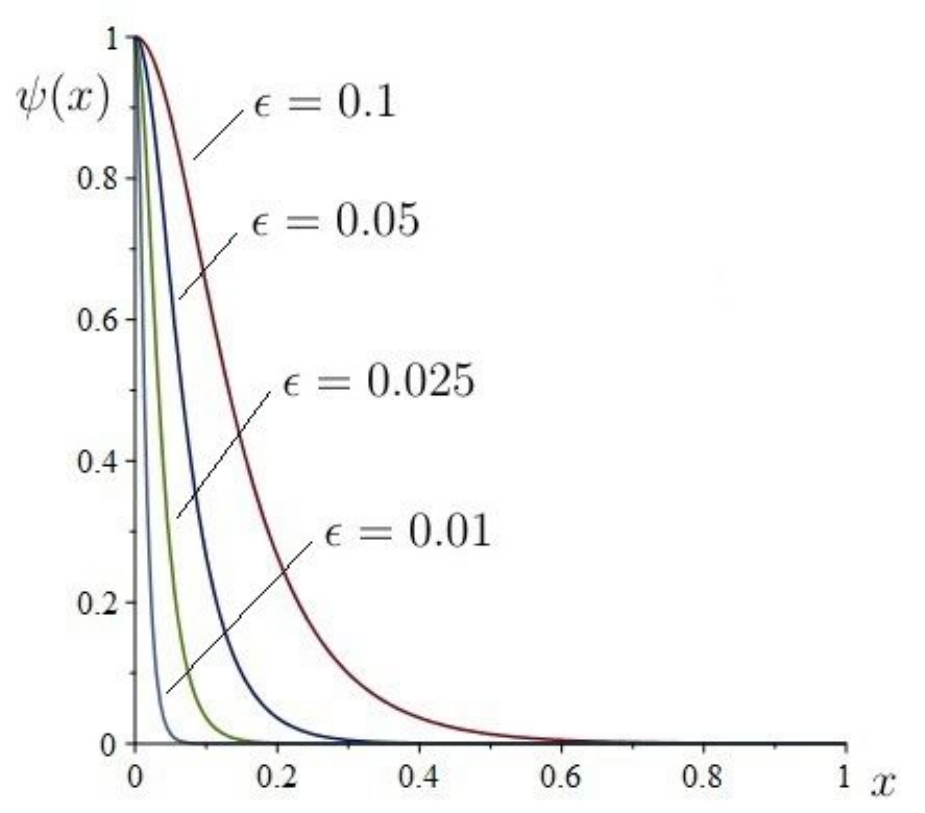}
\caption{Exact solutions for cubic Schr\"{o}dinger equation with $[E - V] = -\frac{1}{2}.$}\label{fig:nonlinear}
\end{center}
\end{figure}

Figure~\ref{fig:nonlinear} shows the exact solutions for different values of $\epsilon$ in the interval $[0, 1]$. We will apply the boundary conditions $\psi(0) = 1$ and $\psi(1) = 0$ within this interval.

This problem exhibits the classic features of a singular perturbation problem (also known as a \emph{boundary layer problem}) in which one explores how the solutions of a boundary value problem change as a parameter like $\epsilon$ here approaches zero. In the case of equation (2.34), it can be seen by inspection that as $\epsilon \rightarrow 0$, the differential equation becomes more and more like an algebraic equation which does not satisfy the boundary condition $\psi(0) = 1$. Therefore, as Figure~\ref{fig:nonlinear} shows, for smaller $\epsilon$-values the exact solution exhibits a sharper `bend' as $x$ approaches the origin from the right, and this can cause problems for approximation. We will want to explore how B-spline collocation at Gaussian points is able to deal with this difficulty. For a book-length treatment of singular perturbation problems, see \cite{omalley}.     

\section{Eigenfunction approach versus eigenvalue approach}

Although we are using the same fundamental radial equation as Shore for our numerical experiments (equation (I.1) in \cite{shore}), our approach is different from Shore's in a way that will now be made clear. Schr\"{o}dinger's radial equation for the hydrogen atom, with the boundary conditions implemented by Shore, is actually an example of a regular Sturm-Liouville problem of the general form
\begin{equation}
\frac{\mathrm{d}}{\mathrm{d} x} \bigg(p(x) \frac{\mathrm{d} \phi}{\mathrm{d} x} \bigg) + \big(q(x) + \lambda w(x)\big) \phi = 0
\end{equation}
\begin{equation*}
A_1 \phi(a) + A_2 \phi^{\prime}(a) = 0
\end{equation*}
\begin{equation*}
B_1 \phi(b) + B_2 \phi^{\prime}(b) = 0
\end{equation*}
for $x \in [a, b]$, where the aim is to find the eigenvalues $\lambda$ and corresponding eigenfunctions $\phi$. For example, one of Shore's implementations of the radial equation for the electron in the hydrogen atom is of the form (2.37) with $a = 0$, $b = 10$, $p(x) = \frac{1}{2}$, $q(x) = \frac{1}{x}$, $w(x) = 1$, $A_1 = B_1 = 1$, and $A_2 = B_2 = 0$. Using cubic spline collocation, Shore implements (2.37) as a matrix generalized eigenvalue problem
\begin{equation}
\bigg[\frac{\mathrm{d}}{\mathrm{d} x} \bigg(p(x) \frac{\mathrm{d} }{\mathrm{d} x} \bigg) + q(x)\bigg]\phi = - \lambda w(x) \phi
\end{equation} 
(cf. equation (VII.6) in \cite{shore}). The eigenvalues for the matrix system (2.38) are easily found numerically using standard methods for generalized eigenvalue problems. Shore's emphasis is very much on finding point estimates for the eigenvalues $\lambda$ in this way, which correspond to the quantum energy levels of the electron in Schr\"{o}dinger's theory of the hydrogen atom. The eigenvectors in the matrix system version of (2.38) could, in principle, then be used to obtain approximations of the eigenfunctions $\phi$ as a by-product, but Shore is not concerned very much with this. 

In contrast to Shore's approach, we seek to study the relative performance of different patterns of collocation points in numerically approximating the  eigenfunctions $\phi$ in (2.37), i.e., the wave functions of Schr\"{o}dinger's equation, in conjunction with different box sizes, mesh sizes and orders of polynomial approximants. We want to do this under `laboratory conditions' in which everything else that can influence approximation accuracy is fully known and controlled for. To this end, we take $\lambda$ as known in (2.37), and thereby convert the Sturm-Liouville problem above into a two-point boundary value problem with only $\phi$ as the unknown, perfectly suited for the machinery in Chapter XV of \cite{deboor2}. For example, to implement the radial equation for the ground state of the electron in the hydrogen atom in (2.10), we convert (2.37) into a two-point BVP by setting  $a = 0$, $b = 10$, $p(x) = \frac{1}{2}$, $q(x) = \frac{1}{x}$, $\lambda = -\frac{1}{2}$, $w(x) = 1$, and by replacing the boundary conditions in (2.37) by $\phi^{\prime}(0) = 2$, $\phi(10) = 0$. We then implement this system using de Boor's B-spline collocation methodology (described in detail in the next chapter), focusing purely on numerically approximating $\phi$.      

Rather than giving us just point estimates of single numbers, each accompanied by a single indicator of approximation error, our approach yields both visually rich and numerically rich approximation outputs consisting of entire wave functions that can be visually compared with known exact solutions, as well as detailed sets of approximation errors for the wave functions at various locations in the breakpoint sequences used in the collocation process. This can provide more detailed insights into the relative performance of different patterns of collocation points. Our approach is also more suitable for extending Shore's framework to nonlinear Schr\"{o}dinger equations, as described in section 2.2. It is not clear how nonlinear Schr\"{o}dinger equations could be studied using Shore's methodology.     

\chapter{B-splines and collocation at Gaussian points}

This chapter provides some necessary background on B-splines and concepts relating to collocation at Gaussian points, as well as outlining the role of some of the relevant MATLAB and Fortran 77 routines originally provided by de Boor in his book \emph{A Practical Guide to Splines} \cite{deboor2}. These have been translated into Maple code for the purposes of this dissertation. We begin in section 3.1 by reviewing the key theory and practical issues relating to piecewise polynomial approximation using B-splines, highlighting the roles of the subroutines INTERV, PPVALU, BSPLVP, BVALUE, BSPLPP and SPLINT. In section 3.2 we then use key ideas from the paper by de Boor and Swartz \cite{deboor} and Chapter XV of \cite{deboor2} to set out our approach to implementing Shore's radial Schr\"{o}dinger equation using B-spline collocation at Gaussian points, focusing in particular on how the use of Gaussian points can reduce approximation errors in this specific context. The key subroutines here are COLPNT, DIFEQU, NEWNOT and COLLOC. 

\section{Piecewise polynomial approximation using B-splines}

A key component of our approach to collocation at Gaussian points based on \cite{deboor2} is the use of B-splines to produce piecewise polynomial approximations to the Schr\"{o}dinger wave functions in our study. Piecewise polynomial (pp) functions generally perform far better as approximants in practical situations than single polynomials (see, e.g., \cite{deboor2}, Chapter II, \cite{powell}, p. 212, \cite{rivlin}, p.104). Splines can be viewed as pp functions with pieces that `blend as smoothly as possible' due to continuity conditions on their derivatives (\cite{deboor2}, p. 105), but de Boor uses the term more inclusively to mean `all linear combinations of B-splines'. B-splines are a numerically convenient set of pp functions used as a basis for all others. 

\subsection{B-splines as a basis for pp function spaces}

Using the same notational conventions as de Boor's book \emph{A Practical Guide to Splines}, a pp function $f$ of order $k$ is defined for $i = 1, \ldots, l$ as

$f(x) = P_i(x)$ if $\xi_i < x < \xi_{i+1}$

where $\xi = \{\xi_1, \ldots, \xi_{l+1}\}$ is a strictly increasing sequence of breakpoints and $P = \{P_1, \ldots, P_l\}$ is any sequence of polynomials of order $k$ (i.e., of degree $< k$). At each breakpoint other than $\xi_1$ and $\xi_{l+1}$, the pp function is (arbitrarily) defined for computational purposes as taking the value from the right, i.e., $f(\xi_i) = f(\xi_i^{+})$ for $i = 2, \ldots, l$. The collection of all pp functions of order $k$ with breakpoint sequence $\xi$ is a linear space of dimension $kl$ denoted by $\Pi_{<k, \xi}$. 

For computational purposes, de Boor represents the pp function $f \in \Pi_{<k, \xi}$ using a structure he calls a ppform, consisting of the integers $k$ and $l$, the breakpoint sequence $\xi$, and the $k \times l$ matrix of the right-derivatives of $f$ at the breakpoints:
\begin{equation*}
C = \big[D^{j-1}f(\xi_i^+)\big]_{j=1; i = 1}^{k \ \ \ l}
\end{equation*}
\begin{equation*}
= 
\begin{bmatrix}
f(\xi_1^+) & f(\xi_2^+) & \cdots & f(\xi_l^+) \\
Df(\xi_1^+) & Df(\xi_2^+) & \cdots & Df(\xi_l^+) \\
\vdots & \vdots & \ddots & \vdots \\
D^{k-1}f(\xi_1^+) & D^{k-1}f(\xi_2^+) & \cdots &  D^{k-1}f(\xi_l^+) 
\end{bmatrix}
\end{equation*}
In our numerical experiments, the output from COLLOC is essentially the transpose of this matrix $C$ for the ppform of the B-spline approximation. It is necessary to process this output further because the required pp function coefficients, say for the $i$-th piece 
\begin{equation*}
P_i(x) = \sum_{j = 1}^k c(i, j)(x - \xi_i)^{j - 1} 
\end{equation*}
for $\xi_i \leq x < \xi_{i+1}$ are not of the form $D^{j-1}f(\xi_i^+)$ as in the matrix $C$, but rather of the form
\begin{equation*}
\frac{D^{j-1}f(\xi_i^+)}{(j - 1)!}
\end{equation*} 
(see \cite{deboor2}, pp. 71-73). We make this adjustment in the post-output processing part of our computer routines, examples of which are provided in Appendices E to G. 

The subroutine PPVALU computes the values of $f$ and its derivatives at a given site $x$ using as inputs the integers $k$ and $l$, a one-dimensional array containing the breakpoints $\xi$, and a two-dimensional array containing the matrix $C$. In our numerical experiments, this output is used within the subroutine DIFEQU to construct approximation errors for our collocation approximations. PPVALU uses the subroutine INTERV to place each site $x$ in the correct place within the breakpoint sequence $\xi$. 

In general, it is necessary to impose continuity conditions on pp functions and their derivatives, of the form
\begin{equation*}
\text{jump}_{\xi_i} D^{j-1} f = 0
\end{equation*}
for $j = 1, \ldots, \nu_i$ and $i = 2, \dots, l$, where the notation means `the jump of the function across the site $\xi_i$', and $\nu = \{\nu_2, \nu_3, \dots, \nu_l\}$ is a set of nonnegative integers with $\nu_i$ counting the number of continuity conditions required at $\xi_i$. (Note that there is no need for elements $\nu_1$ or $\nu_{l+1}$ in this list as continuity conditions are only needed to govern how different pieces of the pp function `meet' at interior breakpoints). For example, $\nu_i = 2$ means that both the function and the first derivative are required to be continuous at $\xi_i$, whereas $\nu_i = 0$ means that there are no continuity conditions at $\xi_i$. These continuity conditions are linear and homogeneous, so the subset of all $f \in \Pi_{<k, \xi}$ satisfying them is a linear subspace of $\Pi_{<k, \xi}$ denoted by $\Pi_{<k, \xi;\nu}$ (see \cite{deboor2}, p. 82). The dimension of $\Pi_{<k, \xi;\nu}$ is 
\begin{equation*}
n = kl - \sum_{i=2}^l \nu_i
\end{equation*}  
B-splines emerge from the desire to have a numerically convenient basis for $\Pi_{<k, \xi;\nu}$. One basis for this space which is not numerically convenient is the `truncated power basis' (see \cite{deboor2}, pp. 82-84) which consists of the double-sequence

$\varphi_{ij}$, $j = \nu_i, \ldots, k-1$ and $i = 2. \ldots, l$

where 
\begin{equation*}
\varphi_{ij} = \frac{(x -\xi_i)_+^j}{j!}
\end{equation*}
and where $(x - \xi_i)_+^j \equiv (\text{max}\{(x - \xi_i), 0\})^j$ is a truncated power function. This is a basis for $\Pi_{<k, \xi;\nu}$ in the sense that every pp function $f \in \Pi_{<k, \xi;\nu}$ can be written in a unique way in the form
\begin{equation*}
f = \sum_{i=1}^l \sum_{j=\nu_i}^{k-1} \alpha_{ij}\varphi_{ij}
\end{equation*}
This basis is not well-suited for numerical work for a number of reasons, particularly because truncated power functions can grow rapidly irrespective of the behaviour of $f$, and also some of the basis functions $\varphi_{ij}$ can become nearly collinear, leading to numerical difficulties (see, e.g., the example in \cite{deboor2}, p. 85). These difficulties can be overcome by using as the basis elements certain \emph{divided differences} of the truncated power functions instead, which have the property that they each have support only over a small interval, vanishing elsewhere. B-splines are basis elements for $\Pi_{<k, \xi;\nu}$ defined in this way.  

To formally introduce B-splines, let $t = \{t_j\}$ be a nondecreasing sequence of numbers (these are called `knots' in the context of splines, and can be viewed as an extension of the breakpoint sequence $\xi$ defined earlier in the sense that $t$ can incorporate the elements of a given $\xi$ but does not have to be strictly increasing, and in principle it can be finite or infinite as required). Then the $j$-th normalised B-spline of order $k$ (i.e., of degree $k-1$) using knot sequence $t$ is denoted by $B_{j,k,t}$, and its value at a site $x \in \mathbf{R}$ is given by
\begin{equation}
B_{j,k,t}(x) = (t_{j+k} - t_j)[t_j, \ldots, t_{j+k}](\cdot - x)_+^{k-1}
\end{equation}
where the notation $[t_j, \ldots, t_{j+k}]g$ denotes the $k$-th divided difference of a function $g$ at the sites $t_j, \ldots, t_{j+k}$ (divided differences are discussed in \cite{deboor2}, Chapter I, and \cite{powell}, Chapter 5), and the dot placeholder notation means that $x$ is regarded as being fixed when calculating the divided difference of the truncated power function, so the latter is being treated as a function of a single variable. This formal definition can be used to generate B-splines of any required order, but it is more convenient to use a recurrence relation (proved in \cite{deboor2}, p. 90) which says that for $k > 1$,
\begin{equation}
B_{j,k,t}(x) = \frac{(x-t_j)B_{j,k-1,t}(x)}{t_{j+k-1} - t_j} + \frac{(t_{j+k} - x)B_{j+1,k-1,t}(x)}{t_{j+k} - t_{j+1}}
\end{equation} 
This relation can be used to generate B-splines by induction, starting from $B_{j,1,t}(x)$, which in turn can be obtained from the formal definition (3.1) above as
\begin{equation*}
B_{j,1,t}(x) = (t_{j+1} - t_j)[t_j, t_{j+1}](\cdot - x)_+^0
\end{equation*}
\begin{equation*}
= (t_{j+1} - t_j) \frac{\{(t_{j+1} - x)_+^0 - (t_j - x)_+^0\}}{(t_{j+1} - t_j)}
\end{equation*}
\begin{equation*}
= (t_{j+1} - x)_+^0 - (t_j - x)_+^0
\end{equation*}
\begin{equation*}
= \left \{
\begin{array}{c c}
1 & \text{if }t_j \leq x < t_{j+1}\\
0 & \text{otherwise } 
\end{array} \right.
\end{equation*}
Note that the B-spline $B_{j,1,t}(x)$ is a piecewise polynomial of order $1$ and has support $[t_j, t_{j+1})$, so it is continuous from the right in accordance with the convention for pp functions stated earlier. By putting $B_{j,1,t}(x)$ into the recurrence relation (3.2), we obtain the B-spline $B_{j,2,t}(x)$ which is a piecewise polynomial of order $2$ with support $[t_j, t_{j+2})$. B-splines of higher order can be found via the recurrence relation (3.2) above in a convenient way using a tableau similar to the one commonly used to work out divided differences of functions. This is discussed in \cite{deboor2}, p. 110, and \cite{powell}, p. 235.  

The Curry-Schoenberg Theorem (proved in \cite{deboor2}, pp. 97-98) shows that the B-splines as defined above constitute a basis for $\Pi_{<k, \xi;\nu}$ under certain conditions. Specifically, the theorem says that the sequence $\{B_{1,k,t}, B_{2,k,t}, \ldots, B_{n,k,t}\}$ is a basis for  $\Pi_{<k, \xi;\nu}$ if:

(i) $\xi = \{\xi_1, \ldots, \xi_{l+1}\}$ is a strictly increasing sequence of breakpoints;

(ii) $\nu = \{\nu_2, \nu_3, \dots, \nu_l\}$ is a set of nonnegative integers with $\nu_i \leq k$ for all $i$;

(iii) $t = \{t_1, \ldots,t_{n+k}\}$ is a nondecreasing sequence with $n = kl - \sum_{i=2}^l \nu_i = \text{dim}\Pi_{<k, \xi;\nu}$;  

(iv) for $i = 2, \ldots, l$, the number $\xi_i$ occurs exactly $k - \nu_i$ times in t;

(v) $t_1 \leq t_2 \leq \ldots t_k \leq \xi_1$ and $\xi_{l+1} \leq t_{n+1} \leq \ldots \leq t_{n+k}$.

These specifications provide the necessary information for generating a knot sequence $t$ from a given breakpoint sequence $\xi$ with the desired amount of `smoothness' (i.e., number of continuity conditions), and we can then construct a B-spline basis using the recurrence relation (3.2) above. The number of continuity conditions at a breakpoint $\xi_i$ is determined by the number of times $\xi_i$ appears in $t$, in the sense that each repetition of $\xi_i$ reduces the number of continuity conditions at that breakpoint by one. If $\xi_i$ appears $k$ times in $t$, this corresponds to imposing no continuity conditions at $\xi_i$. If $\xi_i$ appears $k-1$ times, the function is continuous at $\xi_i$, but not its first or higher derivatives. If $\xi_i$ appears $k-2$ times, the function and its first derivative are continuous at $\xi_i$, but not its second and higher derivatives; and so on. Note that a convenient choice of knot sequence is to make the first $k$ knot points equal to $\xi_1$, and the last $k$ knot points equal to $\xi_{l+1}$, thus imposing no continuity conditions at $\xi_1$ and $\xi_{l+1}$.

\begin{figure}
\begin{center}
\includegraphics[width=6in]{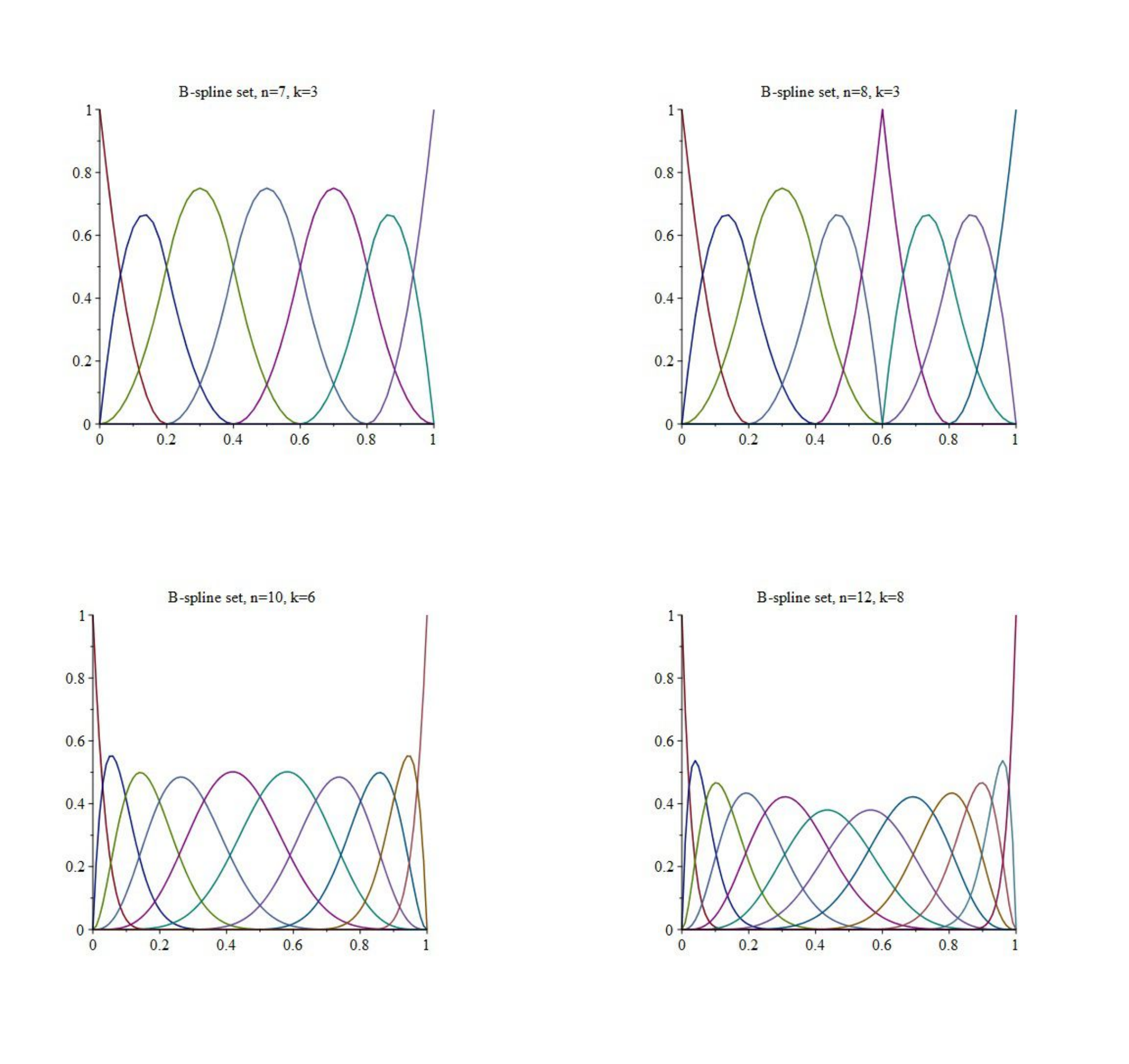}
\caption{B-spline sets for various values of $n$ and $k$.}\label{fig:bsplinesets}
\end{center}
\end{figure}

To illustrate these ideas, we use Maple programs based on the procedure described on page 113 of \cite{deboor2} (an example is provided in Appendix C) which call the subroutines INTERV and BSPLVP to produce B-spline sets with various specifications. These are plotted in Figure~\ref{fig:bsplinesets}. 

The top left plot in Figure~\ref{fig:bsplinesets} shows the quadratic B-spline set of order $3$ with the breakpoint sequence $\xi = \{0, 0.2, 0.4, 0.6, 0.8, 1.0\}$ and corresponding knot sequence $t = \{0, 0, 0, 0.2, 0.4, 0.6, 0.8, 1.0, 1.0, 1.0\}$. We have $k = 3$, $l = 5$, and $\nu = \{\nu_2, \nu_3, \nu_4, \nu_5\} = \{2, 2, 2, 2\}$, so the dimension is $n = 3\times 5 - (2+2+2+2) = 7$. Therefore we expect seven B-splines in this set, which is indeed what the top left plot in Figure~\ref{fig:bsplinesets} shows.

To allow the first derivative at breakpoint $0.6$ to become discontinuous, we repeat this breakpoint once in the knot sequence, so the knot sequence becomes

$t = \{0, 0, 0, 0.2, 0.4, 0.6, 0.6, 0.8, 1.0, 1.0, 1.0\}$

We still have $k = 3$ and $l = 5$, but now $\nu = \{\nu_2, \nu_3, \nu_4, \nu_5\} = \{2, 2, 1, 2\}$, so the dimension is now $n = 3\times 5 - (2+2+1+2) = 8$. Therefore we expect eight B-splines in this set. These are shown in the top right plot in Figure~\ref{fig:bsplinesets}, which also displays the effect of the discontinuous first derivative at $0.6$.   

The lower left part of Figure~\ref{fig:bsplinesets} shows a B-spline set of order $6$, i.e., quintic B-splines. In this case, $k=6$, $l=5$ and $\nu = \{5, 5, 5, 5\}$, so the dimension is $n=10$. The knot sequence $t$ has six repetitions of the breakpoints $0$ and $1.0$. Finally, the lower right part of Figure~\ref{fig:bsplinesets} shows a B-spline set of order $8$, i.e., heptic B-splines. Here, $k=8$, $l=5$ and $\nu = \{7, 7, 7, 7\}$, so the dimension is $n=12$. The knot sequence $t$ has eight repetitions of $0$ and $1.0$ in this case.   

\subsection{B-spline interpolation}

For computatonal purposes, de Boor (\cite{deboor2}, p. 100) uses the Curry-Schoenberg Theorem to represent the pp function $f \in \Pi_{<k, \xi;\nu}$ as a structure he calls a B-form, consisting of the integers $k$ and $n$, the knot sequence $t$, and a set of coefficients $\alpha = \{\alpha_1, \ldots, \alpha_n\}$ of $f$ with respect to the B-spline basis $\{B_{1,k,t}, B_{2,k,t}, \ldots, B_{n,k,t}\}$, such that the value of $f$ at a site $x \in [t_k, t_{n+1}]$ is given by
\begin{equation} 
f(x) = \sum_{i=1}^n \alpha_i B_{i,k,t}(x)
\end{equation}
The subroutine BVALUE computes the values of $f$ and its derivatives at a given site $x$ from its B-form (so it is the analogue of PPVALU for ppforms). In our numerical procedures using COLLOC, the approximate Schr\"{o}dinger wave functions will first be obtained as B-forms. For output purposes, these will then be converted to the ppform described earlier using the subroutine BSPLPP (\cite{deboor2}, pp. 117-120). 

The B-form described above can be used to interpolate a function $g$ at $n$ interpolation sites $\tau = (\tau_1, \ldots, \tau_n)$ by solving a matrix system based on (3.3):
\begin{equation} 
\begin{bmatrix}
B_{1,k,t}(\tau_1) & B_{2,k,t}(\tau_1) & \cdots & B_{n,k,t}(\tau_1) \\
B_{1,k,t}(\tau_2) & B_{2,k,t}(\tau_2) & \cdots & B_{n,k,t}(\tau_2)  \\
\vdots & \vdots & \ddots & \vdots \\
B_{1,k,t}(\tau_n) & B_{2,k,t}(\tau_n) & \cdots & B_{n,k,t}(\tau_n)  
\end{bmatrix}
\begin{bmatrix}
\alpha_1 \\
\alpha_2 \\
\vdots \\
\alpha_n
\end{bmatrix}
=
\begin{bmatrix}
g(\tau_1) \\
g(\tau_2) \\
\vdots \\
g(\tau_n)
\end{bmatrix}
\end{equation}
The knot sequence $t$ determines which B-splines of order $k$ will be involved in the spline approximation and the interpolation sites $\tau$ specify where the spline has to agree with the function $g$. The conditions under which this interpolation procedure will work are given in the Schoenberg-Whitney Theorem (proved in \cite{deboor2}, p. 173). In particular, we require the diagonal elements of coefficient matrix to be nonzero, i.e., $B_{i,k,t}(\tau_i) \neq 0$ for $i = 1, \ldots,n$, which means that each interpolation point $\tau_i$ must lie within the support $[t_i, t_{i+k})$ of the B-spline $B_{i,k,t}$. The knot sequence $t$ needs to be chosen to accommodate this requirement. 

Due to some basic properties of B-splines, the coefficient matrix has a convenient `banded' structure making (3.4) easy to solve by Gaussian elimination without pivoting. The subroutine SPLINT oversees this and provides the B-form coefficients of the approximation $f$ of $g$. BVALUE can then be used with this B-form to evaluate the spline approximation at various points, e.g., for plotting. To illustrate this, we use a Maple program which calls SPLINT and BVALUE (provided in Appendix D) to determine the cubic spline that interpolates the Gauss hypergeometric function $g(x) = {}_2F_1\big([1, 1], [1], x \text{e}^{-x}\big)$ on the interval $[-1, 1]$, with the seven equally spaced interpolation points $\tau = (-1, -2/3, -1/3, 0, 1/3, 2/3, 1)$ and with knot sequence $t = (-1, -1, -1, -1, -1/3, 0, 1/3, 1, 1, 1, 1)$. Note that $n = 7$ and $k = 4$, so the knot sequence $t = \{t_1, \ldots,t_{n+k}\}$ has length $n + k = 11$. 

\begin{figure}
\begin{center}
\includegraphics[width=6in]{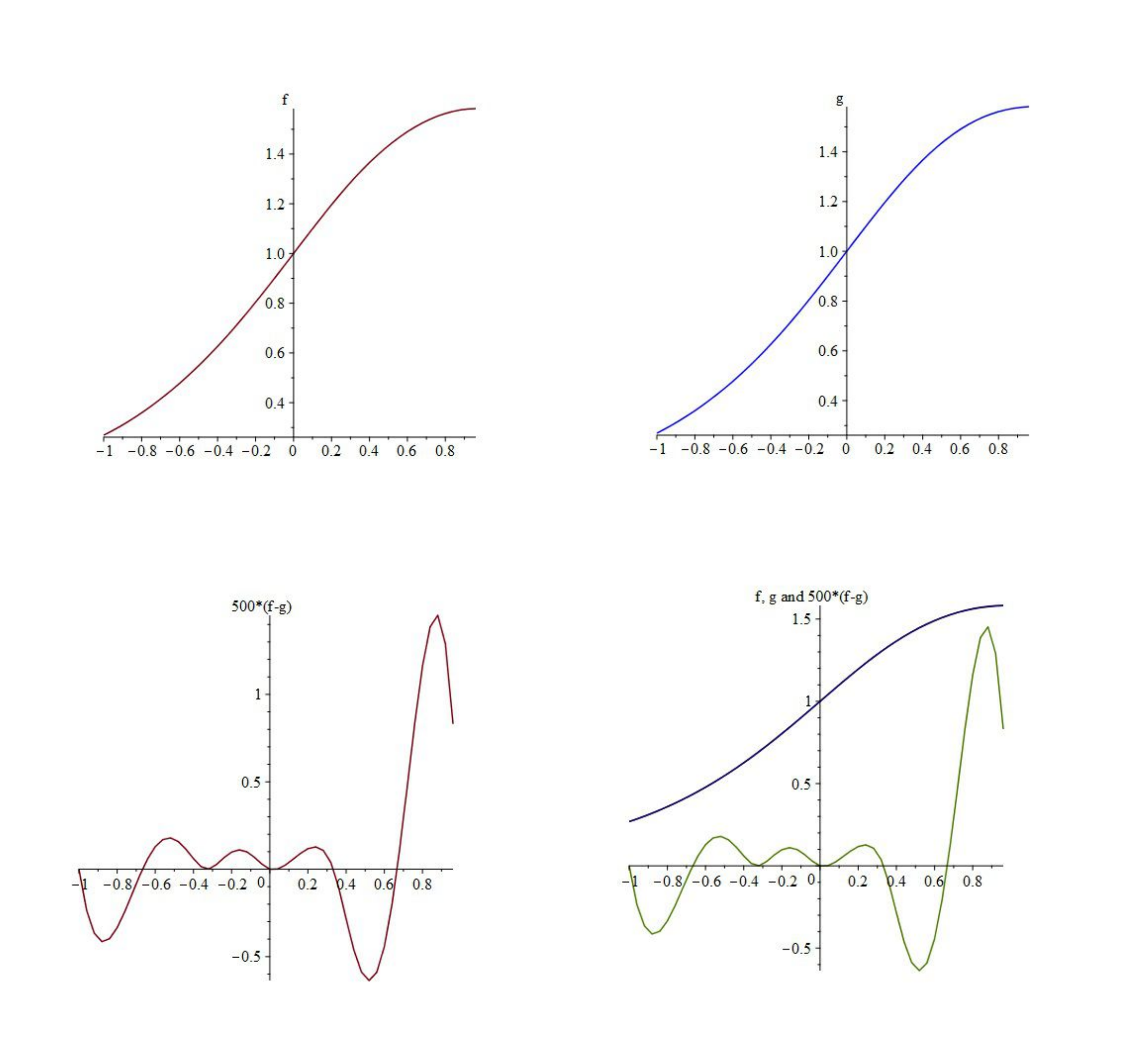}
\caption{Results for cubic spline interpolation of the Gauss hypergeometric function $g(x) = {}_2F_1\big([1, 1], [1], x \text{e}^{-x}\big)$.}\label{fig:hypergeometric}
\end{center}
\end{figure}   
\newpage
The upper part of Figure~\ref{fig:hypergeometric} shows that the cubic spline approximation $f$ is visually almost indistinguishable from the exact function $g$. However, there are some small approximation errors which are plotted in the magnified form $500 \times (f - g)$ in the lower left part of Figure~\ref{fig:hypergeometric}. The lower right part of Figure~\ref{fig:hypergeometric} shows the three plots superimposed.  

\section{Collocation at Gaussian points with Shore's equation}

In using the collocation procedure in \cite{deboor2} to approximate the solution of a second-order ordinary differential equation with boundary conditions on the interval $[a, b]$, the aim is to determine the B-form of a pp function $z \in \Pi_{<k, \xi;\nu}$ which exactly satisfies the differential equation and its boundary conditions at certain sites $\tau = (\tau_1, \ldots, \tau_n)$, where $\tau_1 = a$ and $\tau_n = b$. The form of the differential equation is specified in the subroutine DIFEQU and the collocation sites $(\tau_2, \ldots, \tau_{n-1})$ are constucted from specifications in COLPNT. The subroutine COLLOC oversees the iterative solution of the system using Newton's method, calling on NEWNOT, if required, to seek improvements by making nonlinear adjustments to the relative positions of breakpoints and collocation sites. Note that this collocation process is different from the interpolation procedure described in the previous section, where the pp function is required to match only the \emph{values} of another function $g$ at the interpolation sites.    

\subsection{B-spline collocation using de Boor's subroutines}

All the Schr\"{o}dinger equations in our study are supplied to DIFEQU in the form 
\begin{equation}
v_1(x) z(x) + v_2(x) Dz(x) + v_3(x) D^2 z(x) = v_4(x)
\end{equation}
by varying the specifications of $v_1(x)$, $v_2(x)$, $v_3(x)$ and $v_4(x)$. For example, Shore's radial equation for the ground state of the electron in the hydrogen atom in (2.10) requires the specifications
\begin{equation*}
v_1(x) = \left \{
\begin{array}{c c}
0 & \text{for } x = \tau_1\\
\frac{1}{x} - \frac{1}{2} & \ \text{for } x \in (\tau_2, \ldots, \tau_{n-1}) \\
1  & \text{for } x = \tau_n
\end{array} \right.
\end{equation*}
\begin{equation*}
v_2(x) = \left \{
\begin{array}{c c}
1 & \text{for } x = \tau_1\\
0  & \text{for } x \in (\tau_2, \ldots, \tau_n)
\end{array} \right.
\end{equation*} 
\begin{equation*}
v_3(x) = \left \{
\begin{array}{c c}
0 & \text{for } x \in (\tau_1, \tau_n)\\
\frac{1}{2} & \ \text{for } x \in (\tau_2, \ldots, \tau_{n-1})
\end{array} \right.
\end{equation*} 
\begin{equation*}
v_4(x) = \left \{
\begin{array}{c c}
2 & \text{for } x = \tau_1\\
0 & \ \text{for } x \in (\tau_2, \ldots, \tau_n)
\end{array} \right.
\end{equation*} 
(see DIFEQU in Appendix E), whereas the (linearized) Schr\"{o}dinger equation with cubic nonlinearity in (2.35) requires the specifications  
\begin{equation*}
v_1(x) = \left \{
\begin{array}{c c}
1 & \text{for } x = \tau_1\\
3[z_0(x)]^2 - \frac{1}{2} & \ \text{for } x \in (\tau_2, \ldots, \tau_{n-1}) \\
0  & \text{for } x = \tau_n
\end{array} \right.
\end{equation*}
\begin{equation*}
v_2(x) = 0
\end{equation*}
\begin{equation*}
v_3(x) = \left \{
\begin{array}{c c}
0 & \text{for } x \in (\tau_1, \tau_n)\\
\frac{1}{2}\epsilon^2 & \ \text{for } x \in (\tau_2, \ldots, \tau_{n-1})
\end{array} \right.
\end{equation*} 
\begin{equation*}
v_4(x) = \left \{
\begin{array}{c c}
1 & \text{for } x = \tau_1\\
2[z_0(x)]^3 & \ \text{for } x \in (\tau_2, \ldots, \tau_{n-1}) \\
0 & \text{for } x = \tau_n
\end{array} \right.
\end{equation*} 
where $z_0(x)$ here represents a prior estimate of the solution in the iterative procedure (see DIFEQU in Appendix G). 

Having specified the interval $[a, b]$ (referred to as the `box') and the breakpoints $\xi = (\xi_1, \ldots, \xi_{l+1})$, where $\xi_1 = a$ and $\xi_{l+1} = b$, the pp function approximant $z$ will then have $l$ polynomial pieces (referred to as the `mesh'). The box endpoints $a$ and $b$, and the mesh $l$, have to be supplied to COLLOC, along with the order $k$ of $z$. The program will then calculate $n = kl - 2(l-1)$ as the number of sites in $\tau$, with $\frac{n-2}{l}=k-2$ collocation sites per polynomial piece which it will look for in COLPNT. The knot sequence $t$ will be constructed to be of lengh $n+k$, giving degrees of freedom $\text{length}(t) - k = n$ to match the $n$ conditions represented by the $n-2$ collocation sites $(\tau_2, \ldots, \tau_{n-1})$ together with the boundary conditions at $\tau_1$ and $\tau_n$. Since also $n = \text{dim}\Pi_{<k, \xi;\nu}$, the calculation of the B-form of $z$ 
\begin{equation}
z(x) = \sum_{j=1}^n a_j B_{j,k,t}(x)
\end{equation}
will require the calculation of the $n$ B-splines $(B_{1,k,t}(x), B_{2,k,t}(x), \ldots, B_{n,k,t}(x))$ along with their first and second derivatives at each of the $n$ sites in $\tau$, with continuity conditions $\nu = (\nu_2, \nu_3, \ldots, \nu_l) = (2, 2, \ldots, 2)$ giving $\sum_{i=2}^l \nu_i = 2(l - 1)$. 

It will also be necessary to calculate the values of $v_1(x)$, $v_2(x)$, $v_3(x)$ and $v_4(x)$ in (3.5) at each of the $n$ sites in $\tau$. With these in hand, we can substitute (3.6) into (3.5) to get at each $\tau_j \in \tau$
\begin{equation*}
v_1(\tau_j) z(\tau_j) + v_2(\tau_j) Dz(\tau_j) + v_3(\tau_j) D^2 z(\tau_j) = v_4(\tau_j)
\end{equation*}
$\iff$
\begin{equation*}
v_1(\tau_j) (a_1B_{1,k,t}(\tau_j) + a_2B_{2,k,t}(\tau_j) + \cdots + a_nB_{n,k,t}(\tau_j))
\end{equation*} 
\begin{equation*}
+ v_2(\tau_j) (a_1DB_{1,k,t}(\tau_j) + a_2DB_{2,k,t}(\tau_j) + \cdots + a_nDB_{n,k,t}(\tau_j))
\end{equation*} 
\begin{equation*}
+ v_3(\tau_j) (a_1D^2B_{1,k,t}(\tau_j) + a_2D^2B_{2,k,t}(\tau_j) + \cdots + a_nD^2B_{n,k,t}(\tau_j)) = v_4(\tau_j)
\end{equation*} 
$\iff$
\begin{equation}
a_1(LB_{1,k,t})(\tau_j) + a_2(LB_{2,k,t})(\tau_j) + \cdots + a_n(LB_{n,k,t})(\tau_j) = v_4(\tau_j)
\end{equation} 
where $(LB_{j,k,t}) \equiv v_1B_{j,k,t} + v_2DB_{j,k,t} + v_3D^2B_{j,k,t}$. For all the $n$ sites in $\tau$, (3.7) then represents the matrix system
\begin{equation} 
\begin{bmatrix}
(LB_{1,k,t})(\tau_1) & (LB_{2,k,t})(\tau_1) & \cdots & (LB_{n,k,t})(\tau_1) \\
(LB_{1,k,t})(\tau_2) & (LB_{2,k,t})(\tau_2) & \cdots & (LB_{n,k,t})(\tau_2)  \\
\vdots & \vdots & \ddots & \vdots \\
(LB_{1,k,t})(\tau_n) & (LB_{2,k,t})(\tau_n) & \cdots & L(B_{n,k,t})(\tau_n)  
\end{bmatrix}
\begin{bmatrix}
a_1 \\
a_2 \\
\vdots \\
a_n
\end{bmatrix}
=
\begin{bmatrix}
v_4(\tau_1) \\
v_4(\tau_2) \\
\vdots \\
v_4(\tau_n)
\end{bmatrix}
\end{equation}
The system (3.8) can be solved in a single step for linear Schr\"{o}dinger equations such as (2.10), yielding the B-form coefficients $(a_1, \ldots, a_n)$ for (3.6), but iteration is needed for the (linearized) Schr\"{o}dinger equation with cubic nonlinearity in (2.35). An initial B-form $z_0(x)$ is used to specify $v_1(x)$, $v_2(x)$, $v_3(x)$ and $v_4(x)$ at each $\tau_j \in \tau$ and the system (3.8) is then solved to get an updated B-form for $z(x)$. The process is then repeated with the updated B-form and this continues until the B-forms converge, i.e., until $\text{max}\{|z_{r+1}(\tau_j)-z_r(\tau_j)|: \tau_j \in \tau \} < 0.000001$. Note that (3.8) yields only a zero vector if the right-hand side vector consists entirely of zeros, which explains the comments in Chapter 2 about only obtaining trivial solutions when both the function values and first derivatives are zero at the boundaries.     

\subsection{Specifying the collocation sites as Gaussian points}

In COLPNT, the $k-2$ interior collocation sites within each subinterval $[\xi_i, \xi_{i+1}]$ of the breakpoint sequence $\xi$ are specified as a fixed set of points $\rho_j$, $j = 2, \ldots, k-1$, within the interval $[-1, 1]$, such that 
\begin{equation*}
-1 < \rho_2 < \rho_3 < \cdots < \rho_{k-1} < 1
\end{equation*}
This set of points is then mapped uniformly to each $[\xi_i, \xi_{i+1}]$ using the formula 
\begin{equation}
\tau_{(i-1)(k-2) + j} = \frac{(1-\rho_j)\xi_i}{2} +  \frac{(1+\rho_j)\xi_{i+1}}{2}
\end{equation}
yielding a total of $n-2$ interior collocation sites $(\tau_2, \ldots, \tau_{n-1})$. By default, COLPNT chooses the points $\rho_j$ to be the zeros of the Legendre polynomial of degree $k-2$ (called `Gaussian points', since they are the same as the sites used for Gauss quadrature). Theorem 4.1 in \cite{deboor} shows that this choice of collocation points can significantly reduce the size of approximation errors by introducing a Legendre polynomial into their Green's function integral. Some polynomial components of the Green's function integral which are of lower degree than this Legendre polynomial will then vanish, since Legendre polynomials are orthogonal to polynomials of lower degree. The effect of this will be particularly significant at the boundaries of each subinterval $[\xi_i, \xi_{i+1}]$, producing a phenomenon called `superconvergence' there. 

To give a flavour of how this can work with Shore's equation, consider the radial equation for the ground state of the electron in the hydrogen atom in (2.10), which we will rewrite here as 

$(Ly)(x) = 0$, $Dy(0) = 2, y(10) = 0$

where $(Ly)(x) \equiv D^2y(x) + \big[\frac{2}{x} - 1\big]y(x)$. This has the exact solution
\begin{equation*}
y(x) = 2x\text{e}^{-x} 
\end{equation*}
We seek to solve this system by collocation, which means finding a pp function $y_{\triangle} \in \Pi_{<k, \xi;\nu}$ such that 

$(Ly_{\triangle})(\tau_i) = 0$, $Dy_{\triangle}(\tau_1) = 2, y_{\triangle}(\tau_n) = 0$

with $i = 2, \ldots, n-1$, $\tau_1 = 0$ and $\tau_n = 10$. For the purposes of this illustration, we will take the first two elements of the breakpoint sequence $\xi$ to be $\xi_1 = 0$ and $\xi_2 = 0.1$, so the first subinterval is $[\xi_1, \xi_2]\equiv [0, 0.1]$, and we will assume that we want to collocate at two interior sites, $\tau_2$ and $\tau_3$, as well as at the right-hand boundary of the subinterval. Thus, there are four collocation points in this setup, namely $0$, $\tau_2$, $\tau_3$ and $0.1$. Suppose further that we consider as an approximation of $y_{\triangle}$ in this subinterval the function $z(x)$ which, to second-order, is a linear interpolant of $y(x)$ passing through the two interior collocation sites $\tau_2$ and $\tau_3$ of the form
\begin{equation}
z(x) = \frac{2\tau_3(1-\tau_3)(x - \tau_2) - 2\tau_2(1-\tau_2)(x - \tau_3)}{\tau_3 - \tau_2}
\end{equation}   
Then since $y(x) = 2x\text{e}^{-x} = 2x(1 - x) + O(x^3)$, the approximation error at a site $x \in [0, 0.1]$, $x \neq \tau_2, \tau_3$, is
\begin{equation}
(y - z)(x) \approx  2x(1 - x) - z(x) = -2(x - \tau_2)(x - \tau_3)
\end{equation}
Now, the true approximation error $(y - y_{\triangle})(x)$ will satisfy a differential equation

$(L(y - y_{\triangle}))(x) = h(x)$, $D(y - y_{\triangle})(0) = 0, (y - y_{\triangle})(0.1) = 0$

for $x \in [0, 0.1]$, $x \neq \tau_2, \tau_3$, where the form of $h(x)$ depends on the form of the approximant $y_{\triangle}(x)$. This problem has a Green's function $G(x, u)$ and its solution can therefore be written as 
\begin{equation}
(y - y_{\triangle})(x) = \int_0^{0.1} du G(x, u)h(u)
\end{equation}
Suppose we now take $y_{\triangle} \approx z(x)$ in order to obtain $h(x)$ in (3.12), where $z(x)$ is the linear interpolant in (3.10) above. Then using (3.11) we have 
\begin{equation}
h(x) \approx (L(y - z))(x) = 2(x - \tau_2)(x - \tau_3) + g(x)
\end{equation}
where $g(x)$ is a function involving $x$, $\tau_2$ and $\tau_3$. Putting (3.13) into (3.12) we get
\begin{equation}
(y - y_{\triangle})(x) \approx 2\int_0^{0.1} du G(x, u)(u - \tau_2)(u - \tau_3) + \int_0^{0.1} du G(x, u)g(u)
\end{equation} 
We may now be able to reduce the size of the approximation error in (3.14) by choosing the points $\rho_2$ and $\rho_3$, as they appear in COLPNT, to be the zeros of the quadratic Legendre polynomial, i.e., $\rho_2 = -\frac{1}{\sqrt{3}}$ and $\rho_3 = \frac{1}{\sqrt{3}}$. Given any polynomial $q(x)$ of degree $1$, we will then have in the interval [-1, 1]:
\begin{equation*}
\int_{-1}^{1} du q(u)(u-\rho_2)(u-\rho_3) = \int_{-1}^{1} du q(u)\bigg(u-\frac{1}{\sqrt{3}}\bigg)\bigg(u+ \frac{1}{\sqrt{3}}\bigg) = 0
\end{equation*}
(cf. equation (4.13) in Theorem 4.1 in \cite{deboor}, p. 600). These Gaussian points will then be mapped by formula (3.9) above, with $\xi_1 = 0$ and $\xi_2 = 0.1$, to the interior collocation sites
\begin{equation*}
\tau_2 = \frac{(1+\rho_2)}{20} = \frac{\sqrt{3}-1}{20\sqrt{3}} = 0.02113248654 
\end{equation*}
and 
\begin{equation*}
\tau_2 = \frac{(1+\rho_3)}{20} = \frac{\sqrt{3}+1}{20\sqrt{3}} = 0.07886751345
\end{equation*}
Given any polynomial $q(x)$ of degree $1$, we will then have in the interval [0, 0.1]:
\begin{equation*}
\int_0^{0.1} du q(u)(1-\tau_2)(1-\tau_3) = \int_0^{0.1} du q(u)\bigg(u-\frac{\sqrt{3}-1}{20\sqrt{3}}\bigg)\bigg(u- \frac{\sqrt{3}+1}{20\sqrt{3}}\bigg) = 0
\end{equation*}
Therefore the quadratic $(x - \tau_2)(x - \tau_3) = \big(x-\frac{\sqrt{3}-1}{20\sqrt{3}}\big)\big(x- \frac{\sqrt{3}+1}{20\sqrt{3}}\big)$ in the first integral in (3.14), arising purely from specifying the collocation sites $\rho_2$ and $\rho_3$ as Gaussian points in COLPNT, can now reduce the size of the approximation error by making linear components of $G(x, u)$ vanish. This example is rather contrived, but Theorem 4.1 in \cite{deboor} shows that this idea applies more generally in linear and nonlinear collocation problems.  

As well as using the Gaussian points provided by default in COLPNT for our numerical experiments, we will also amend COLPNT to enable us to explore equally-spaced collocation points. In addition, we will explore nonlinearly distributed collocation points by calling the NEWNOT subroutine from COLLOC. The algorithm carried out by NEWNOT is described in detail in Chapter XII of \cite{deboor2}. NEWNOT works by examining the $(k-1)$-th derivative of the pp function approximation, which will always be a piecewise constant function for a pp function of order $k$, to identify any large `jumps' in this derivative at the interior breakpoints of $\xi$. If any such jump is identified, the program will alter the positions of the breakpoints so that more of the breakpoints are placed near the jump. Since the collocation sites are uniformly distributed within each subinterval of the breakpoint sequence $\xi$, this has the effect of accumulating more collocation sites near the areas where large jumps occur in the $(k-1)$-th derivative, hopefully improving the approximation accuracy there. Shore \cite{shore} and other authors were trying to achieve essentially the same thing when they re-distributed their collocation sites nonlinearly so that, for example, more collocation sites occurred near the nucleus of the hydrogen atom where the Schr\"{o}dinger wave functions tend to oscillate most sharply. Using NEWNOT in our numerical experiments is therefore an effective way to try to replicate the use of nonlinearly distributed collocation sites in the atomic theory literature.    

\chapter{Numerical results for electron wave functions in hydrogen}

In this chapter we report results for electron wave functions in the hydrogen atom. Section 4.1 reports results for different energy levels but with no angular momentum. Section 4.2 reports results with nonzero angular momentum. 

\section{Results for equations with zero angular momentum} 

\subsection{Ground state}

For the ground state electron wave function, we seek to approximate the exact solution (2.11) of the differential equation (2.10). Figure~\ref{fig:R10Solution} indicates that the box needs to have a right-hand endpoint of at least 10 (representing a distance of ten Bohr radii away from the atomic nucleus) to accommodate the right-hand boundary condition that the wave function should converge to zero at infinity. We therefore first try to implement Shore's equation (2.10) with box $[0, 10]$ and various combinations of mesh (i.e., number of divisions of the box into subintervals) and numbers of collocation sites per subinterval. The modified versions of the subroutines COLPNT and DIFEQU for this problem, and also the Maple code used for post-output processing after calling COLLOC, are provided in Appendix~\ref{fifth-appendix}.

\begin{figure}
\begin{center}
\includegraphics[width=6in]{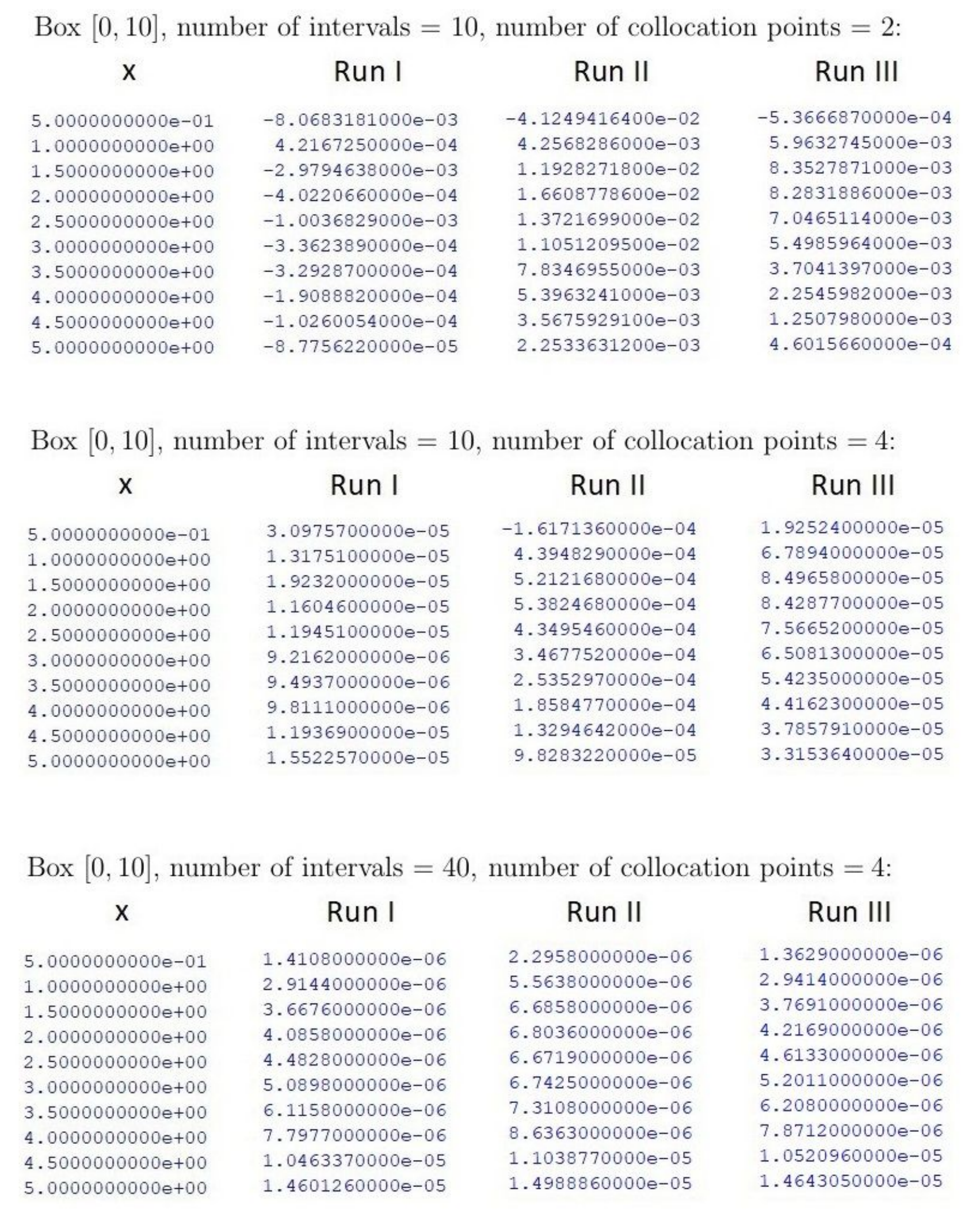}
\caption{Approximation errors at selected points $x$ for Run I (Gaussian points), Run II (equally spaced points), and Run III (nonlinear points), with box $[0, 10]$.}\label{fig:R10errorsBox10}
\end{center}
\end{figure}

For each combination of box size, mesh and number of collocation points, we conducted three runs as follows: Run I using Gaussian collocation points; Run II using equally spaced collocation points; Run III using nonlinearly spaced collocation points (produced by the NEWNOT procedure). Approximation errors at selected points were recorded for each of these runs. These are displayed in Figure~\ref{fig:R10errorsBox10}. Corresponding plots of the exact solution, the B-spline approximation and the two superimposed are shown in Figure~\ref{fig:data1010K2plots}, Figure~\ref{fig:data1010K4plots} and Figure~\ref{fig:data4010K4plots}. 

\begin{figure}
\begin{center}
\includegraphics[width=6in]{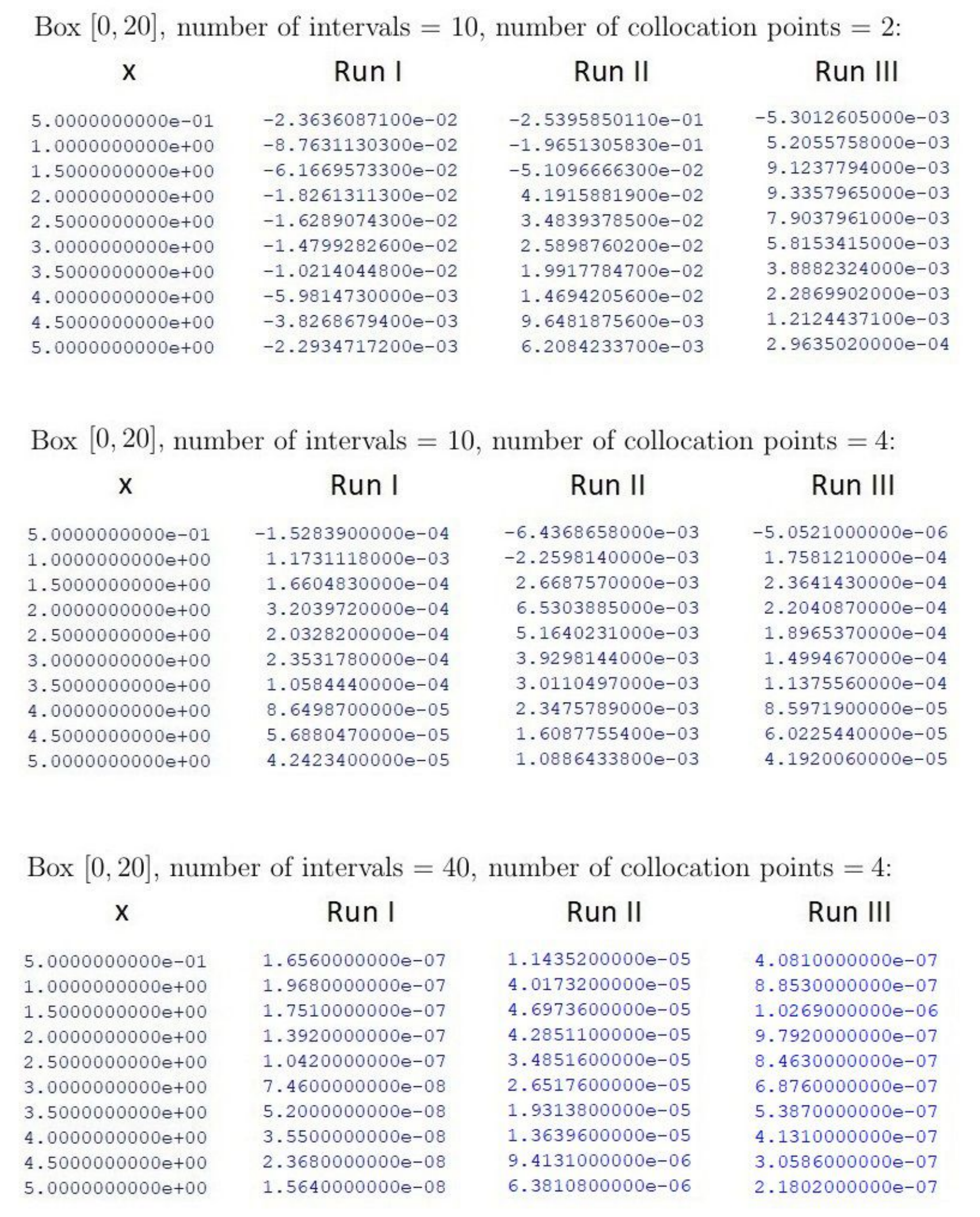}
\caption{Approximation errors at selected points $x$ for Run I (Gaussian points), Run II (equally spaced points), and Run III (nonlinear points), with box $[0, 20]$.}\label{fig:R10errorsBox20}
\end{center}
\end{figure}

\begin{figure}
\begin{center}
\includegraphics[width=6in]{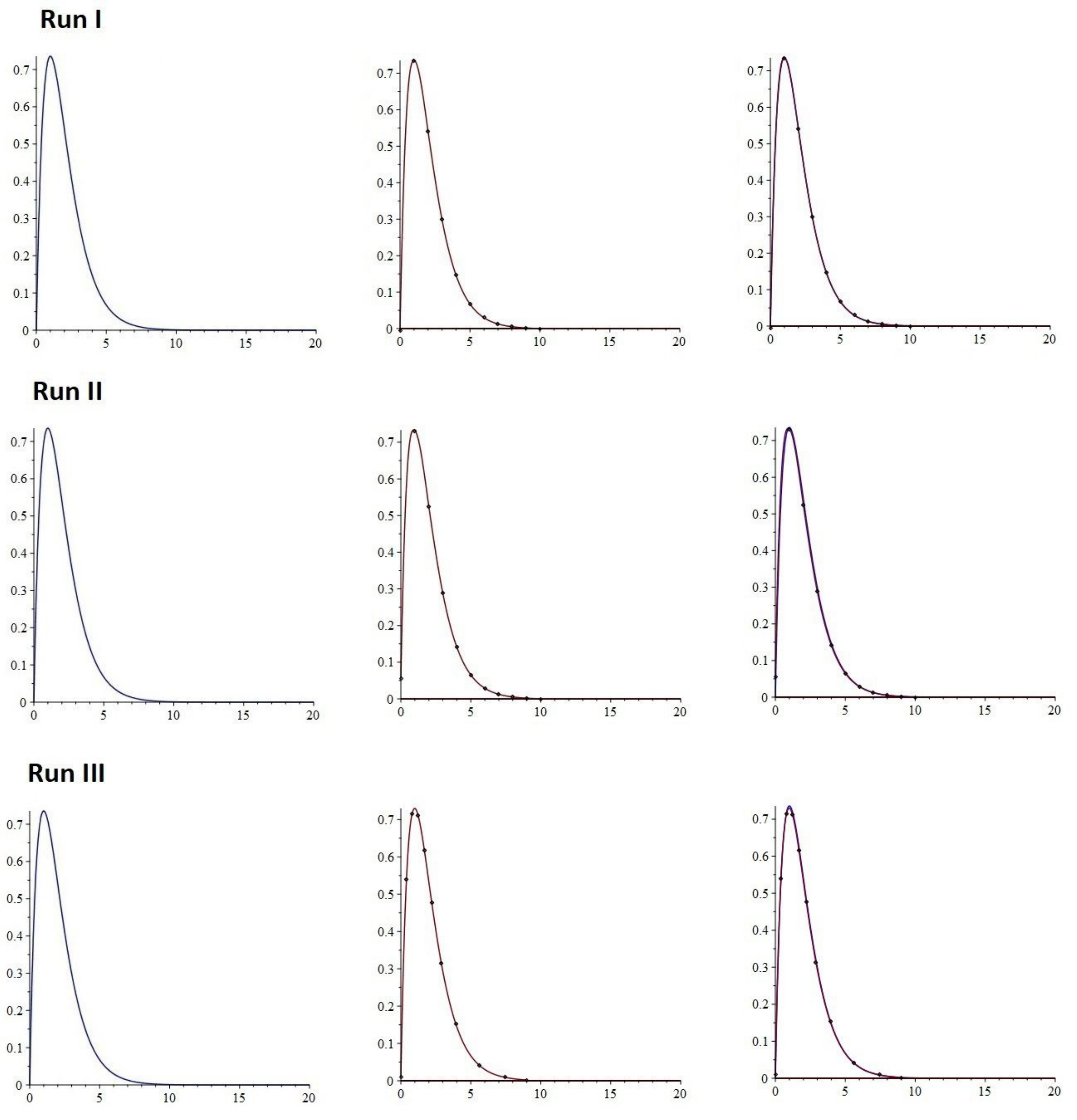}
\caption{Exact solution, B-spline approximation and the two superimposed for Run I (Gaussian points), Run II (equally spaced points), and Run III (nonlinear points), with box $[0, 10]$, 10 intervals, 2 collocation sites per interval.}\label{fig:data1010K2plots}
\end{center}
\end{figure}

\begin{figure}
\begin{center}
\includegraphics[width=6in]{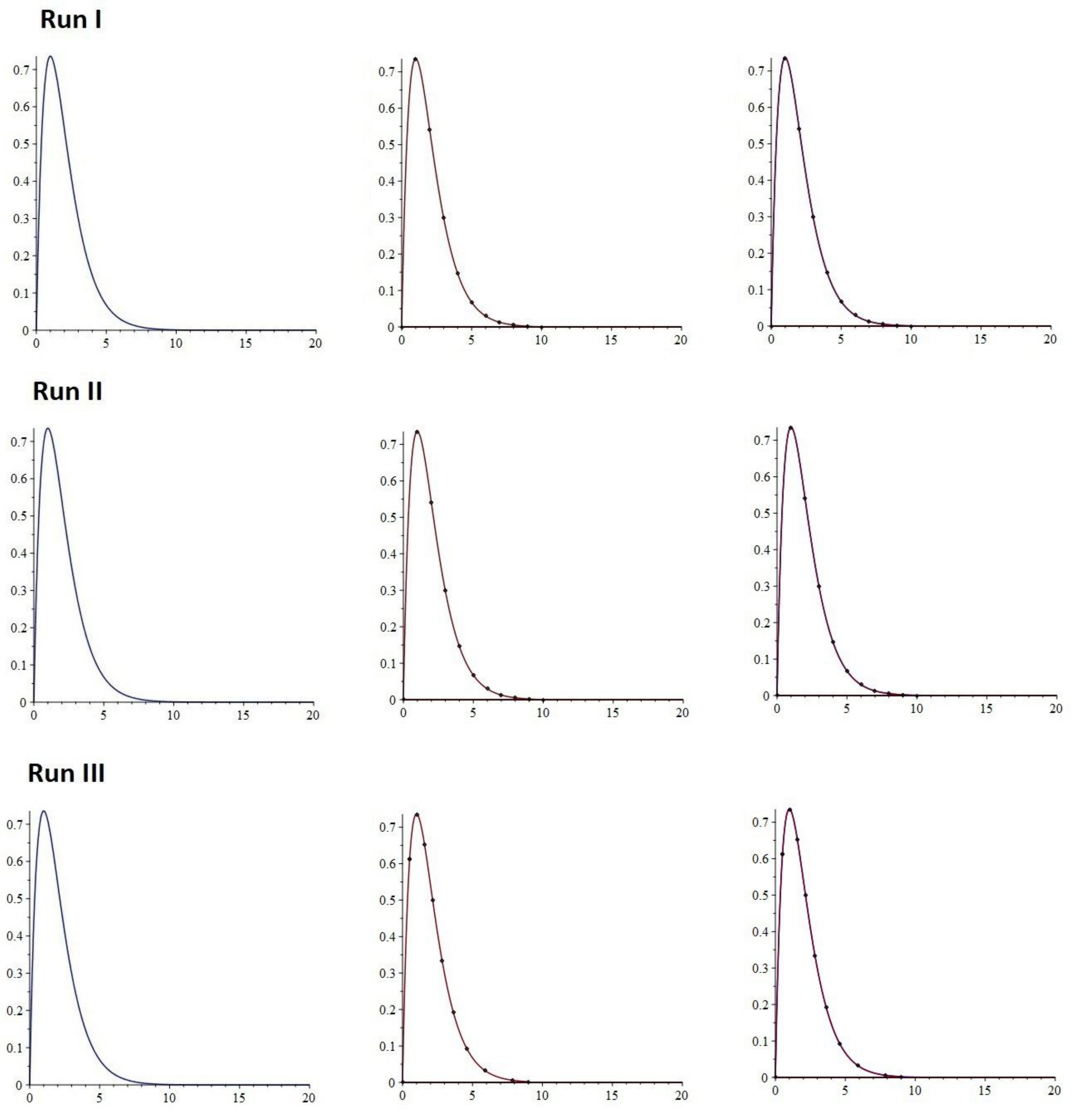}
\caption{Exact solution, B-spline approximation and the two superimposed for Run I (Gaussian points), Run II (equally spaced points), and Run III (nonlinear points), with box $[0, 10]$, 10 intervals, 4 collocation sites per interval.}\label{fig:data1010K4plots}
\end{center}
\end{figure}

\begin{figure}
\begin{center}
\includegraphics[width=6in]{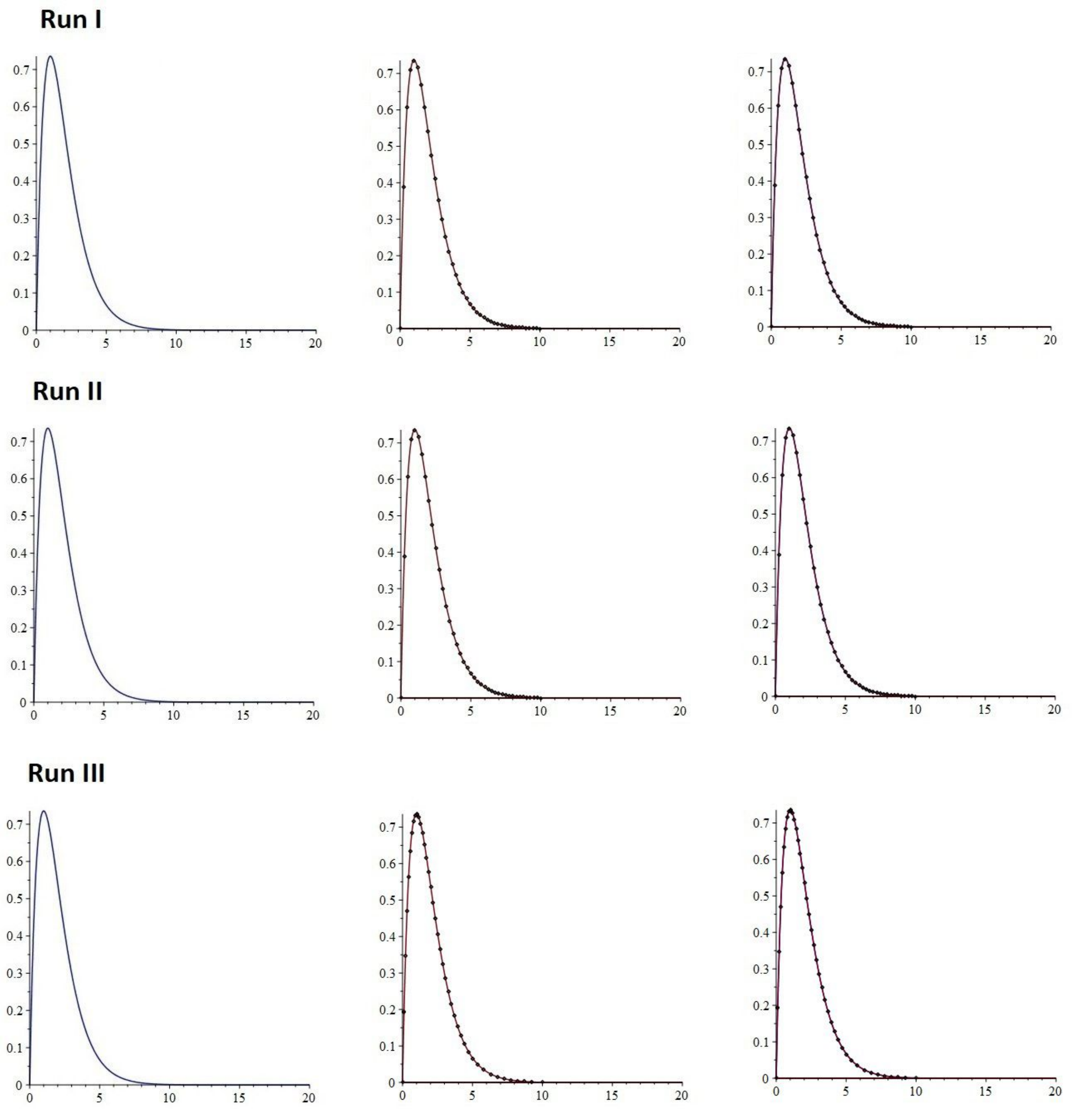}
\caption{Exact solution, B-spline approximation and the two superimposed for Run I (Gaussian points), Run II (equally spaced points), and Run III (nonlinear points), with box $[0, 10]$, 40 intervals, 4 collocation sites per interval.}\label{fig:data4010K4plots}
\end{center}
\end{figure}

\begin{figure}
\begin{center}
\includegraphics[width=6in]{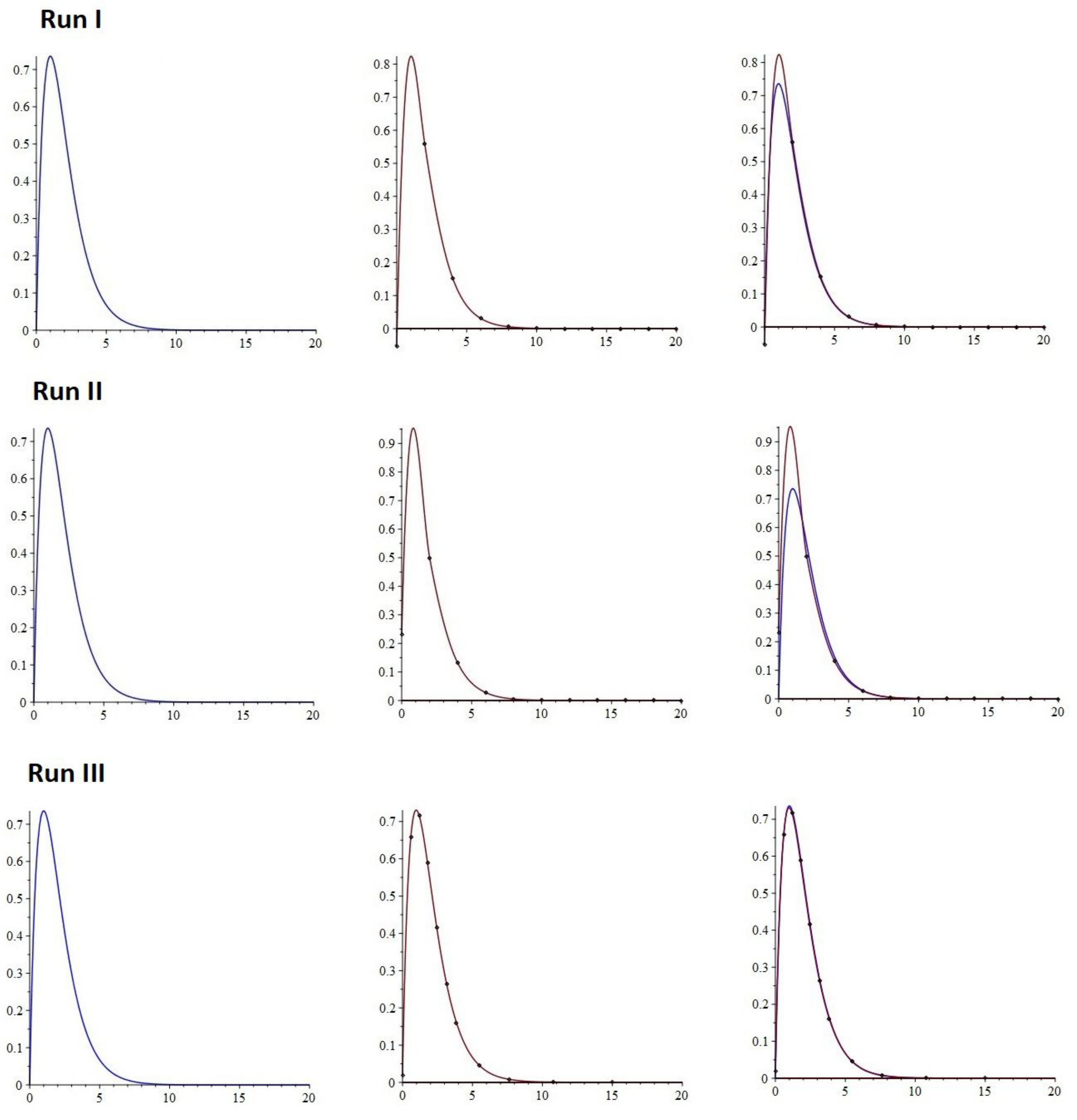}
\caption{Exact solution, B-spline approximation and the two superimposed for Run I (Gaussian points), Run II (equally spaced points), and Run III (nonlinear points), with box $[0, 20]$, 10 intervals, 2 collocation sites per interval.}\label{fig:data1020K2plots}
\end{center}
\end{figure}

\begin{figure}
\begin{center}
\includegraphics[width=6in]{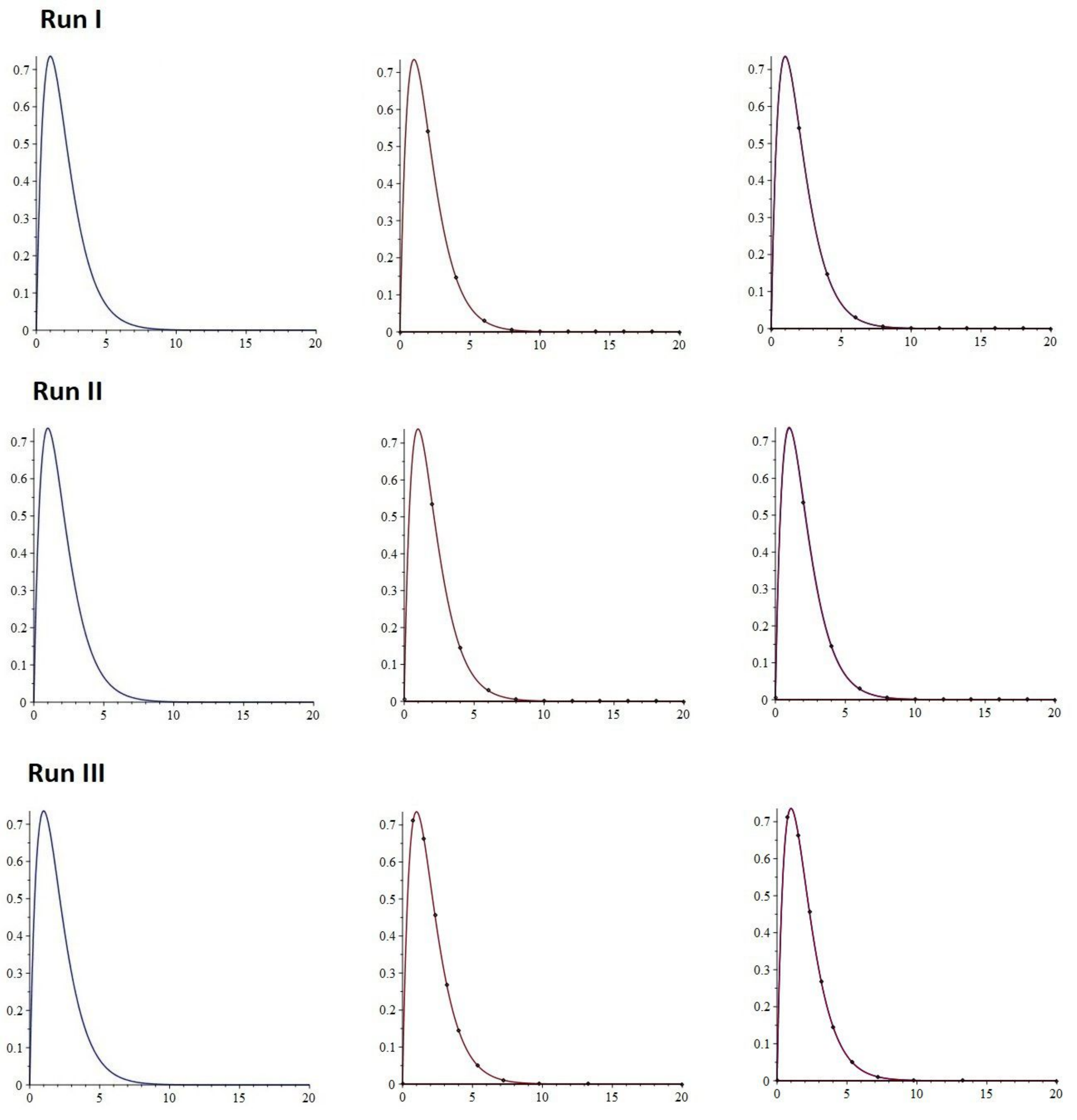}
\caption{Exact solution, B-spline approximation and the two superimposed for Run I (Gaussian points), Run II (equally spaced points), and Run III (nonlinear points), with box $[0, 20]$, 10 intervals, 4 collocation sites per interval.}\label{fig:data1020K4plots}
\end{center}
\end{figure}

\begin{figure}
\begin{center}
\includegraphics[width=6in]{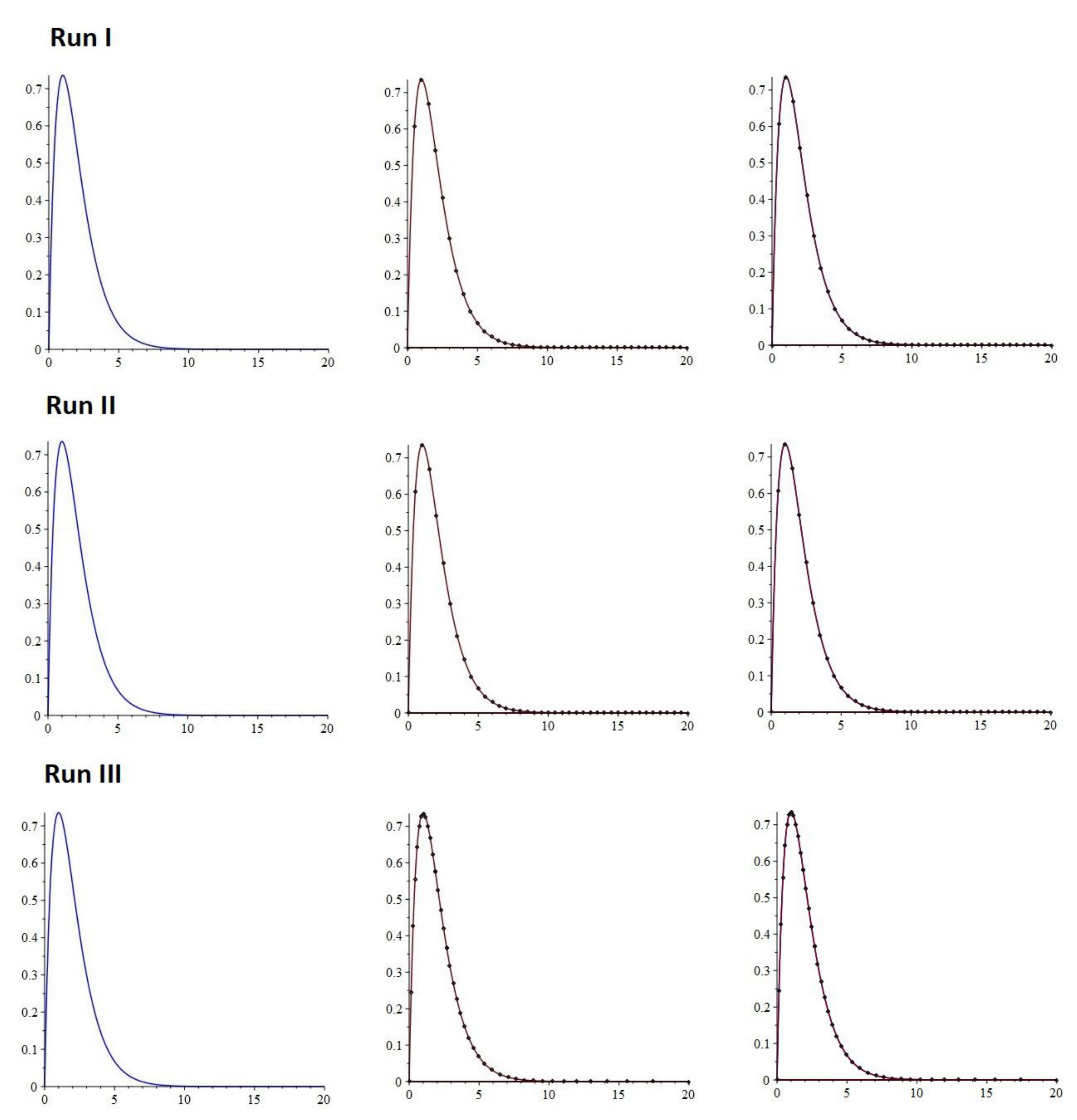}
\caption{Exact solution, B-spline approximation and the two superimposed for Run I (Gaussian points), Run II (equally spaced points), and Run III (nonlinear points), with box $[0, 20]$, 40 intervals, 4 collocation sites per interval.}\label{fig:data4020K4plots}
\end{center}
\end{figure}

To examine the effects of changing box size, we repeated these experiments with boxes of various sizes. The results for box $[0, 20]$ are reported here, as these capture the key features. The approximation errors for various combinations of mesh and numbers of collocation sites with box $[0, 20]$ are reported in Figure~\ref{fig:R10errorsBox20}, and corresponding plots of the exact solution, the B-spline approximation and the two superimposed are shown in Figure~\ref{fig:data1020K2plots}, Figure~\ref{fig:data1020K4plots} and Figure~\ref{fig:data4020K4plots}. 

In the case of box $[0, 10]$, Figures 4.3 to 4.5 show that all the approximations are visually almost indistinguishable from the exact solution, even when using only two collocation sites per subinterval. However, the approximation errors in Figure~\ref{fig:R10errorsBox10} show that equally spaced collocation points (Run II) perform consistently less well than collocation at Gaussian points (Run I) or collocation at nonlinearly distributed points produced by NEWNOT (Run III). It is also clear that collocation at Gaussian points is not noticeably inferior to collocation at nonlinearly distributed points, and actually produces slightly more accurate results with 10 subintervals and two or four collocation sites. The pattern of measurement errors also shows that significant improvements in accuracy were obtained  when the number of collocation sites was increased from two to four, and there was another significant improvement when the number of subintervals was quadrupled, from 10 subintervals to 40 subintervals.  

Changing the box size from $[0, 10]$ to $[0, 20]$ produced a noticeable \emph{worsening} of approximation accuracy in the case of 10 subintervals and two collocation sites per subinterval, as is evident from Figure~\ref{fig:data1020K2plots}. This was a surprise because the emphasis in the literature tends to be on ensuring the box size is not too small. 

However, our results show that making the box size \emph{too large} in relation to the mesh can also cause problems for approximation accuracy. The consistent picture that emerged from numerous additional experiments with different box sizes is that the mesh needs to be as fine as possible relative to the box size for greatest accuracy. It is clear from the approximation errors in Figure~\ref{fig:R10errorsBox20} that the largest approximation errors again occurred for the equally spaced collocation points, and that in the case of 10 subintervals and two or four collocation sites per subinterval, Gaussian collocation points produced larger approximation errors than the nonlinearly distributed collocation points created by the NEWNOT procedure. The superimposed plots in Figure~\ref{fig:data1020K2plots} show that in the case of 10 subintervals and two collocation sites per subinterval, neither Gaussian points nor equally spaced points produced very satisfactory approximations, while the approximation using nonlinearly spaced points is already amost indistinguishable from the exact solution at this stage. 

Increasing the number of collocation sites from two to four, still using 10 subintervals, produced a significant improvement in results. Figure~\ref{fig:data1020K4plots} shows that all the approximations become visually indistiguinshable from the exact solution when this single change is made. Again, the consistent picture that emerged from numerous additional experiments is that the number of collocation sites per interval needs to be as large as possible for greatest accuracy. Ideally, therefore, for greatest accuracy one would like to have as fine a mesh as possible and as many collocation sites per subinterval as possible, but there is a limit to how much these can be improved. For example, it was not possible to have a combination of 40 subintervals and six or more collocation sites per subinterval here, as attempts to implement such combinations led to matrix sizes for the collocation equations that were larger than those accommodated by the relevant subroutines in de Boor's package of programs.   

Nevertheless, to see how the approximations were affected by using a mesh with a significantly larger number of subintervals and polynomial approximations of higher order as determined by a higher number of collocation sites per subinterval, we implemented Shore's equation with box $[0, 20]$, a mesh of 40 intervals, and four collocation sites per interval. The polynomial pieces were quintic in this case. We again conducted three runs, Run I using Gaussian collocation points, Run II using equally spaced collocation points and Run III using nonlinearly spaced collocation points produced by the NEWNOT procedure. Approximation errors at the same points as in the previous experiments are recorded in the third table in Figure~\ref{fig:R10errorsBox20}, and plots of the exact solution, the B-spline approximation and the two superimposed for this final experiment are shown in Figure~\ref{fig:data4020K4plots}.

In this case, all three runs produced approximations which are visually indistinguishable from the exact solution. However, although the approximation errors are again largest for equally spaced collocation points, we now find that collocation at Gaussian points produces smaller approximation errors than collocation at nonlinearly spaced points. This is a reversal of the situation in the previous experiments with box $[0, 20]$ and confirms that for certain combinations of box size, mesh and order of polynomial approximants, collocation at Gaussian points is capable of producing more accurate results than the nonlinearly distributed points produced by NEWNOT. Interestingly, the results here were also more accurate for Run I and Run III than the corresponding results for box $[0, 10]$ with 40 subintervals and four collocation sites per subinterval.  

\subsection{Excited states}

The minimum required box sizes increase rapidly as we move into the excited states of the electron in the hydrogen atom. For the first excited state, corresponding to the principal quantum number $\tilde{n} = 2$, we seek to approximate the exact solution (2.13) of the differential equation (2.12). Figure~\ref{fig:R20Solution} indicates that, already, the box needs to have a right-hand endpoint about three times larger than in the ground state, around 30 (representing a distance of thirty Bohr radii away from the atomic nucleus) to accommodate the right-hand boundary condition that the wave function should converge to zero at infinity. 

\begin{figure}
\begin{center}
\includegraphics[width=6in]{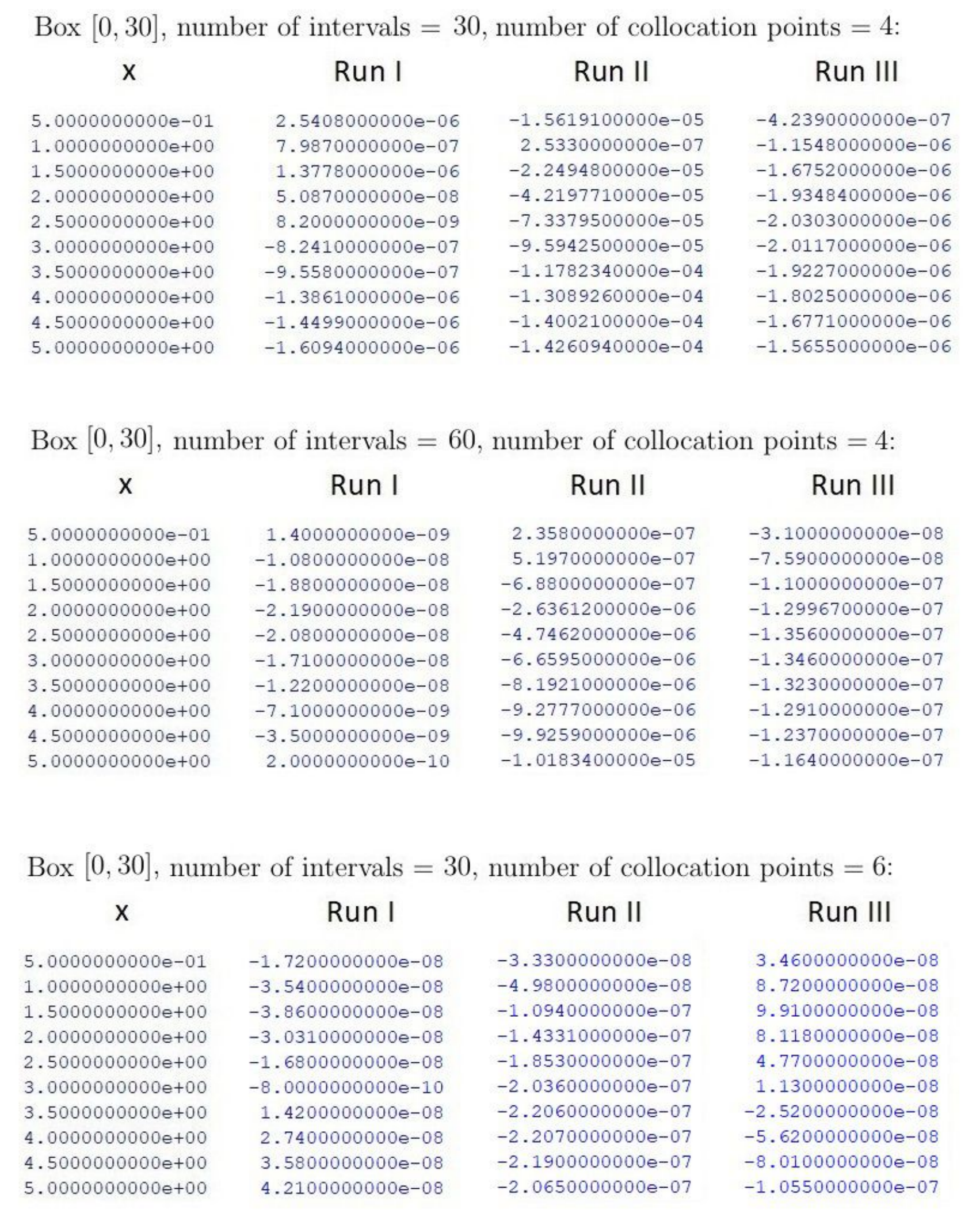}
\caption{Approximation errors at selected points $x$ for Run I (Gaussian points), Run II (equally spaced points), and Run III (nonlinear points), with box $[0, 30]$.}\label{fig:R20errorsBox30}
\end{center}
\end{figure}

\begin{figure}
\begin{center}
\includegraphics[width=6in]{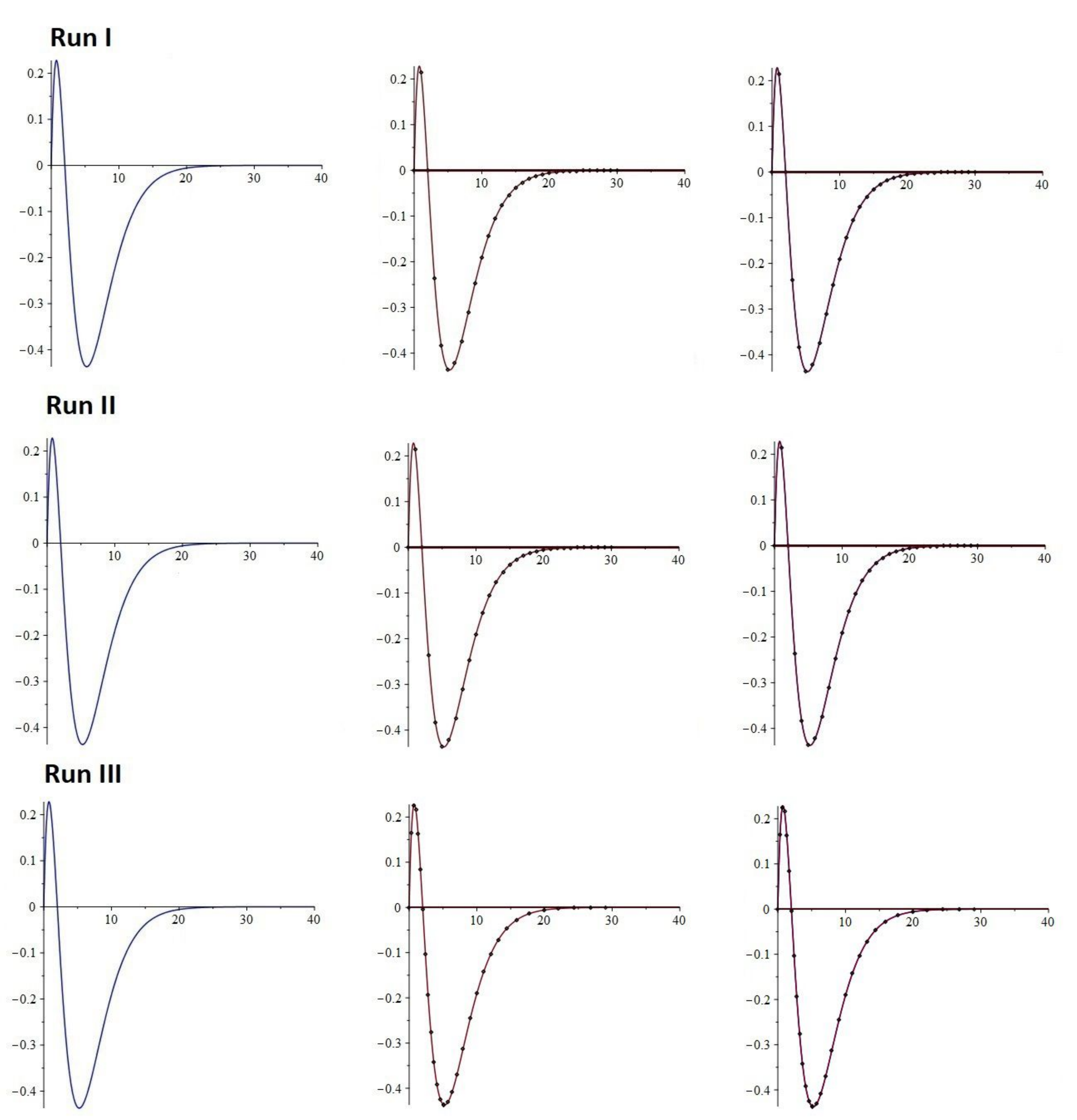}
\caption{Exact solution, B-spline approximation and the two superimposed for Run I (Gaussian points), Run II (equally spaced points), and Run III (nonlinear points), with box $[0, 30]$, 30 intervals, 4 collocation sites per interval.}\label{fig:data3030K4plots}
\end{center}
\end{figure}

\begin{figure}
\begin{center}
\includegraphics[width=6in]{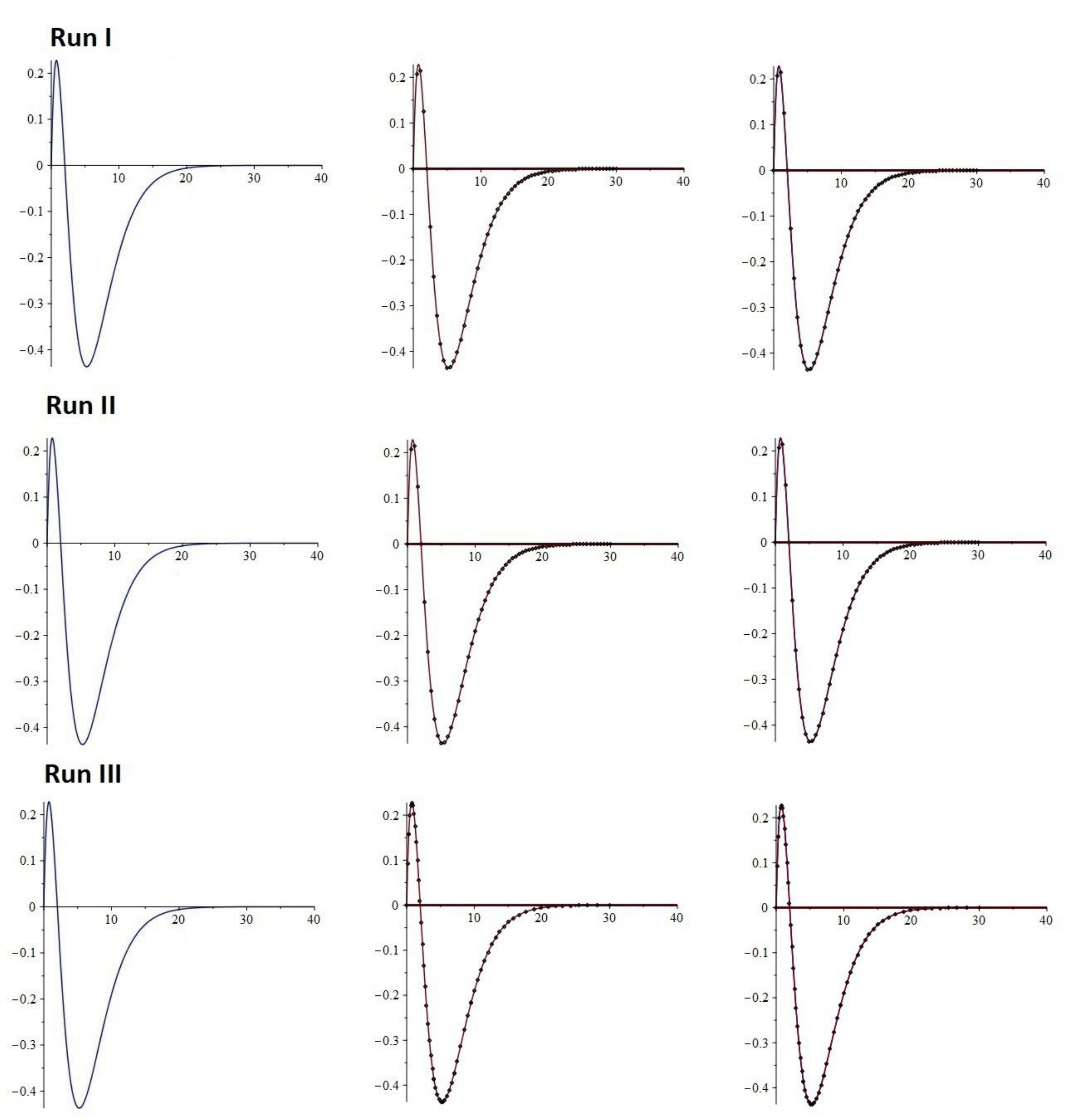}
\caption{Exact solution, B-spline approximation and the two superimposed for Run I (Gaussian points), Run II (equally spaced points), and Run III (nonlinear points), with box $[0, 30]$, 60 intervals, 4 collocation sites per interval.}\label{fig:data6030K4plots}
\end{center}
\end{figure}

\begin{figure}
\begin{center}
\includegraphics[width=6in]{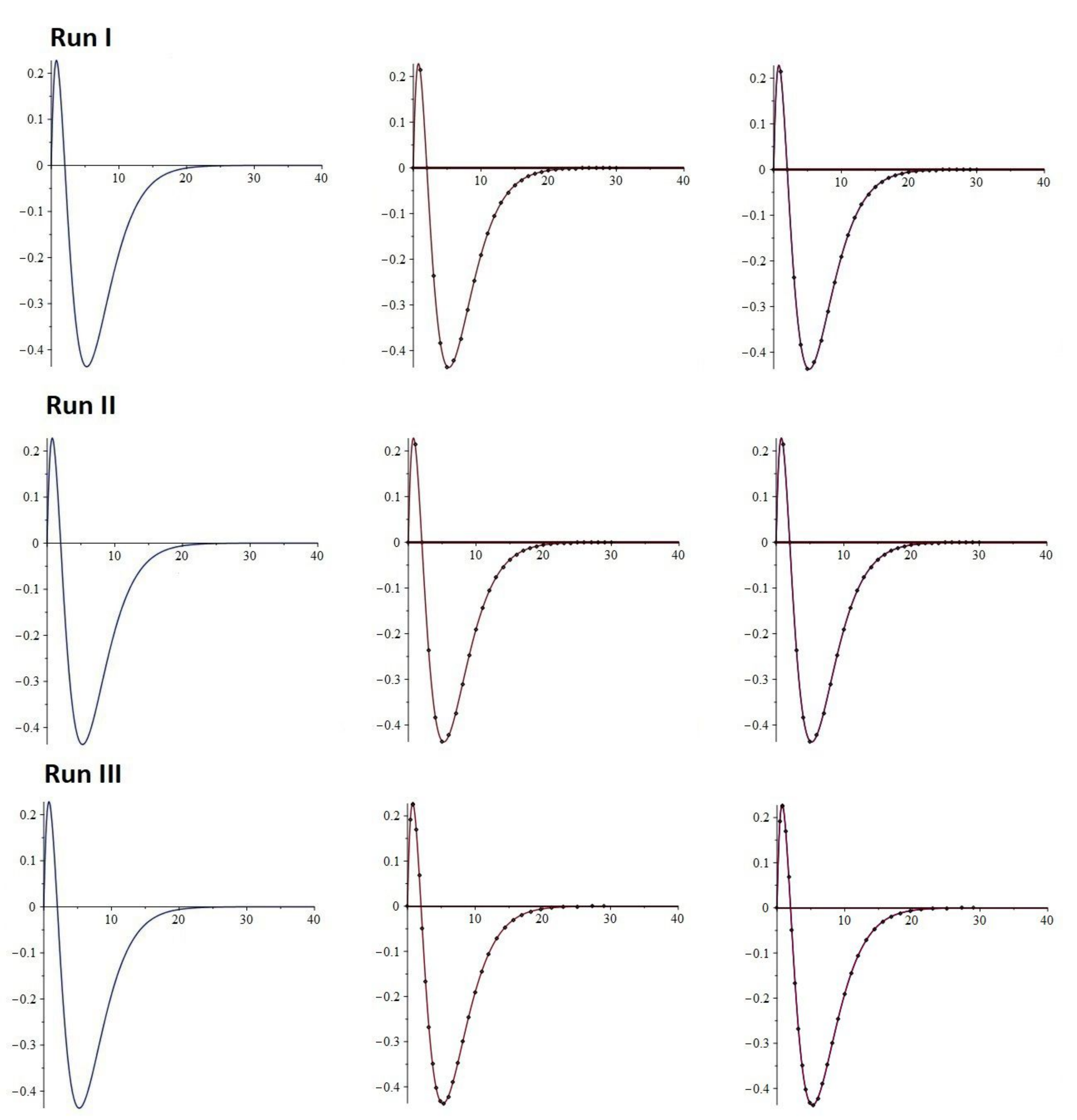}
\caption{Exact solution, B-spline approximation and the two superimposed for Run I (Gaussian points), Run II (equally spaced points), and Run III (nonlinear points), with box $[0, 30]$, 30 intervals, 6 collocation sites per interval.}\label{fig:data3030K6plots}
\end{center}
\end{figure}

In order to compare the improvements in accuracy obtained by increasing the number of subintervals (i.e., increasing the number of polynomial pieces in the approximation) versus increasing the order of each of the polynomial pieces (i.e. increasing the number of collocation sites per subinterval),  we report here the results of three experiments implementing Shore's equation (2.12) with box $[0, 30]$: one with 30 subintervals and four collocation sites per subinterval; one with 60 subintervals and four collocation sites per subinterval (i.e., doubling the number of polynomial pieces, keeping the number of collocation sites the same); and one with 30 subintervals but six collocation sites per subinterval (i.e., inceasing the order of the polynomial pieces from quintics to heptics, while keeping the number of polynomial pieces the same). The approximation errors in each experiment for Run I using Gaussian collocation points, Run II using equally spaced collocation points and Run III using nonlinearly spaced collocation points are displayed in Figure~\ref{fig:R20errorsBox30}. Corresponding plots of the exact solution, the B-spline approximation and the two superimposed are shown in Figure~\ref{fig:data3030K4plots}, Figure~\ref{fig:data6030K4plots} and Figure~\ref{fig:data3030K6plots}. 

Figures 4.10 to 4.12 show that all the approximations are visually almost indistinguishable from the exact solution, even when using only 30 subintervals and four collocation sites per subinterval. However, as in previous experiments, the approximation errors in Figure~\ref{fig:R20errorsBox30} show that equally spaced collocation points (Run II) performed consistently less well than Gaussian collocation points (Run I) or collocation at nonlinearly distributed points produced by NEWNOT (Run III). It is also again clear that collocation at Gaussian points performed just as well or better than collocation at nonlinearly distributed points in these experiments. 

The measurement errors show that significant improvements in accuracy were obtained when the number of subintervals (i.e., number of polynomial pieces) was doubled from 30 to 60 keeping the number of collocation sites the same. However, similar improvements were obtained when the number of collocation sites was increased from four to six, keeping the number of polynomial pieces the same. There seems to be little to choose between these two approaches in terms of increasing the accuracy of approximations here. 
\begin{figure}
\begin{center}
\includegraphics[width=6in]{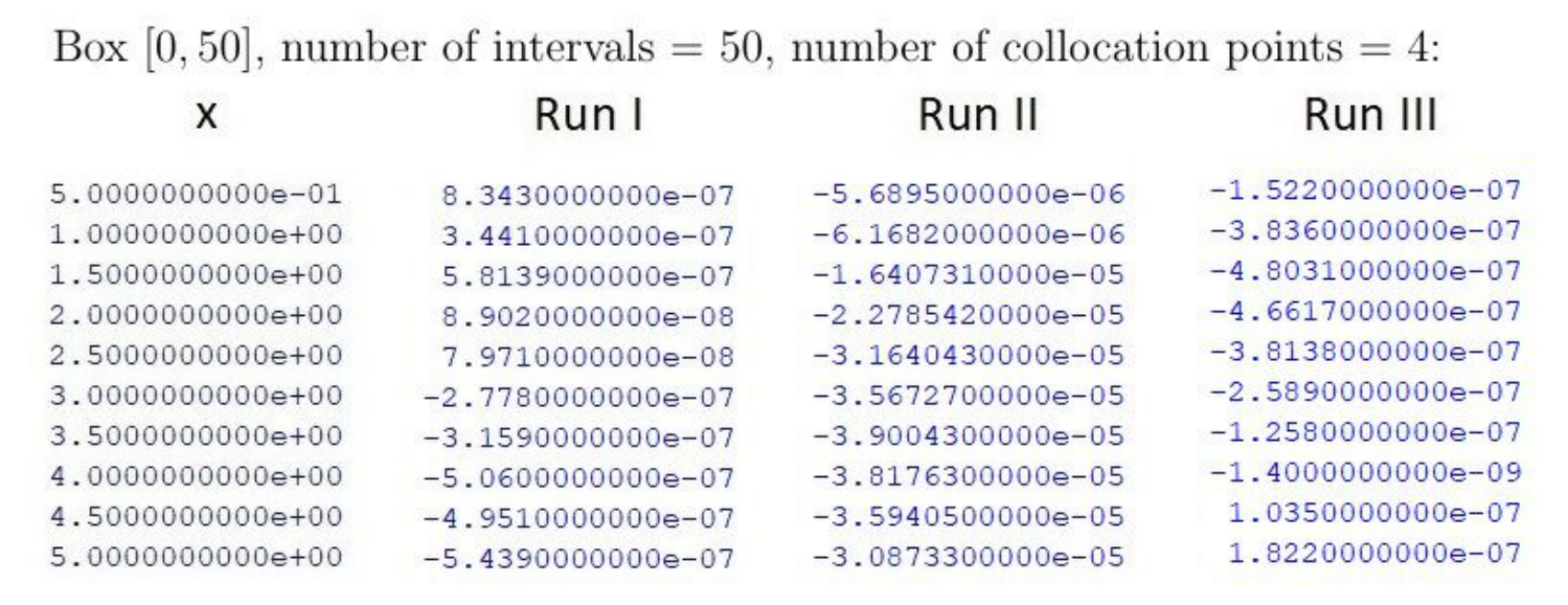}
\caption{Approximation errors at selected points $x$ for Run I (Gaussian points), Run II (equally spaced points), and Run III (nonlinear points), with box $[0, 50]$.}\label{fig:R30errorsBox50}
\end{center}
\end{figure}  

A further sharp increase in box size is required when we move to the second excited state, corresponding to principal quantum number $\tilde{n} = 3$. Here we are trying to approximate the exact solution (2.15) of the differential equation (2.14). Figure~\ref{fig:R30Solution} indicates that the box now needs to have a right-hand endpoint around 50, representing a distance of fifty Bohr radii away from the atomic nucleus.  

\begin{figure}
\begin{center}
\includegraphics[width=6in]{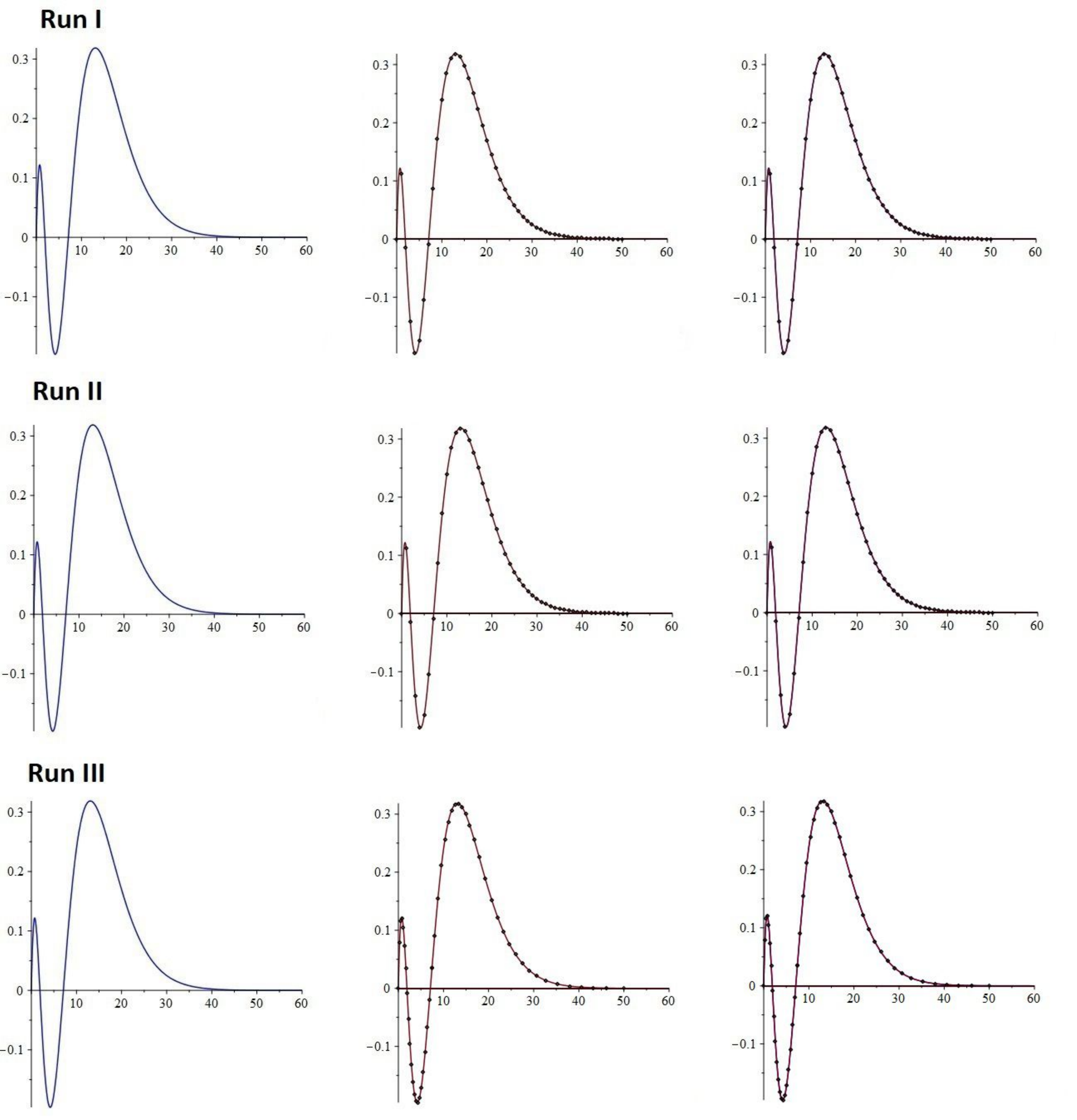}
\caption{Exact solution, B-spline approximation and the two superimposed for Run I (Gaussian points), Run II (equally spaced points), and Run III (nonlinear points), with box $[0, 50]$, 50 intervals, 4 collocation sites per interval.}\label{fig:data5050K4plots}
\end{center}
\end{figure}

We report here the results of an experiment to approximate the exact solution for the second excited state with box $[0, 50]$, 50 subintervals and four collocation sites per subinterval. Approximation errors are recorded in Figure~\ref{fig:R30errorsBox50} for Run I using Gaussian collocation points, Run II using equally spaced points and Run III using nonlinearly spaced points. Plots of the exact solution, the B-spline approximation and the two superimposed for this experiment are shown in Figure~\ref{fig:data5050K4plots}.

We observe similar patterns to those in the previous experiments. All the approximations are visually close to the exact solution, but the approximation errors in Figure~\ref{fig:R30errorsBox50} show that equally spaced collocation points perform less well than Gaussian points or nonlinearly distributed points. The performance of Gaussian collocation points is more or less on a par with collocation at nonlinearly distributed points in terms of approximation accuracy.

\section{Results for equations incorporating angular momentum} 

As discussed in subsection 2.1.3, the inclusion of angular momentum in the radial Schr\"{o}dinger equations posed a numerical difficulty causing the COLLOC procedure to find only trivial solutions. In the case $\tilde{n} = 2$, $l = 1$, this had to be overcome by transforming the original differential equation (2.17) into differential equation (2.22) instead, for which de Boor's methodology is able to provide nontrivial solutions. The transformation can then easily be reversed using the resulting output to obtain the desired approximations of the exact solution (2.18). Therefore here we report our approximation of (2.19) from the differential equation (2.22), from which we obtained the desired approximation of (2.18) using $F(x) = \tilde{F}\big(\frac{x}{2}\big)$. The modified version of the subroutine DIFEQU for this problem, and also the Maple code used for post-output processing after calling COLLOC, are provided in Appendix~\ref{sixth-appendix}. 

\begin{figure}
\begin{center}
\includegraphics[width=6in]{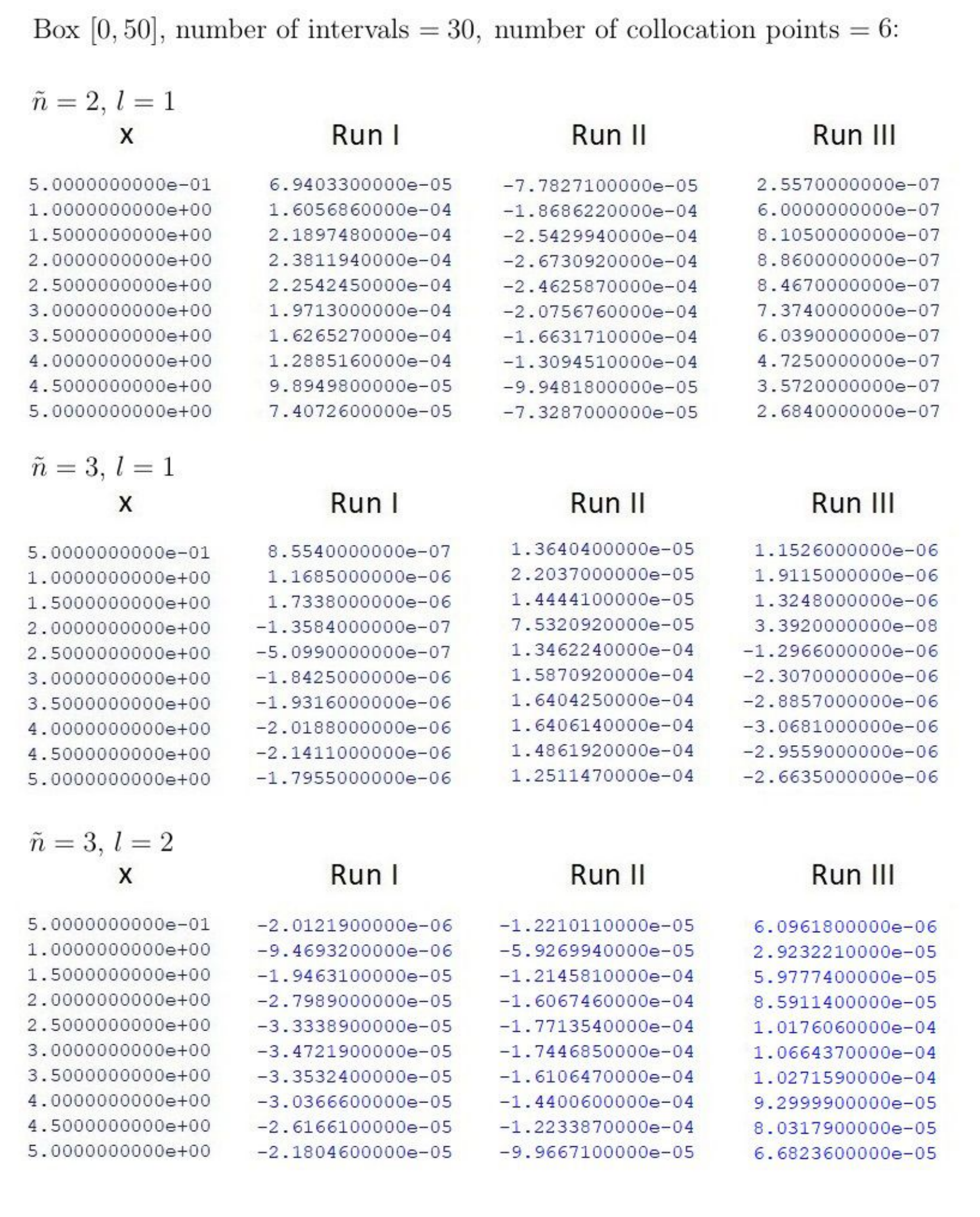}
\caption{Approximation errors at selected points $x$ for Run I (Gaussian points), Run II (equally spaced points), and Run III (nonlinear points), with box $[0, 50]$.}\label{fig:FTildeerrorsBox50}
\end{center}
\end{figure}

\begin{figure}
\begin{center}
\includegraphics[width=6in]{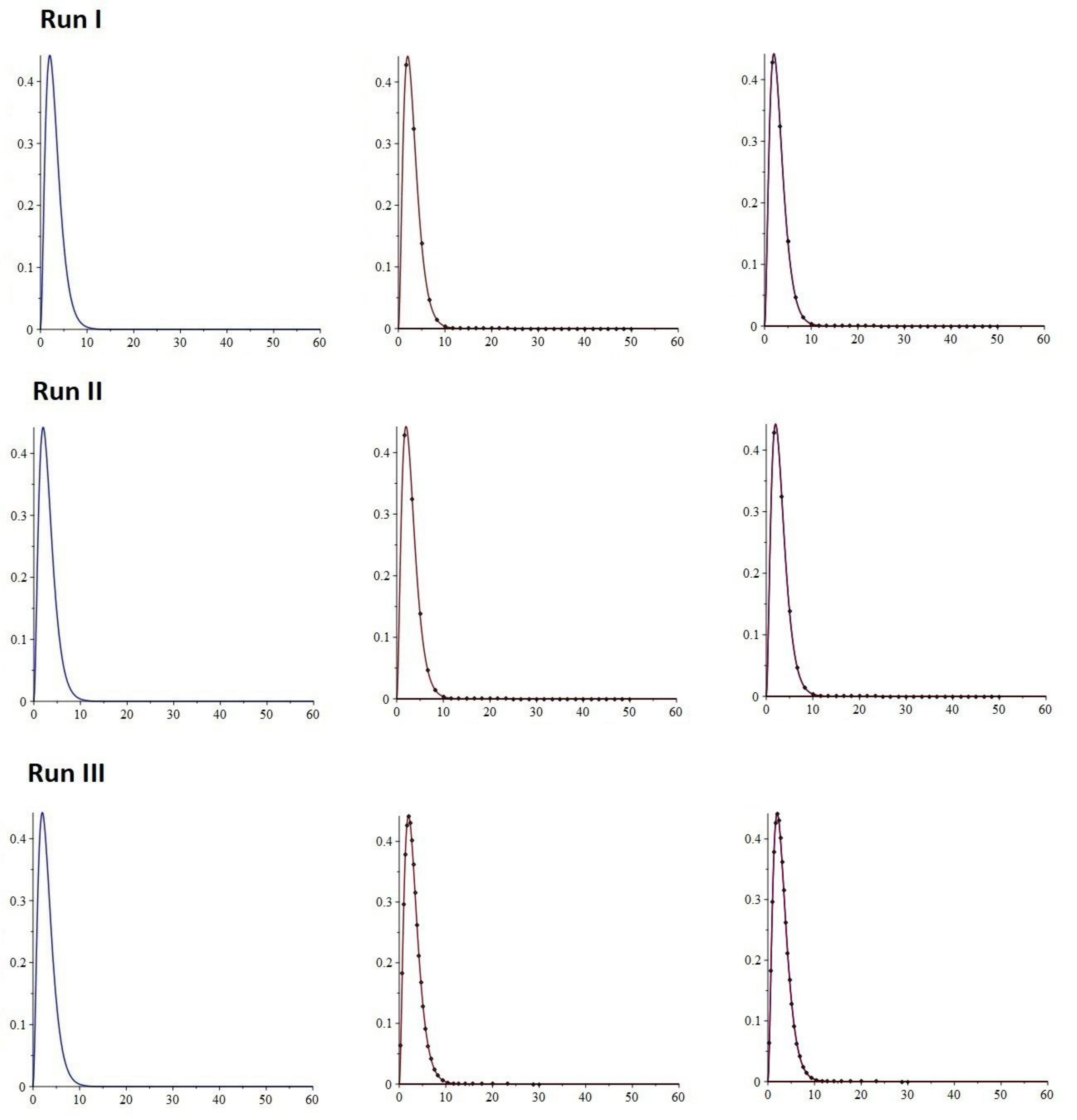}
\caption{Exact solution (2.19), B-spline approximation and the two superimposed for Run I (Gaussian points), Run II (equally spaced points), and Run III (nonlinear points), with box $[0, 50]$, 30 intervals, 6 collocation sites per interval.}\label{fig:data5030K6plots}
\end{center}
\end{figure}

\begin{figure}
\begin{center}
\includegraphics[width=6in]{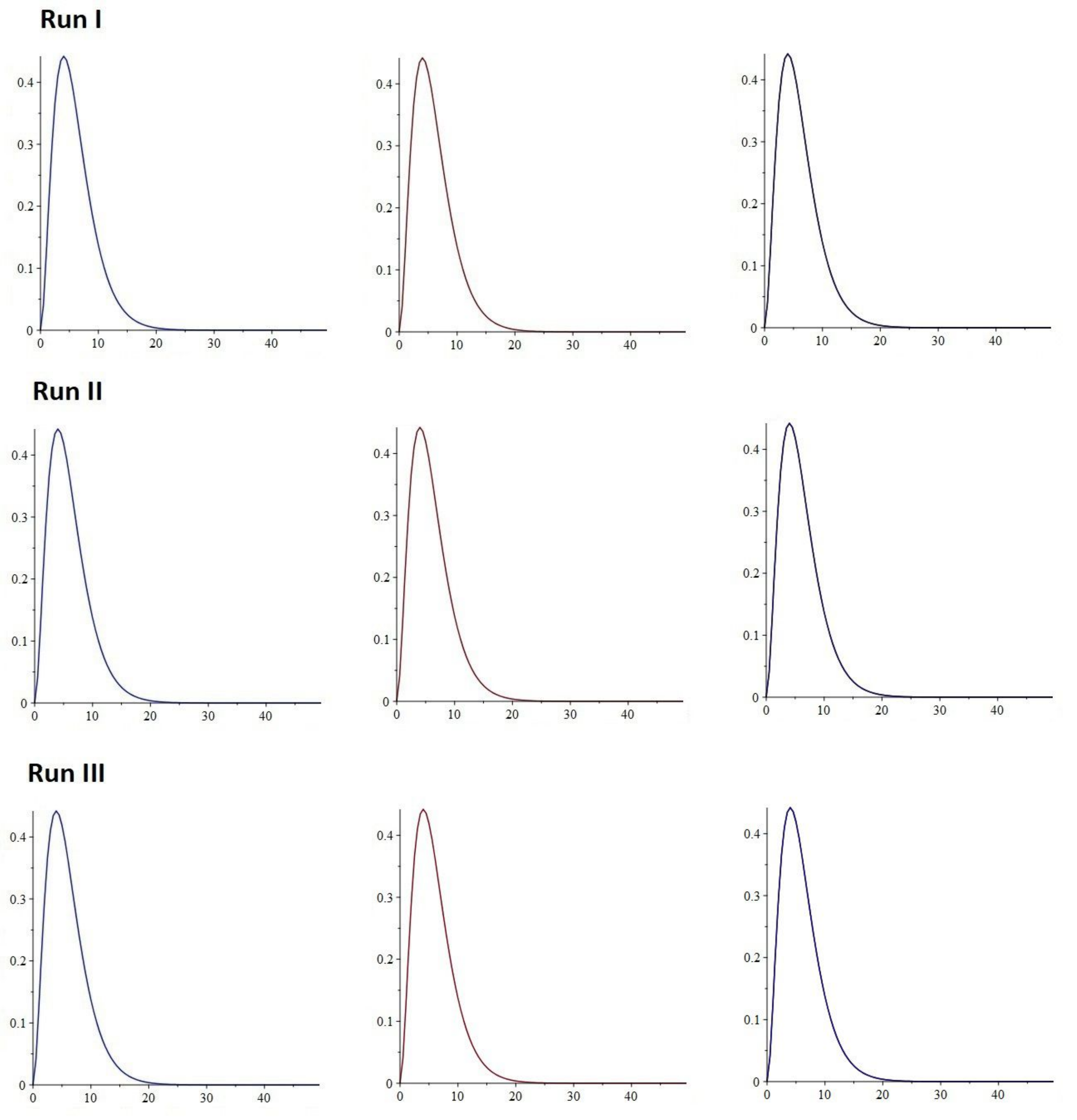}
\caption{Exact solution (2.18) for $\tilde{n} = 2$, $l = 1$, derived B-spline approximation and the two superimposed for Run I (Gaussian points), Run II (equally spaced points), and Run III (nonlinear points), with box $[0, 50]$, 30 intervals, 6 collocation sites per interval.}\label{fig:R21plots}
\end{center}
\end{figure}

\begin{figure}
\begin{center}
\includegraphics[width=6in]{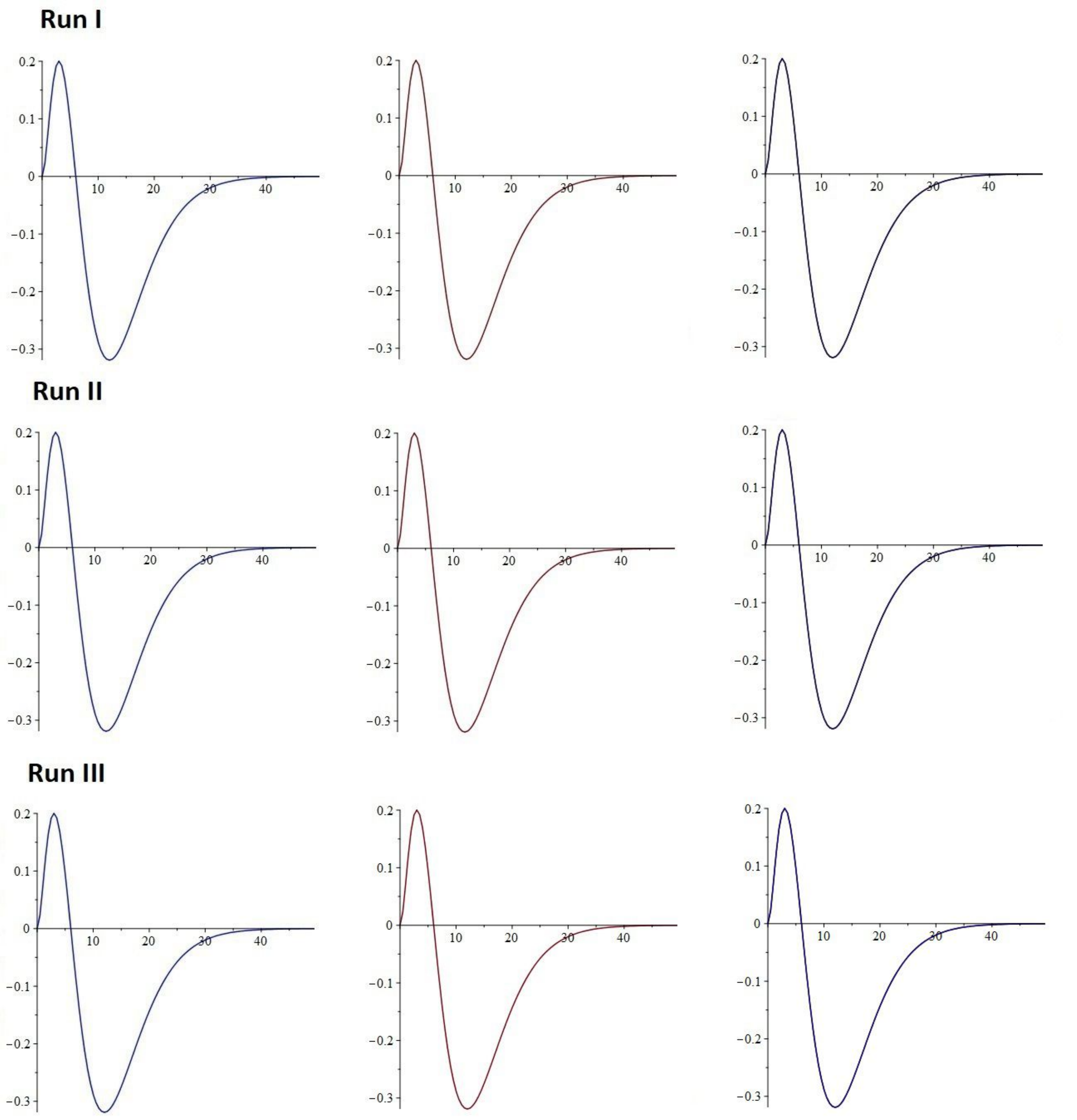}
\caption{Exact solution (2.25) for $\tilde{n} = 3$, $l = 1$, derived B-spline approximation and the two superimposed for Run I (Gaussian points), Run II (equally spaced points), and Run III (nonlinear points), with box $[0, 50]$, 30 intervals, 6 collocation sites per interval.}\label{fig:R31plots}
\end{center}
\end{figure}

\begin{figure}
\begin{center}
\includegraphics[width=6in]{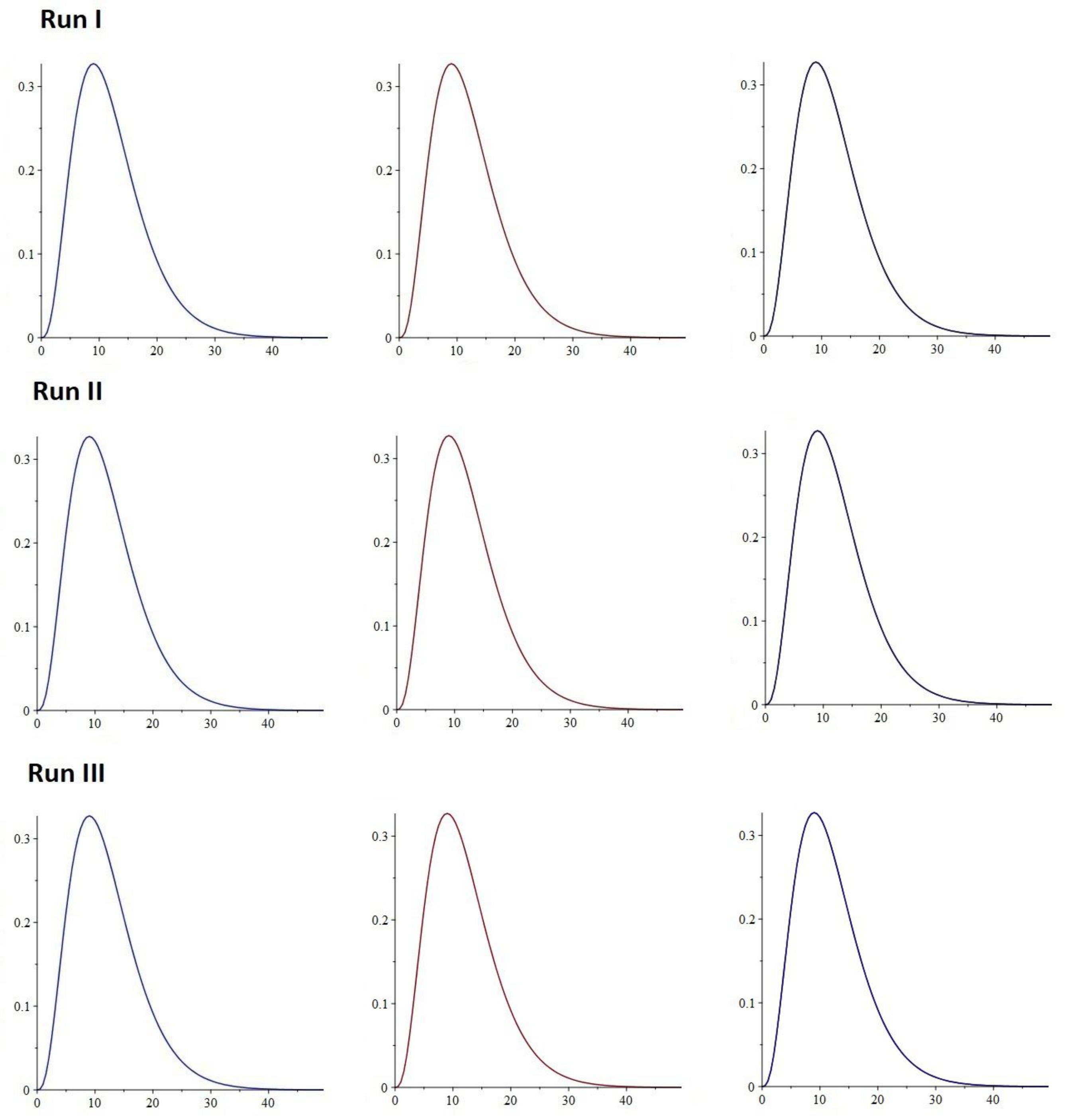}
\caption{Exact solution (2.27) for $\tilde{n} = 3$, $l = 2$, derived B-spline approximation and the two superimposed for Run I (Gaussian points), Run II (equally spaced points), and Run III (nonlinear points), with box $[0, 50]$, 30 intervals, 6 collocation sites per interval.}\label{fig:R32plots}
\end{center}
\end{figure}

We used box $[0, 50]$, 30 subintervals and 6 collocation sites per subinterval as this combination gave the most accurate results for all three runs. We repeated the same kind of approach for the cases $\tilde{n} = 3$, $l = 1$ and $\tilde{n} = 3$, $l = 2$. Approximation errors for all three cases with nonzero angular momentum are recorded in Figure~\ref{fig:FTildeerrorsBox50} for Run I using Gaussian collocation points, Run II using equally spaced points and Run III using nonlinearly spaced points. Plots of the exact solution (2.19), the B-spline approximation and the two superimposed for the $\tilde{n} = 2$, $l = 1$ experiment are shown in Figure~\ref{fig:data5030K6plots}. Plots of the derived approximations of (2.18) are shown in Figure~\ref{fig:R21plots}. Finally, plots of the derived approximations of exact solutions (2.25) and (2.27) for the cases $\tilde{n} = 3$, $l = 1$ and $\tilde{n} = 3$, $l = 2$, respectively, are shown in Figure~\ref{fig:R31plots} and Figure~\ref{fig:R32plots}.  

Figures 4.16 to 4.19 show that all the approximations are visually almost indistinguishable from the corresponding exact solutions, but differences in performance between the different patterns of collocation points become clear when looking at the approximation errors in Figure~\ref{fig:FTildeerrorsBox50}. 

As in previous experiments, the approximation errors in Figure~\ref{fig:FTildeerrorsBox50} show that equally spaced collocation points (Run II) performed consistently less well than Gaussian collocation points (Run I) or collocation at nonlinearly distributed points produced by NEWNOT (Run III). It is also clear that collocation at Gaussian points performed just as well or better than collocation at nonlinearly distributed points in the experiments for the cases $\tilde{n} = 3$, $l = 1$ and $\tilde{n} = 3$, $l = 2$, confirming again that for certain combinations of box size, mesh and order of polynomial approximants in quantum systems, collocation at Gaussian points is capable of producing more accurate results than the nonlinearly distributed points produced by NEWNOT. 

However, the first table in Figure~\ref{fig:FTildeerrorsBox50} shows that there is a significant reversal in the case $\tilde{n} = 2$, $l = 1$, with both Gaussian collocation points and equally spaced points performing relatively poorly compared to the high approximation accuracy achieved with nonlinearly distributed collocation points produced by NEWNOT. This is reminiscent of the situation encountered earlier in the equations without angular momentum with box size $[0, 20]$, 10 subintervals and two collocation sites per subinterval, in which both Gaussian collocation points and equally spaced points produced larger approximation errors than nonlinearly distributed points produced by NEWNOT (see the first table in Figure~\ref{fig:R10errorsBox20}). The situation is even more pronounced here. Additional experiments both in the previous section and here showed that this significantly better performance by nonlinearly distributed collocation points compared to Gaussian points tends to occur sometimes in situations in which the mesh is relatively coarse (i.e., too few subintervals) compared to the box size.  


\chapter{Numerical results for the nonlinear Schr\"{o}dinger equation}

In this chapter we report results for the nonlinear Schr\"{o}dinger equation with different values for the perturbation parameter, $\epsilon$. Section 5.1 reports results for $\epsilon = 0.1$, $\epsilon = 0.05$ and $\epsilon = 0.025$. Section 5.2 reports results for $\epsilon = 0.01$, $\epsilon = 0.005$ and $\epsilon = 0.001$. 

\section{Results for $\epsilon = 0.1$, $\epsilon = 0.05$ and $\epsilon = 0.025$} 

Here we seek to approximate the exact solution (2.36) of the cubic Schr\"{o}dinger equation (2.34) with box $[0, 1]$, 20 subintervals and 6 collocation sites per subinterval. Therefore we are using 20 polynomial pieces, each of order 8, i.e., the polynomials are heptics. The modified version of the subroutine DIFEQU for this problem, and also the Maple code used for post-output processing after calling COLLOC, are provided in Appendix~\ref{seventh-appendix}.

\begin{figure}
\begin{center}
\includegraphics[width=6in]{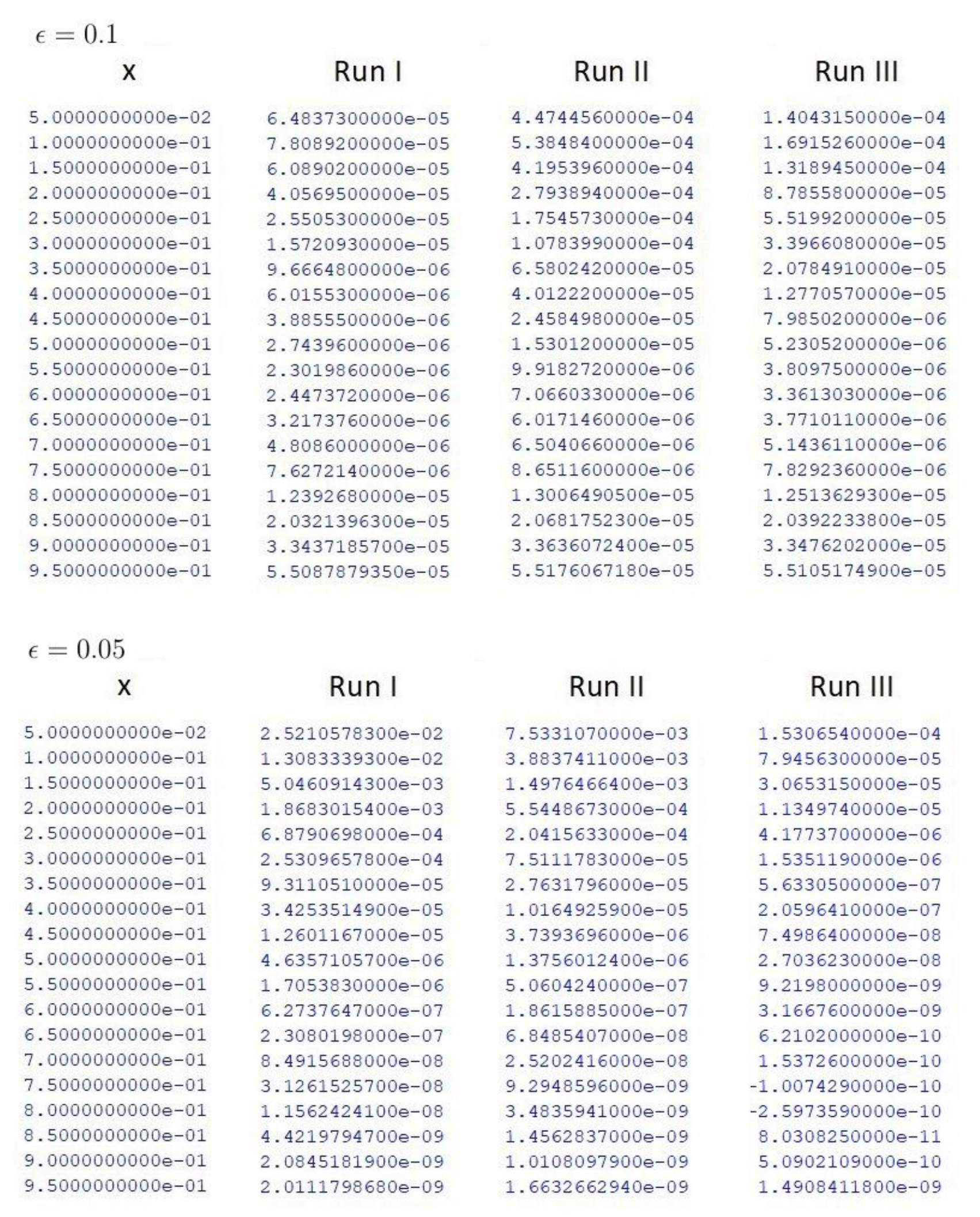}
\caption{Approximation errors at selected points $x$ for Run I (Gaussian points), Run II (equally spaced points), and Run III (nonlinear points), with box $[0, 1]$.}\label{fig:EPSILON01005errors}
\end{center}
\end{figure}

\begin{figure}
\begin{center}
\includegraphics[width=6in]{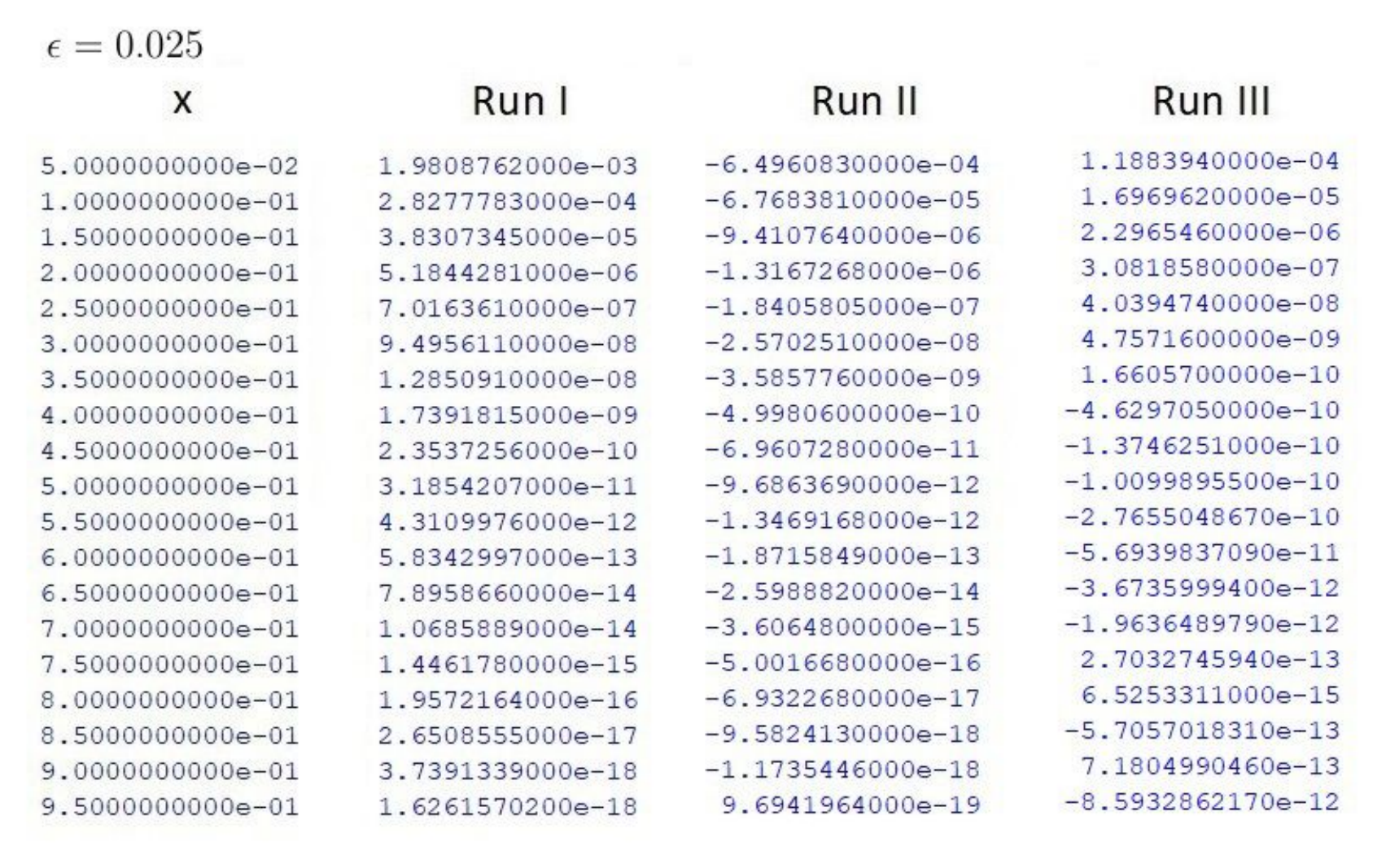}
\caption{Approximation errors at selected points $x$ for Run I (Gaussian points), Run II (equally spaced points), and Run III (nonlinear points), with box $[0, 1]$.}\label{fig:EPSILON0025errors}
\end{center}
\end{figure}

\begin{figure}
\begin{center}
\includegraphics[width=6in]{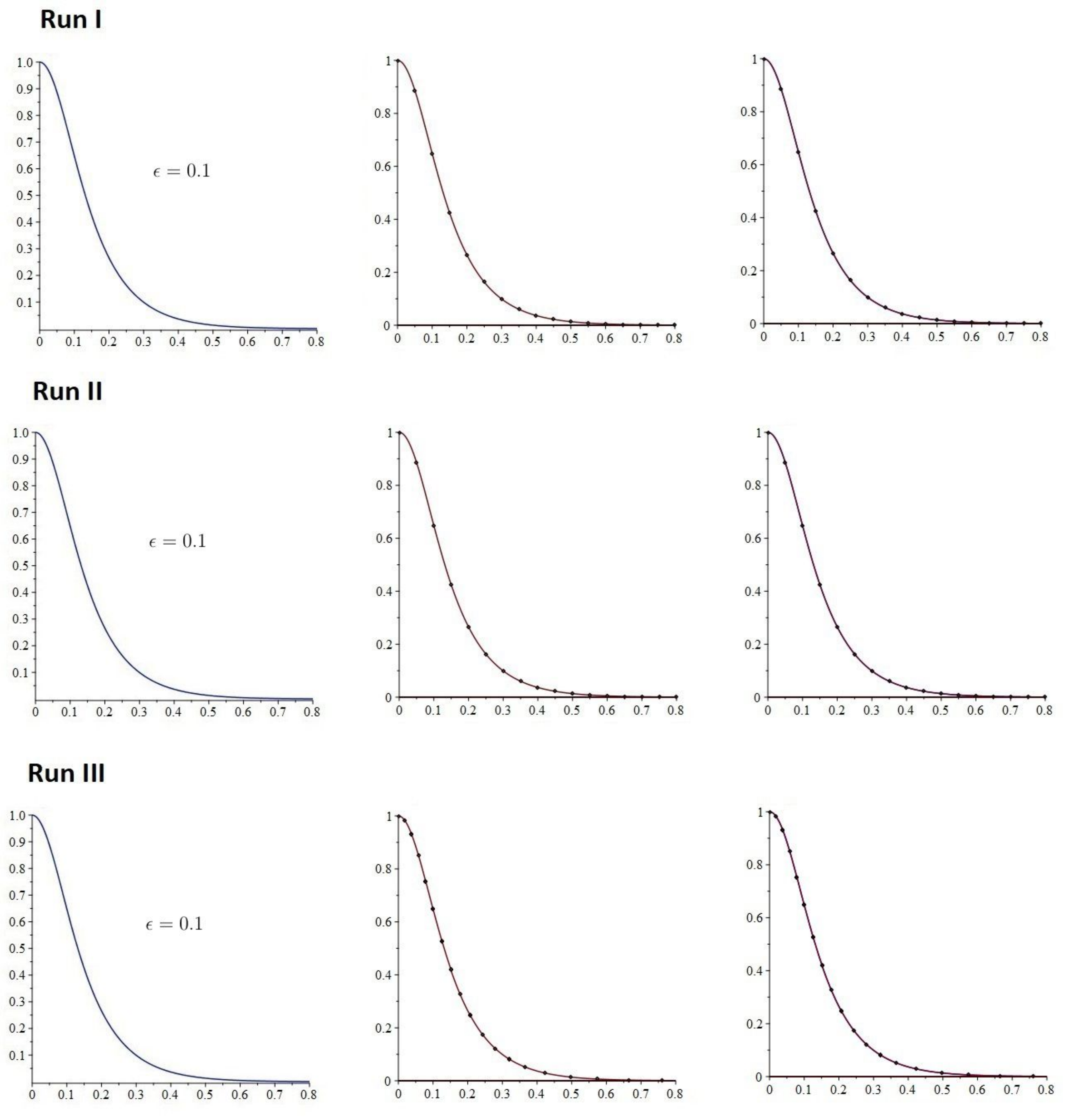}
\caption{Exact solution, B-spline approximation and the two superimposed for Run I (Gaussian points), Run II (equally spaced points), and Run III (nonlinear points), with $\epsilon = 0.1$, box $[0, 1]$, 20 intervals, 6 collocation sites per interval.}\label{fig:EPSILON01plots}
\end{center}
\end{figure}

\begin{figure}
\begin{center}
\includegraphics[width=6in]{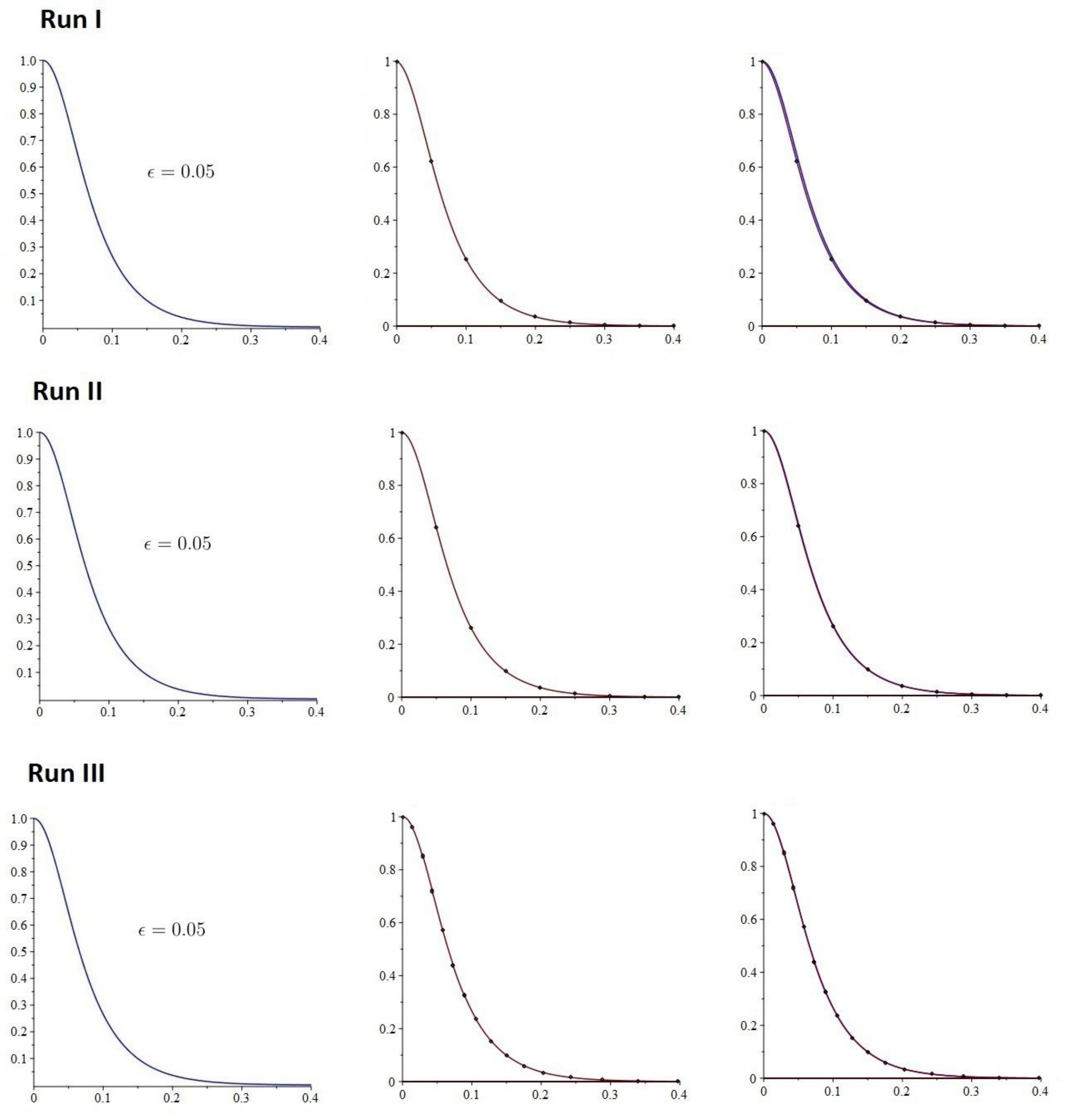}
\caption{Exact solution, B-spline approximation and the two superimposed for Run I (Gaussian points), Run II (equally spaced points), and Run III (nonlinear points), with $\epsilon = 0.05$, box $[0, 1]$, 20 intervals, 6 collocation sites per interval.}\label{fig:EPSILON005plots}
\end{center}
\end{figure}

\begin{figure}
\begin{center}
\includegraphics[width=6in]{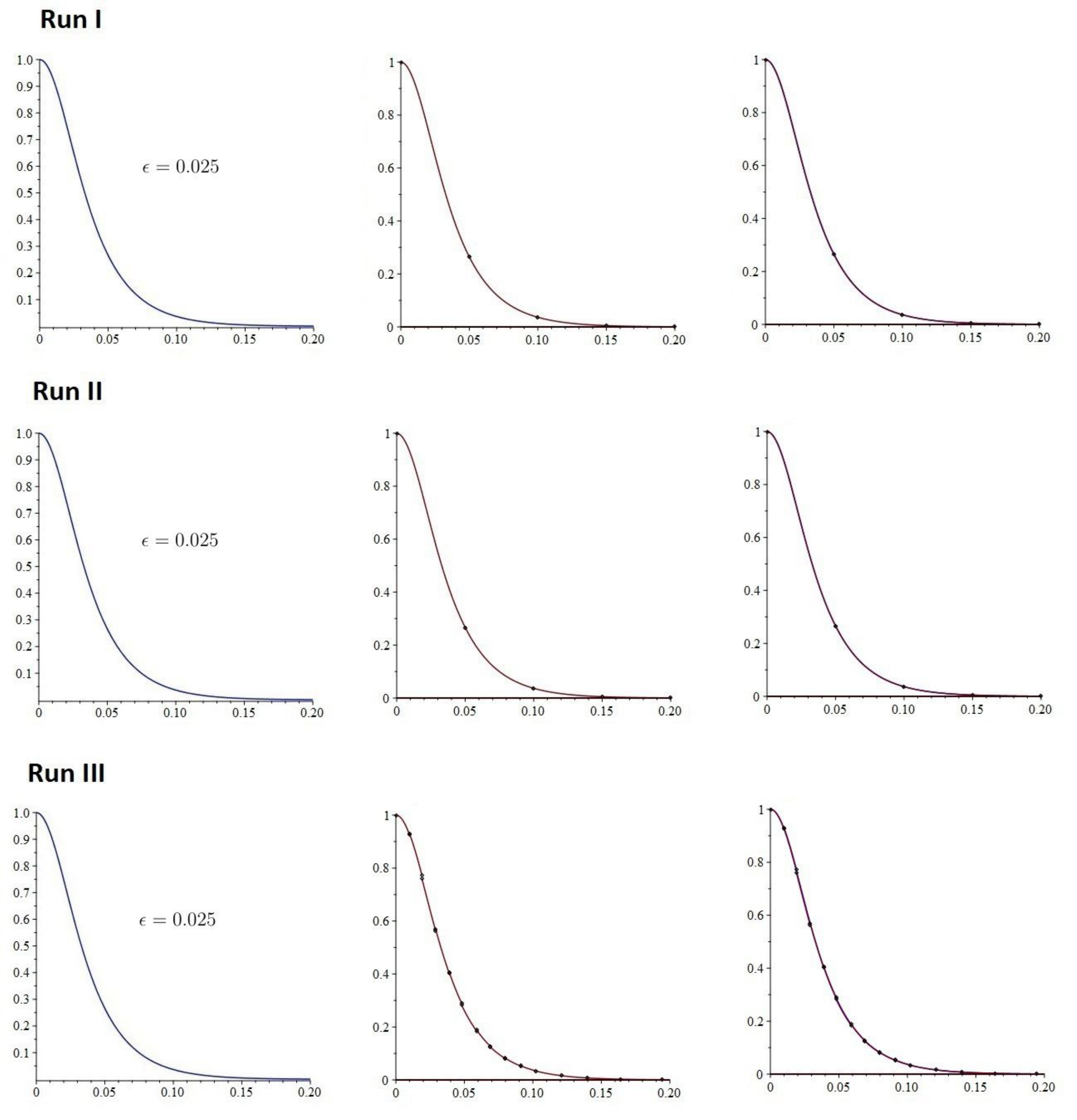}
\caption{Exact solution, B-spline approximation and the two superimposed for Run I (Gaussian points), Run II (equally spaced points), and Run III (nonlinear points), with $\epsilon = 0.025$, box $[0, 1]$, 20 intervals, 6 collocation sites per interval.}\label{fig:EPSILON0025plots}
\end{center}
\end{figure}

For each value of $\epsilon$, we conducted three runs: Run I using Gaussian collocation points; Run II using equally spaced collocation points; Run III using nonlinearly spaced collocation points (produced by the NEWNOT procedure). Approximation errors at selected points were recorded for each of these runs. These are displayed for  $\epsilon = 0.1$ and $\epsilon = 0.05$ in Figure~\ref{fig:EPSILON01005errors}, and for $\epsilon = 0.025$ in Figure~\ref{fig:EPSILON0025errors}. Corresponding plots of the exact solution, the B-spline approximation and the two superimposed are shown in Figure~\ref{fig:EPSILON01plots}, Figure~\ref{fig:EPSILON005plots} and Figure~\ref{fig:EPSILON0025plots}. 


Figures 5.3 to 5.5 show that all the approximations here are visually almost indistinguishable from the corresponding exact solutions. Note that in order to show inaccuracies as clearly as possible, all the plots shown in this chapter are `zoomed in' to the point where the exact solution converges to zero. This point moves closer and closer to the origin as $\epsilon \rightarrow 0$, as was shown in Figure~\ref{fig:nonlinear}. 

The approximation errors in Figure~\ref{fig:EPSILON01005errors} show that collocation at Gaussian points (Run I) produced slightly more accurate approximations than collocation at equally or nonlinearly spaced points for $\epsilon = 0.1$, while collocation using nonlinearly distributed collocation points (Run III) produced signficantly more accurate approximations than the other two configurations for $\epsilon = 0.05$. In the case of $\epsilon = 0.025$, equally and nonlinearly spaced points (Runs II and III) seem to perform approximately as well as each other, and both seem marginally better than collocation at Gaussian points. 

Therefore, the picture that emerges in this section is that we are able to obtain relatively good approximations to the exact solutions of the cubic Schr\"{o}dinger equations with perturbation parameters $\epsilon = 0.1$, $\epsilon = 0.05$ and $\epsilon = 0.025$, and there does not seem to be too much to choose between the three patterns of collocation points in Runs I, II and III in terms of there being one which is consistently better than the others. 

\section{Results for $\epsilon = 0.01$, $\epsilon = 0.005$ and $\epsilon = 0.001$} 

The accuracy of our approximations begins to deteriorate rapidly as we continue to reduce the size of the perturbation parameter. Here we again seek to approximate the exact solution (2.36) of the cubic Schr\"{o}dinger equation (2.34) with box $[0, 1]$ and 20 heptic polynomial pieces, but this time with the much smaller perturbation parameter values $\epsilon = 0.01$, $\epsilon = 0.005$ and $\epsilon = 0.001$. 

\begin{figure}
\begin{center}
\includegraphics[width=6in]{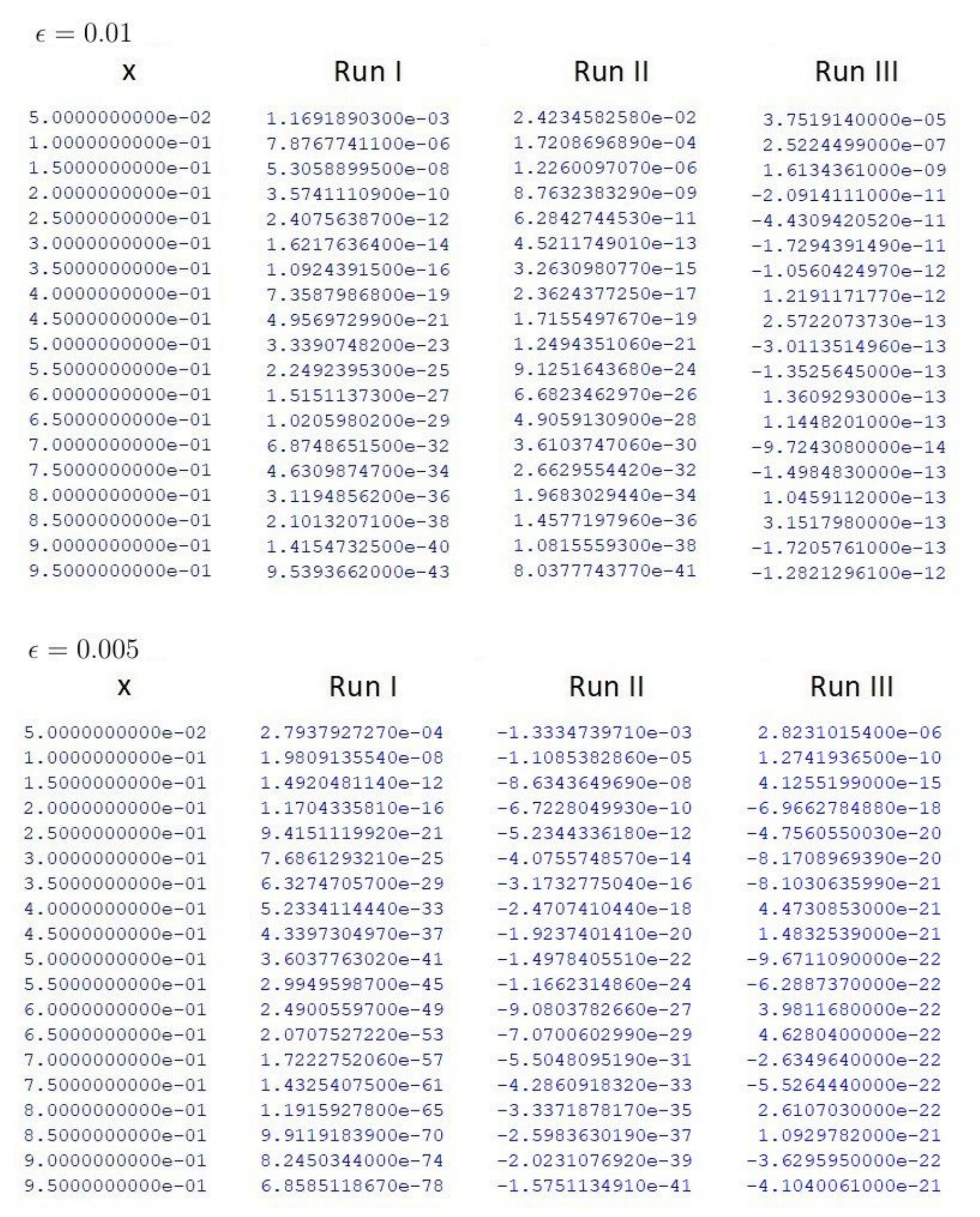}
\caption{Approximation errors at selected points $x$ for Run I (Gaussian points), Run II (equally spaced points), and Run III (nonlinear points), with box $[0, 1]$.}\label{fig:EPSILON0010005errors}
\end{center}
\end{figure}

\begin{figure}
\begin{center}
\includegraphics[width=6in]{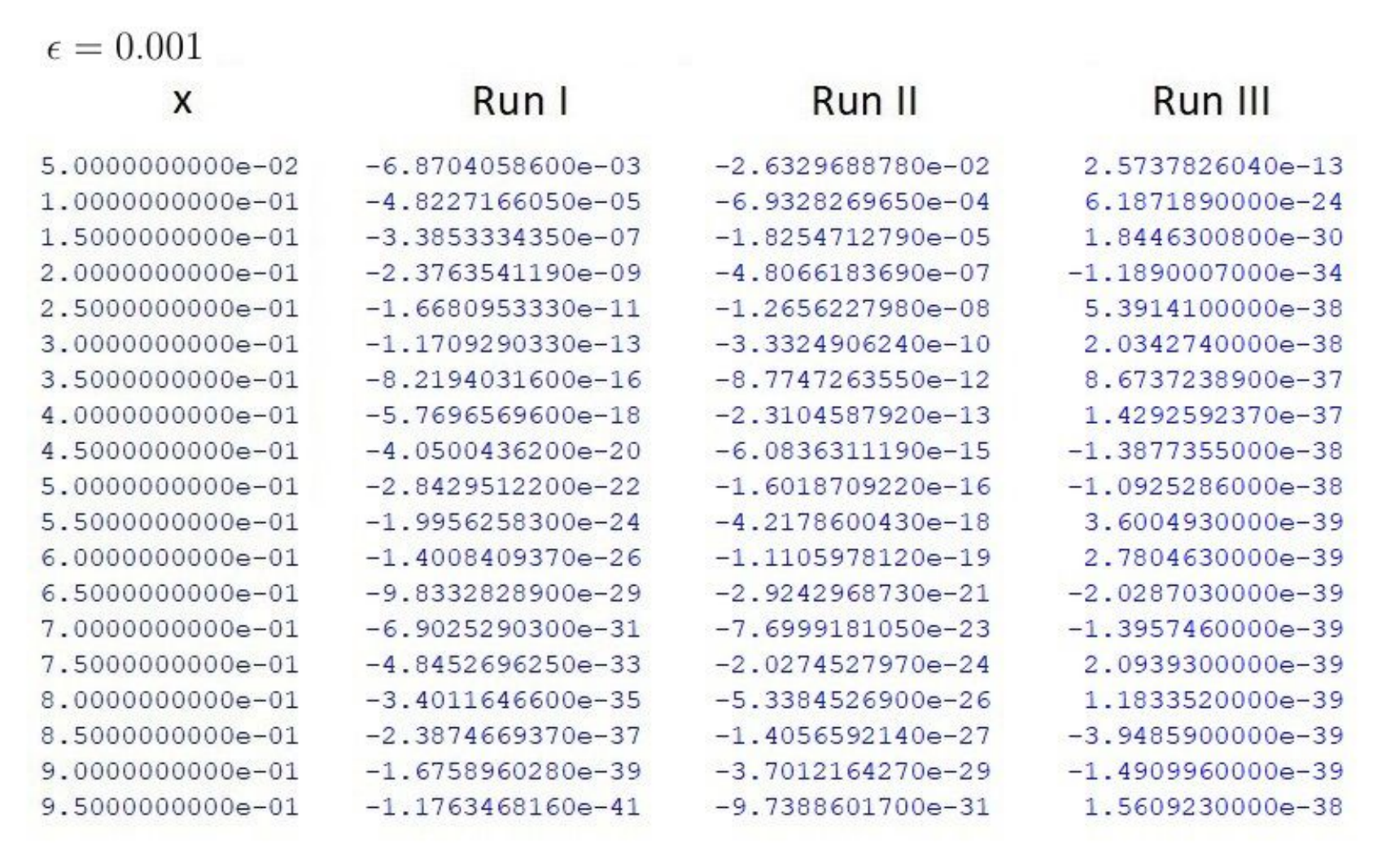}
\caption{Approximation errors at selected points $x$ for Run I (Gaussian points), Run II (equally spaced points), and Run III (nonlinear points), with box $[0, 1]$.}\label{fig:EPSILON0001errors}
\end{center}
\end{figure}

\begin{figure}
\begin{center}
\includegraphics[width=6in]{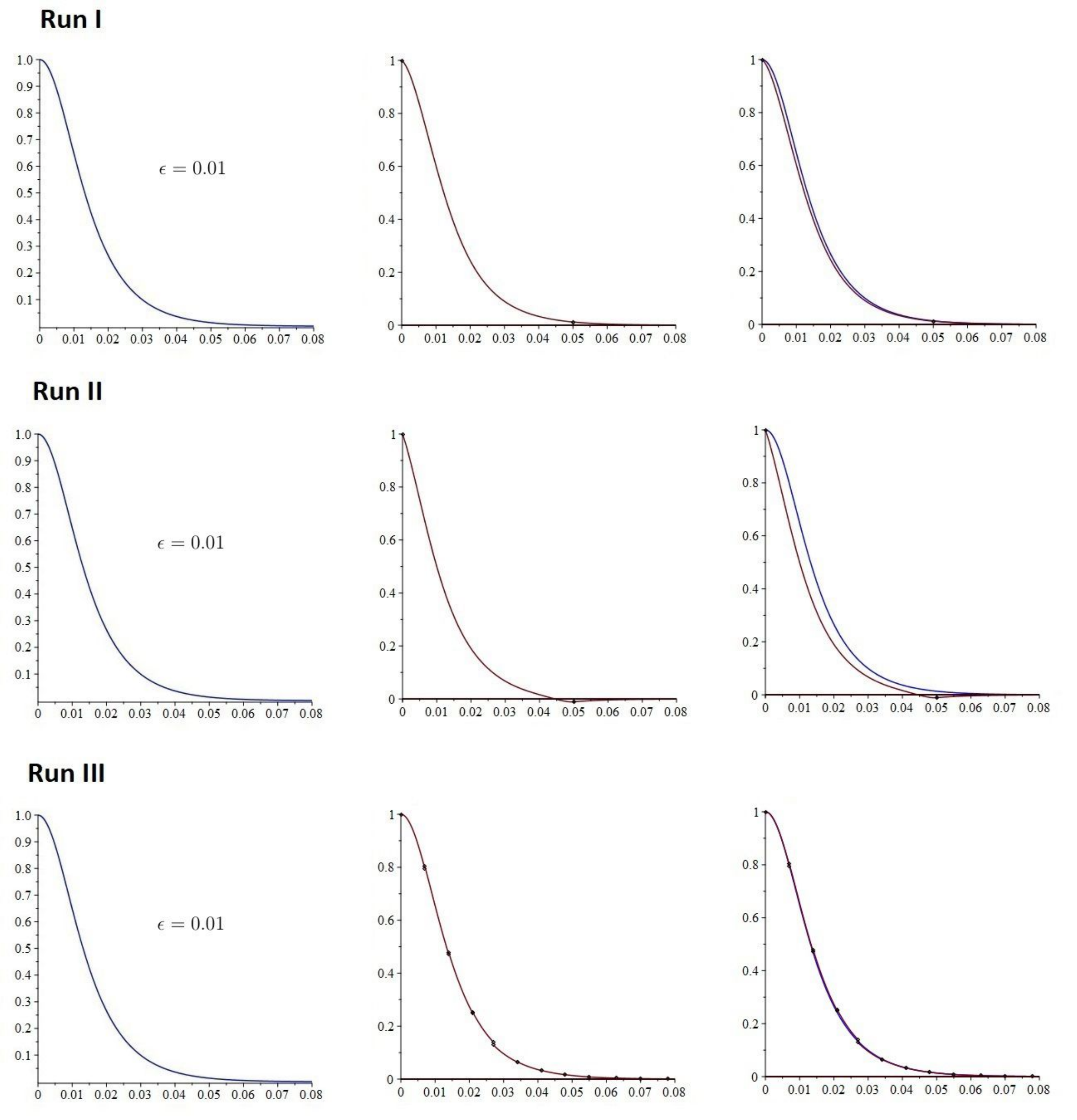}
\caption{Exact solution, B-spline approximation and the two superimposed for Run I (Gaussian points), Run II (equally spaced points), and Run III (nonlinear points), with $\epsilon = 0.01$, box $[0, 1]$, 20 intervals, 6 collocation sites per interval.}\label{fig:EPSILON001plots}
\end{center}
\end{figure}

\begin{figure}
\begin{center}
\includegraphics[width=6in]{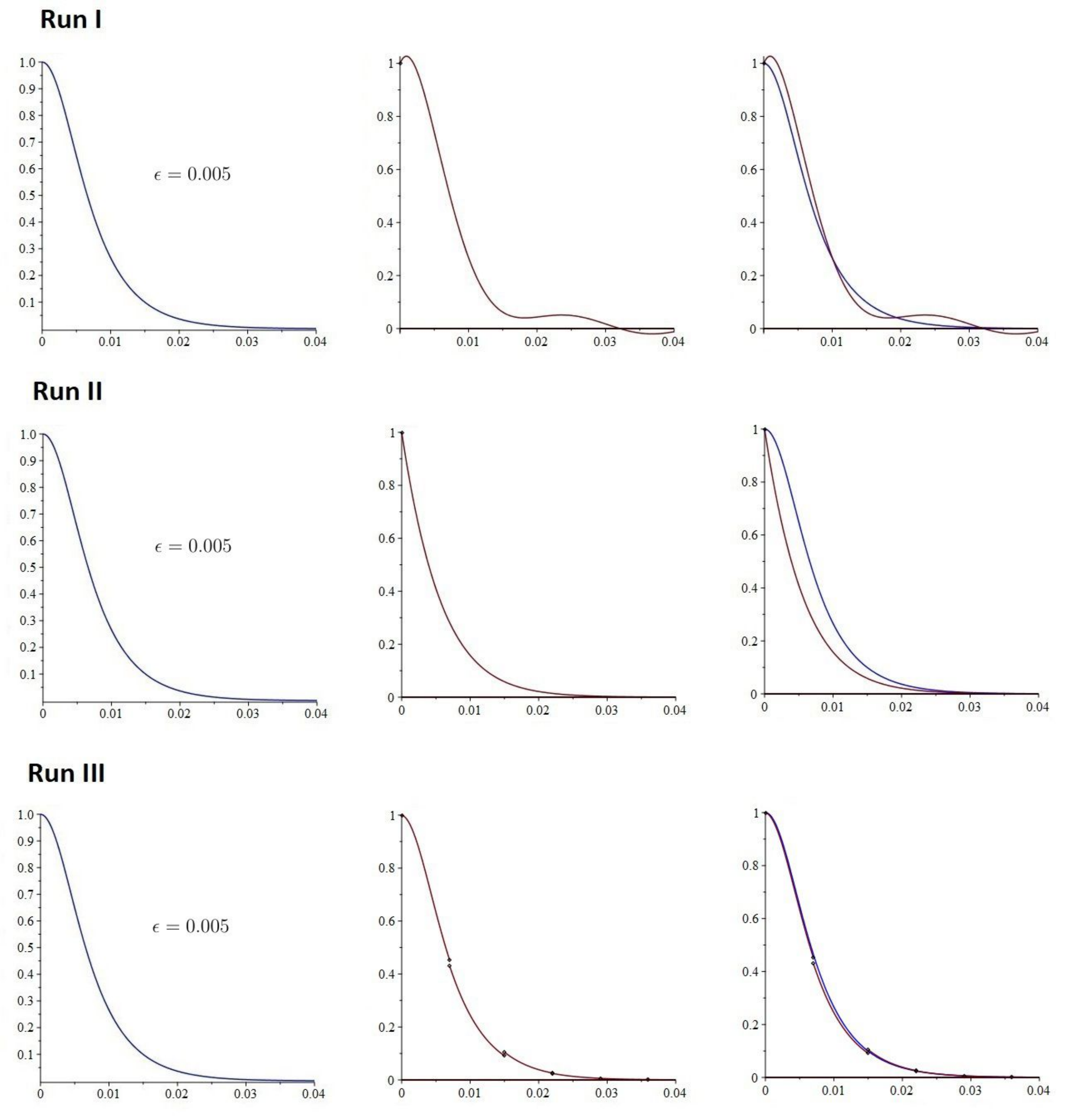}
\caption{Exact solution, B-spline approximation and the two superimposed for Run I (Gaussian points), Run II (equally spaced points), and Run III (nonlinear points), with $\epsilon = 0.005$, box $[0, 1]$, 20 intervals, 6 collocation sites per interval.}\label{fig:EPSILON0005plots}
\end{center}
\end{figure}

\begin{figure}
\begin{center}
\includegraphics[width=6in]{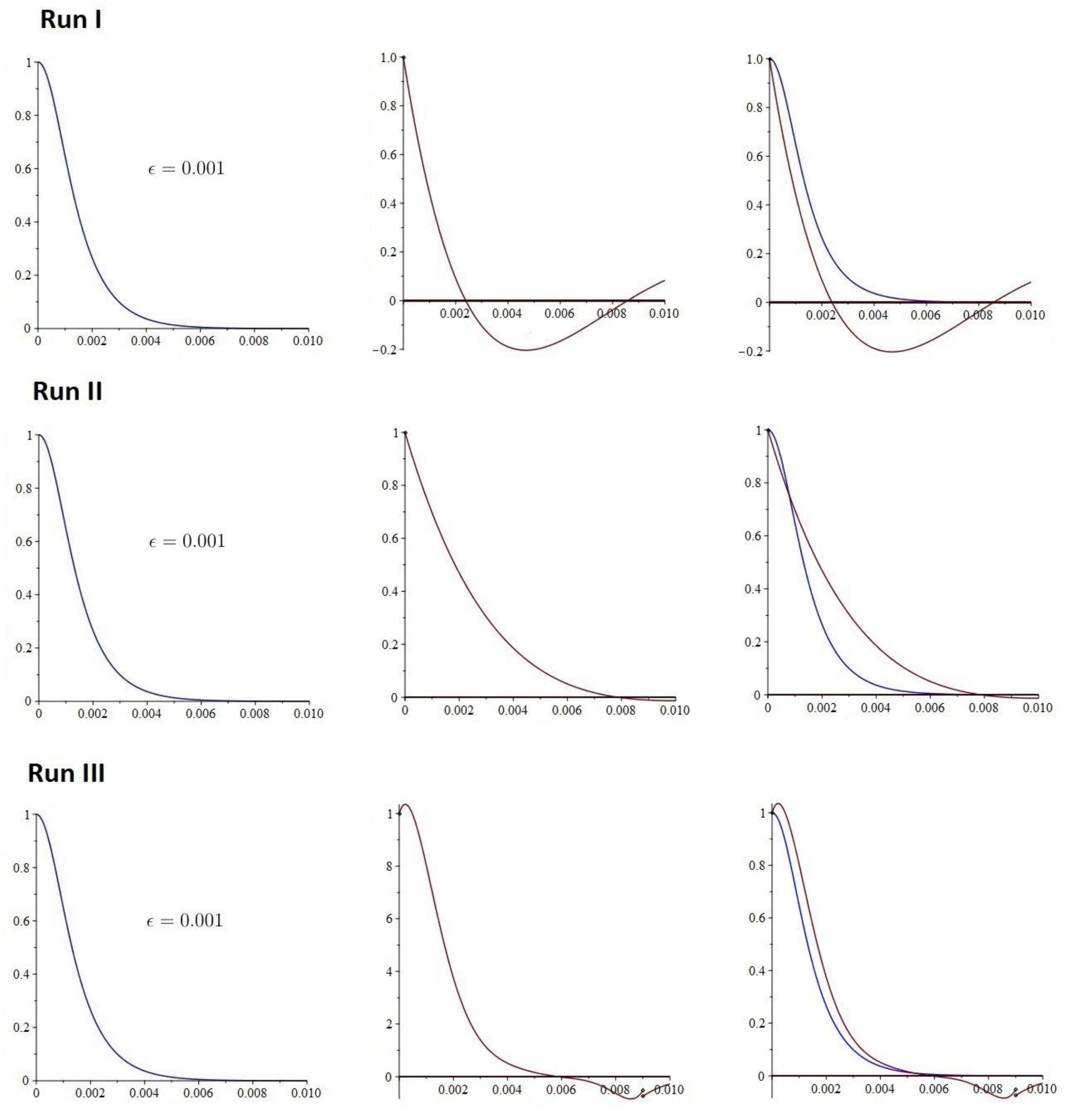}
\caption{Exact solution, B-spline approximation and the two superimposed for Run I (Gaussian points), Run II (equally spaced points), and Run III (nonlinear points), with $\epsilon = 0.001$, box $[0, 1]$, 20 intervals, 6 collocation sites per interval.}\label{fig:EPSILON0001plots}
\end{center}
\end{figure}

Approximation errors at selected points for Run I using Gaussian collocation points, Run II using equally spaced collocation points and Run III using nonlinearly spaced collocation points are displayed for $\epsilon = 0.01$ and $\epsilon = 0.005$ in Figure~\ref{fig:EPSILON0010005errors}, and for $\epsilon = 0.001$ in Figure~\ref{fig:EPSILON0001errors}. Corresponding plots of the exact solution, the B-spline approximation and the two superimposed are shown in Figure~\ref{fig:EPSILON001plots}, Figure~\ref{fig:EPSILON0005plots} and Figure~\ref{fig:EPSILON0001plots}. 

It is immediately apparent from Figures 5.8 to 5.10 that there are now quite serious divergences between our approximations and the corresponding exact solutions. There is also now a clear difference in performance between collcation at Gaussian points and equally spaced points on the one hand (Runs I and II) and collocation at nonlinearly distributed points produced by NEWNOT (Run III) on the other, with only the latter remaining visually close to the corresponding exact solutions in the cases $\epsilon = 0.01$ and $\epsilon = 0.005$. As can be seen in  Figure~\ref{fig:EPSILON0001plots}, collocation at Gaussian points seems to fail catastrophically in the case $\epsilon = 0.001$, with equally spaced collocation points also performing very poorly. Nonlinearly distributed points produced by NEWNOT perform better than the other two configurations in this case, though divergence between the approximation and the exact solution is now clearly visible even with this approach. 

The superiority of nonlinearly distributed collocation points over the other two configurations is also apparent from the tables of approximation errors in Figure~\ref{fig:EPSILON0010005errors} and Figure~\ref{fig:EPSILON0001errors}. Near the `boundary layer' at $x = 0$ in particular, the approximation errors for nonlinearly distributed collocation points are orders of magnitude smaller than for the other two configurations, presumably because NEWNOT is able to concentrate more of the collocation sites around this region where they are needed most. 

\begin{figure}
\begin{center}
\includegraphics[width=6in]{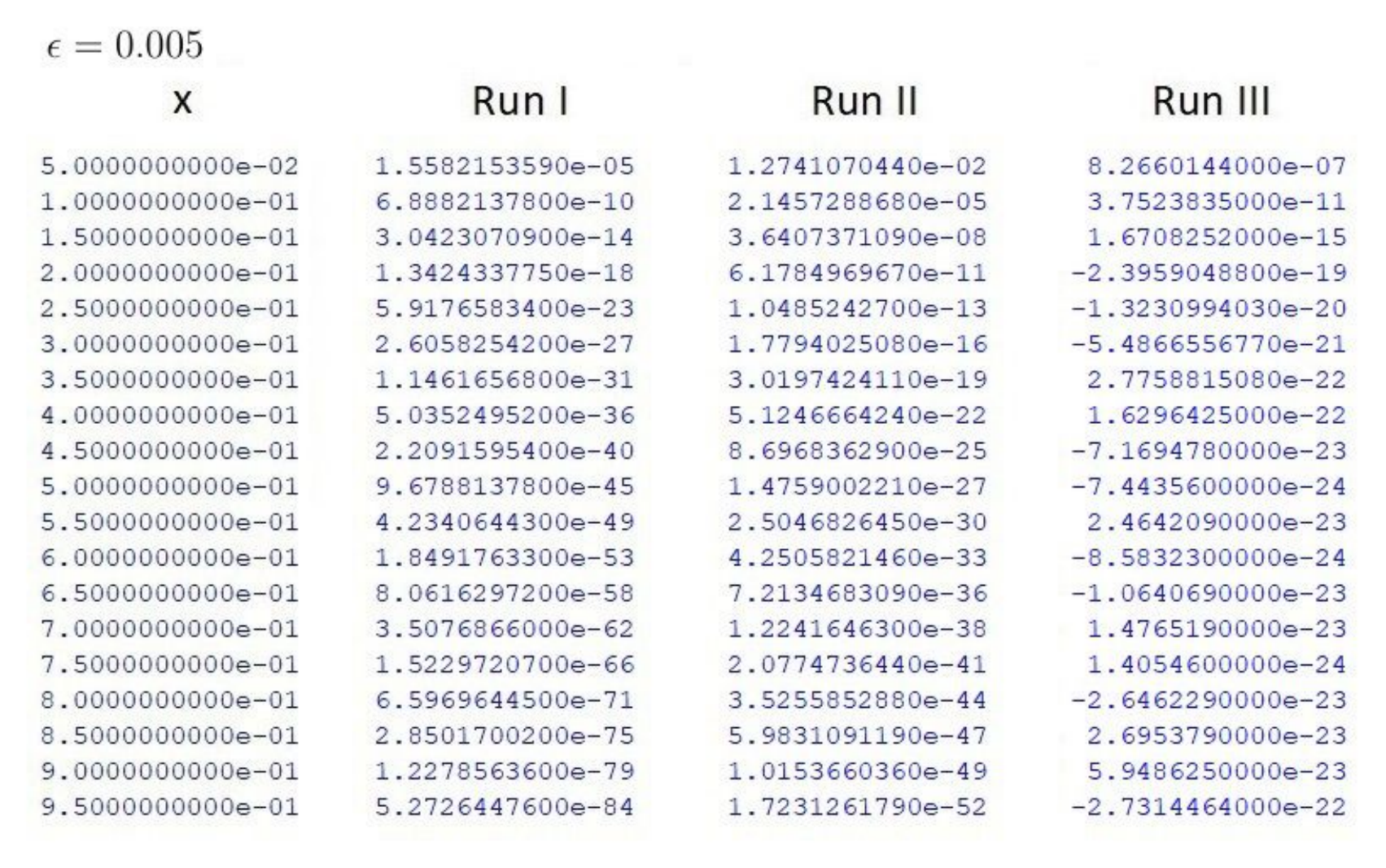}
\caption{Approximation errors at selected points $x$ for Run I (Gaussian points), Run II (equally spaced points), and Run III (nonlinear points), with box $[0, 1]$, 20 intervals, 8 collocation sites per interval.}\label{fig:EPSILON0005K8errors}
\end{center}
\end{figure}

\begin{figure}
\begin{center}
\includegraphics[width=6in]{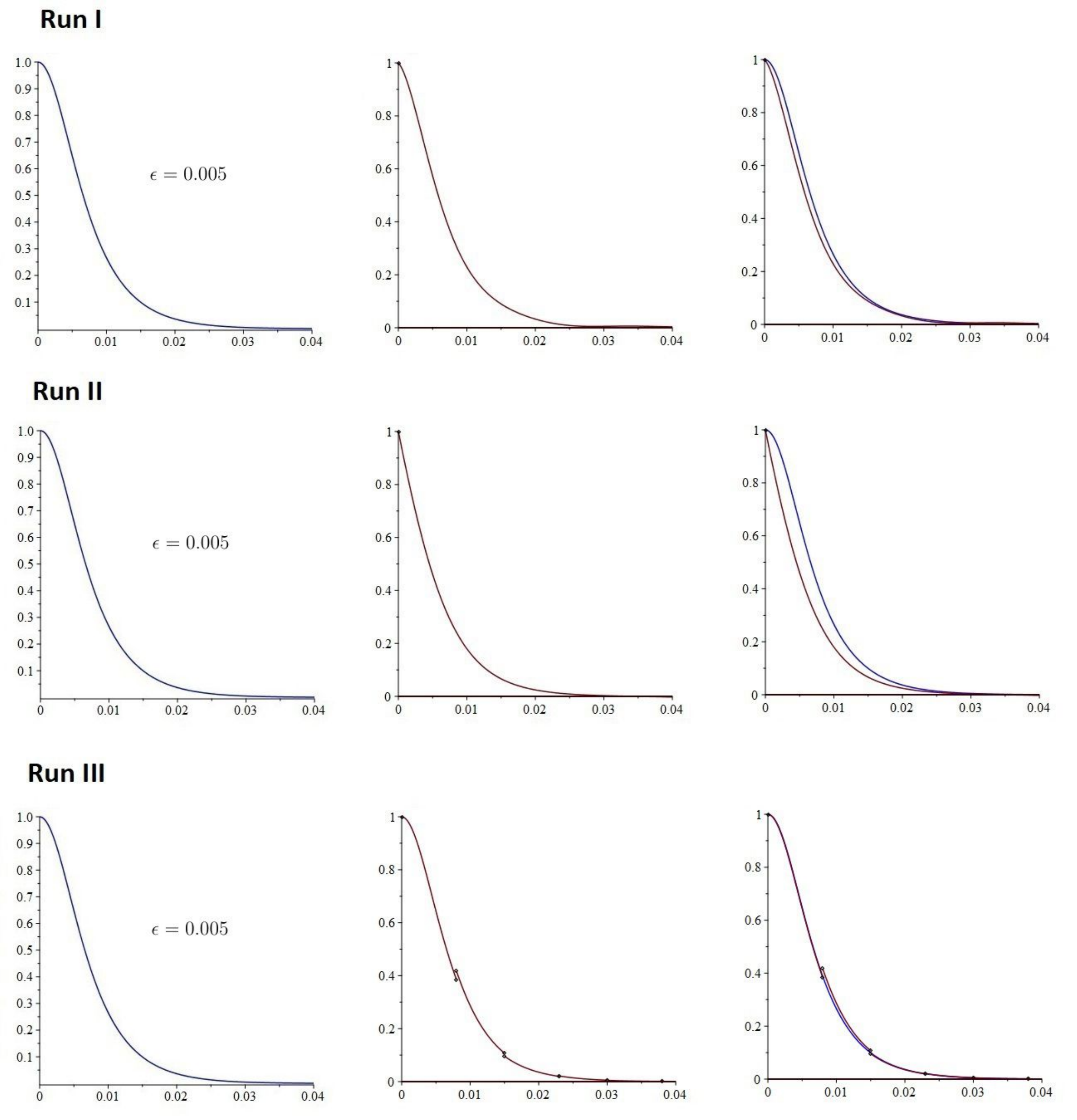}
\caption{Exact solution, B-spline approximation and the two superimposed for Run I (Gaussian points), Run II (equally spaced points), and Run III (nonlinear points), with $\epsilon = 0.005$, box $[0, 1]$, 20 intervals, 8 collocation sites per interval.}\label{fig:EPSILON0005K8plots}
\end{center}
\end{figure}

To see if any improvements could be made to the approximations in this section, we also experimented with higher numbers of polynomial pieces and higher numbers of collocation sites per subinterval. Only moderate improvements were possible, as shown by sample results for an experiment with $\epsilon = 0.005$, box $[0, 1]$, 20 subintervals and 8 collocation sites per subinterval (i.e., nonic polynomial pieces), reported in Figure~\ref{fig:EPSILON0005K8errors} and Figure~\ref{fig:EPSILON0005K8plots}.   



\chapter{Conclusions}

Treating the radial equation in Shore's paper as a two-point BVP rather than as a regular Sturm-Liouville problem, and thereby focusing on approximating its eigenfunctions rather than its eigenvalues, we have thoroughly investigated the relative performance of equally spaced collocation points, Gaussian collocation points and nonlinearly distributed collocation points in approximating Schr\"{o}dinger wave functions for the hydrogen atom. We were able to expand our exploration by extending the framework in Shore's original paper to include radial equations with nonzero angular momentum, using novel transformations of these equations to enable de Boor's methodology to be applied to them. We also succeeded in extending the basic framework in Shore's study to a nonlinear Schr\"{o}dinger equation with cubic nonlinearity, enabling us to explore the relative performance of the three different patterns of collocation sites in this setting as well. These investigations have yielded numerous insights not only into the relative performance of Gaussian collocation points, but also into the numerical effects of changing box sizes, meshes, and orders of polynomial approximants in conjunction with the different patterns of collocation sites, as well as into the overall applicability and limitations of de Boor's B-spline collocation methodology in the case of Schr\"{o}dinger's equation.

With regard to the electron wave functions for the hydrogen atom, a clear and consistent result is that equally spaced collocation points perform less well than either Gaussian points or nonlinearly distributed points. Equally spaced collocation points are sometimes used in the atomic theory literature so this result is of relevance in assessing the suitability of this approach. It is also clear that Gaussian points can be successfully applied in the hydrogen atom context. Our results confirm that there are combinations of box sizes, mesh sizes and orders of polynomial approximants for which Gaussian points yield better results than either of the other two configurations. We did encounter some situations, typically in which the mesh was relatively coarse for the given box size, when nonlinearly distributed collocation points performed better than Gaussian points. Otherwise, the performance of Gaussian points was either better or more or less on a par with that of nonlinearly distributed points. One might therefore have expected Gaussian collocation points to appear more often in the atomic theory literature.  

We found the situation to be different in the case of the perturbed nonlinear Schr\"{o}dinger equation, which is actually a boundary layer problem of the type exemplified in Chapter XV of de Boor's book. As the size of the perturbation parameter was reduced in our numerical experiments, nonlinearly distributed collocation points produced by the NEWNOT subroutine began to significantly outperform both equally spaced and Gaussian collocation points, eventually by orders of magnitude. This is, perhaps, not too surprising as the example in Chapter XV of de Boor's book and COLLOC's ability to call on NEWNOT seem to have been tailored to cater for the kind of boundary layer problem which we encountered with the cubic Schr\"{o}dinger equation.     

On the basis of our numerical results overall, it seems likely that Gaussian collocation points can perform at least as well as nonlinearly distributed points, and possibly better, in situations where the Schr\"{o}dinger wave functions being approximated do not exhibit excessively sudden oscillations or changes in curvature, and where the mesh and number of collocation sites per subinterval are adequate for the box size. Mostly, these favourable conditions seemed to be the prevailing ones in the case of the hydrogen atom. In less favourable situations, nonlinearly distributed collocation points might outperform Gaussian points due to the greater flexibility in being able to concentrate the collocation sites in difficult regions, thereby improving the quality of the approximation there. This clearly became a significant advantage in the case of the nonlinear Schr\"{o}dinger equation.     

With regard to the effects of changing the box size, it was surprising to find that in some situations an increase in box size led to a worsening of approximation accuracy, probably because the mesh then became too coarse relative to the larger interval. In the cases of equally spaced and Gaussian collocation points, it will not have been possible to re-distribute collocation sites to compensate for this effect, so these approaches tended to perform less well than nonlinearly distributed points in these situations. The emphasis in the atomic theory literature is almost always on ensuring that the box size is not too small. Our results show that it is also necessary to ensure that the box size does not become too large relative to the mesh being used.

Not too surprisingly, we found that finer meshes and larger numbers of collocation sites per subinterval produced greater approximation accuracy. In our experiments we did not find that either one of these was particularly more effective than the other in improving accuracy. On the contrary, we found that there was not much to choose between them in this respect. It did come as a surprise, however, that with the relatively large box sizes required for the excited states of the electron in the hydrogen atom, it was not possible to increase both the number of subintervals and the number of collcation sites per subinterval \emph{together} to a greater extent. In exploring the limits of this, we found that it was not possible in some cases to have a combination of more than forty subintervals with six or more collocation sites per subinterval, as this led to matrix sizes for the collocation equations that were larger than those accommodated by de Boor's package of programs. This was an unexpected limitation.      

Another interesting issue is that de Boor's collocation methodology, as exemplified in Chapter XV of his book, is unable to produce nontrivial results when the column vector on the right-hand side of the matrix system (3.8) is a zero vector. For the purposes of our numerical experiments using different patterns of collocation sites, we had to rely on our pre-existing knowledge of the eigenvalues and exact solutions of Schr\"{o}dinger's radial equation to be able to implement the equations as two-point BVPs with a nonzero vector on the right-hand side of (3.8). We were then able to focus on the numerical performance of different patterns of collocation sites in approximating the \emph{eigenfunctions} of Schr\"{o}dinger's equation. This produced visually and numerically rich outputs which enabled more detailed assessments of numerical performance to be made than if we had focused on estimating individual \emph{eigenvalues}, as Shore did in his paper. However, if our objective had been to solve for \emph{both} eigenvalues \emph{and} eigenfunctions in Schr\"{o}dinger's equation as if they were both unknown, we would not have been able to employ the two-point BVP approach in Chapter XV of de Boor's book. This distinction between our approach and Shore's approach became much clearer as a result of the detailed study of de Boor's methodology for the purposes of this dissertation.

There is scope for extending our study in a number of interesting directions. We have only focused on time-independent Schr\"{o}dinger equations in this dissertation. It is possible to use collocation approaches with the full time-dependent Schr\"{o}dinger equation as well, and indeed this is explored using Shore's methodology in \cite{odero}. It would be interesting to see if our two-point BVP approach using de Boor's methodology could be extended to time-dependent Schr\"{o}dinger equations. Another avenue for extending our approach is to consider two-dimensional problems, for example, the helium atom. The application of B-splines to this and other many-body problems is discussed in \cite{sapirstein}, and again there is scope for exploring how de Boor's methodology could be applied here. Our numerical experiments in this dissertation have involved only negative energy systems. Ideally we would have liked to explore the applicability of our methods to positive energy scenarios as well, i.e., scattering problems. Shore successfully applied his approach to scattering from an Eckart potential in \cite{shore}, focusing on obtaining estimates of reflection and transmission probabilities. It would be an interesting and challenging exercise to see if de Boor's approach could be applied to approximating the wave functions for scattering problems, as these are generally complex-valued with both real and imaginary components. Finally, there are many other areas of physics and nonlinear science in which there do not seem to have been any applications of B-spline methods so far. For example, there do not appear to be any applications of B-splines in the context of general relativity.         


\appendix

\chapter{Derivation of electron wave function in hydrogen}\label{first-appendix}

In this note I try to provide a thorough derivation of the electron's wave function in the hydrogen atom, bringing out the mathematical details clearly. The exposition is guided by a number of texts including \cite{beiser}, \cite{brehm}, \cite{griffiths}, \cite{jain}, and \cite{boas}.
    
In general, four quantum numbers are needed to fully describe atomic electrons in many-electron atoms. These four numbers and their permissible values are:

Principal quantum number $\tilde{n} = 1, 2, 3, \ldots$

Orbital quantum number $l = 0, 1, 2, \ldots, (\tilde{n} - 1)$

Magnetic quantum number $m_l = 0, \pm 1, \pm 2, \ldots, \pm l$

Spin magnetic quantum number $m_s = -\frac{1}{2}, +\frac{1}{2}$

The principal quantum number determines the electron's energy, the orbital quantum number its orbital angular-momentum magnitude, the magnetic quantum number its orbital angular-momentum direction, and the spin magnetic quantum number its spin direction.

I have noticed that it is often not explained clearly why, for example, the orbital quantum number cannot exceed the principal quantum number minus one, or why the magnitude of the magnetic quantum number cannot exceed that of the orbital quantum number. I want to bring out details like this clearly. The time-independent Schrödinger equation for the hydrogen atom only involves the first three quantum numbers. I will not discuss the spin magnetic quantum number here.

\section{Schrödinger's wave equation for the electron in the hydrogen atom}
In Cartesian coordinates, Schrödinger's three-dimensional equation for the electron in the hydrogen atom is
\begin{equation*}
\frac{\partial^2 \psi}{\partial x^2} + \frac{\partial^2 \psi}{\partial y^2} + \frac{\partial^2 \psi}{\partial z^2} + \frac{2m_e}{\hbar^2}(E - U) \psi = 0
\end{equation*}
where $m_e$ denotes the electron mass. The potential energy $U$ is the electric potential energy of a charge $-e$ given that it is at distance $r$ from another charge $+e$, namely
\begin{equation*}
U = -\frac{e^2}{4 \pi \epsilon_0 r}
\end{equation*}
It is necessary to change variables in Schrödinger's equation since the potential energy is a function of radial distance $r$ rather than the Cartesian coordinate variables $x$, $y$ and $z$. Given the spherical symmetry of the atom, it is sensible to proceed by changing the variables in Schrödinger's equation to those of spherical polar coordinates (rather than changing the $r$ variable in $U$ to Cartesian coordinates using $r = \sqrt{x^2 + y^2 + z^2}$). Only the variables in the Laplacian part of Schrödinger's equation need to be changed, so we can use a standard approach to changing variables in Laplace's equation (see \cite{boas}, p. 228) to get
\begin{equation*}
\frac{1}{r^2} \frac{\partial }{\partial r}\bigg( r^2 \frac{\partial \psi}{\partial r}\bigg) + \frac{1}{r^2 \sin \theta}\frac{\partial }{\partial \theta}\bigg( \sin \theta \frac{\partial \psi}{\partial \theta} \bigg) + \frac{1}{r^2 \sin^2 \theta}\frac{\partial^2 \psi}{\partial \phi^2} + \frac{2m_e}{\hbar^2}(E - U) \psi = 0
\end{equation*}
I will now temporarily simplify things by using the representation of the square of the angular momentum operator in spherical polar coordinates (see \cite{jain}, p. 207), namely
\begin{equation*}
L^2 = - \hbar^2 \bigg( \frac{1}{\sin \theta}\frac{\partial}{\partial \theta} \bigg( \sin \theta \frac{\partial}{\partial \theta}\bigg) + \frac{1}{\sin^2 \theta}\frac{\partial^2}{\partial \phi^2} \bigg)
\end{equation*} 
\begin{equation*}
= - \hbar^2 r^2 \bigg( \frac{1}{r^2 \sin \theta}\frac{\partial}{\partial \theta} \bigg( \sin \theta \frac{\partial}{\partial \theta}\bigg) + \frac{1}{r^2 \sin^2 \theta}\frac{\partial^2}{\partial \phi^2} \bigg)
\end{equation*}
Using this to replace the two middle terms in Schrödinger's equation and rearranging we get
\begin{equation*}
\frac{1}{r^2} \frac{\partial }{\partial r}\bigg( r^2 \frac{\partial \psi}{\partial r}\bigg) + \frac{2m_e}{\hbar^2}(E - U) \psi = \frac{L^2}{\hbar^2 r^2} \psi
\end{equation*}
This equation can now be solved by the usual separation of variables approach. We assume that the $\psi$ function can be expressed as a product
\begin{equation*}
\psi(r, \theta, \phi) = R(r) Y(\theta, \phi)
\end{equation*}
and then substitute this back into the wave equation to get
\begin{equation*}
\frac{Y}{r^2} \frac{d}{ d r}\bigg( r^2 \frac{d R}{d r}\bigg) + \frac{2m_e}{\hbar^2}(E - U) R Y = \frac{R}{\hbar^2 r^2} L^2 Y
\end{equation*}
Dividing through by $\frac{R Y}{r^2}$ we get
\begin{equation*}
\frac{1}{R} \frac{d}{d r}\bigg( r^2 \frac{d R}{d r}\bigg) + \frac{2m_e r^2}{\hbar^2}(E - U) = \frac{1}{\hbar^2 Y} L^2 Y
\end{equation*}
Since the left-hand side of this equation depends only on $r$ while the right-hand side depends only on $\theta$ and $\phi$, both sides must be equal to some constant which we can call $\lambda$. Setting the left and right-hand sides equal to $\lambda$ in turn and rearranging slightly we finally get the radial equation
\begin{equation*}
\frac{1}{r^2} \frac{d}{d r}\bigg( r^2 \frac{d R}{d r}\bigg) + \bigg[ \frac{2m_e}{\hbar^2}(E - U) - \frac{\lambda}{r^2} \bigg] R = 0
\end{equation*}
and the angular equation
\begin{equation*}
L^2 Y = \lambda \hbar^2 Y
\end{equation*}
We can now apply separation of variables again to the angular equation. Rewriting the operator $L^2$ in full the angular equation becomes
\begin{equation*}
- \hbar^2 \bigg( \frac{1}{\sin \theta}\frac{\partial}{\partial \theta} \bigg( \sin \theta \frac{\partial Y}{\partial \theta}\bigg) + \frac{1}{\sin^2 \theta}\frac{\partial^2 Y}{\partial \phi^2} \bigg) = \lambda \hbar^2 Y
\end{equation*}
which simplifies to
\begin{equation*}
\frac{1}{\sin \theta}\frac{\partial}{\partial \theta} \bigg( \sin \theta \frac{\partial Y}{\partial \theta}\bigg) + \frac{1}{\sin^2 \theta}\frac{\partial^2 Y}{\partial \phi^2} + \lambda Y = 0
\end{equation*}
We assume that the $Y$ function can be written as the product
\begin{equation*}
Y(\theta, \phi) = \Theta (\theta) \Phi (\phi)
\end{equation*}
Substituting this into the angular equation gives
\begin{equation*}
\frac{\Phi}{\sin \theta}\frac{d}{d \theta} \bigg( \sin \theta \frac{d \Theta}{d \theta}\bigg) + \frac{\Theta}{\sin^2 \theta}\frac{d^2 \Phi}{d \phi^2} + \lambda Y \Theta \Phi = 0
\end{equation*}
Multiplying through by $\frac{\sin^2 \theta}{\Theta \Phi}$ and rearranging we get
\begin{equation*}
-\frac{1}{\Phi} \frac{d^2 \Phi}{d \phi^2} = \frac{\sin^2 \theta}{\Theta}\bigg[ \frac{1}{\sin \theta}\frac{d}{d \theta}\bigg(\sin \theta \frac{d \Theta}{d \theta}\bigg) + \lambda \Theta\bigg]
\end{equation*}
Since the left-hand side of this equation depends only on $\phi$ while the right-hand side depends only on $\theta$, both sides must be equal to some constant which we can provisionally call $k$. Setting the left and right-hand sides equal to $k$ in turn and rearranging we get
\begin{equation*}
\frac{d^2 \Phi}{d \phi^2} + k \Phi = 0
\end{equation*}
and
\begin{equation*}
\frac{1}{\sin \theta}\frac{d}{d \theta}\bigg(\sin \theta \frac{d \Theta}{d \theta}\bigg) + \bigg( \lambda - \frac{k}{\sin^2 \theta}\bigg) \Theta = 0
\end{equation*} 
We now have three ordinary differential equations that need to be solved, one for $\Phi$, one for $\Theta$ and one for $R$. We will solve each of them in turn.
\section{The equation for $\Phi$}
The equation for $\Phi$ is a straightforward second-order differential equation with auxiliary equation
\begin{equation*}
\zeta^2 + k = 0
\end{equation*}
implying $\zeta = \pm \sqrt{-k}$ if $k < 0$ and $\zeta  = \pm i \sqrt{k}$ if $k > 0$. Therefore it has a general solution of the form
\begin{equation*}
\Phi(\phi) = Ae^{\sqrt{-k} \phi} + Be^{- \sqrt{-k} \phi}
\end{equation*}
if $k < 0$ and 
\begin{equation*}
\Phi(\phi) = Ae^{i \sqrt{k} \phi} + Be^{- i \sqrt{k} \phi} 
\end{equation*}
if $k > 0$, where $A$ and $B$ are arbitrary constants. Now, the azimuth angle $\phi$ can take any value in $(-\infty, \infty)$ but the function $\Phi$ must take a single value at each point in space (since this is a required property of the quantum wave function which $\Phi$ is a constituent of). It follows that the function $\Phi$ must be periodic since it must take the same value at $\phi$ and $\phi + 2\pi$ for any given $\phi$. This imposes two constraints on the form of the general solution: (1) it cannot consist only of exponential functions with real arguments since these are not periodic (thus ruling out the first general solution above and thereby implying that the separation constant $k$ must be nonnegative); (2) $\sqrt{k}$ must be an integer. Given these constraints, it is customary in quantum mechanics to denote $\pm \sqrt{k}$ by the letter $m$ (it is called the magnetic quantum number) and to specify the separation constant in the angular equations as $m^2$, which guarantees its nonnegativity. We then state the general solution of the equation for $\Phi$ as
\begin{equation*}
\Phi(\phi) = Ae^{i m \phi} + Be^{- i m \phi}
\end{equation*}
In principle this allows two independent solutions but we only need one of them for any given electron wave function. Either of the two terms in the general solution would satisfy $\Phi^{\prime \prime} = -m^2 \Phi$, so we will take only the first as is customary. We can therefore state the general solution of the equation for $\Phi$ for a given magnetic quantum number $m$ as
\begin{equation*}
\Phi(\phi) \propto e^{i m \phi}
\end{equation*}

\section{The equation for $\Theta$}
Given that we now know the separation constant for the angular equations is either zero or a positive square number $k = m^2$, we can write the equation for $\Theta$ as
\begin{equation*}
\frac{1}{\sin \theta}\frac{d}{d \theta}\bigg(\sin \theta \frac{d \Theta}{d \theta}\bigg) + \bigg( \lambda - \frac{m^2}{\sin^2 \theta}\bigg) \Theta = 0
\end{equation*}
Expanding the first term we get
\begin{equation*}
\frac{1}{\sin \theta} \cos \theta \frac{d \Theta}{d \theta} + \frac{d^2 \Theta}{d \theta^2} + \bigg( \lambda - \frac{m^2}{\sin^2 \theta}\bigg) \Theta = 0
\end{equation*}
I am now going to multiply and divide the first two terms by $ \sin^2 \theta$ to get
\begin{equation*}
\sin^2 \theta \bigg(\frac{\cos \theta}{\sin^3 \theta} \frac{d \Theta}{d \theta} + \frac{1}{\sin^2 \theta} \frac{d^2 \Theta}{d \theta^2} \bigg) + \bigg( \lambda - \frac{m^2}{\sin^2 \theta}\bigg) \Theta = 0
\end{equation*}
$\iff$
\begin{equation*}
\sin^2 \theta \bigg(- \frac{\cos \theta}{\sin^3 \theta} \frac{d \Theta}{d \theta} + \frac{1}{\sin^2 \theta} \frac{d^2 \Theta}{d \theta^2} + \frac{2 \cos \theta}{\sin^3 \theta} \frac{d \Theta}{d \theta}\bigg) + \bigg( \lambda - \frac{m^2}{\sin^2 \theta}\bigg) \Theta = 0
\end{equation*}
Now we can make the change of variable $x = \cos \theta$ which implies $dx = - \sin \theta d \theta$ and therefore
\begin{equation*}
\frac{d \theta}{dx} = - \frac{1}{\sin \theta}
\end{equation*}
\begin{equation*}
\frac{d \Theta}{d x} = \frac{d \Theta}{d \theta} \frac{d \theta}{d x} = - \frac{1}{\sin \theta} \frac{d \Theta}{d \theta}
\end{equation*}
\begin{equation*}
\frac{d^2 \Theta}{d x^2} = \frac{d}{d \theta} \bigg[ - \frac{1}{\sin \theta} \frac{d \Theta}{d \theta} \bigg] \frac{d \theta}{d x} = - \frac{\cos \theta}{\sin^3 \theta} \frac{d \Theta}{d \theta} + \frac{d^2 \Theta}{d \theta^2}
\end{equation*}
Using these in the amended form of the $\Theta$ equation together with the fact that $\sin^2 \theta$ = $ 1 - x^2$, the $\Theta$ equation becomes
\begin{equation*}
(1 - x^2) \bigg(\frac{d^2 \Theta}{d x^2} - \frac{2x}{1 - x^2} \frac{d \Theta}{d x} \bigg) + \bigg(\lambda - \frac{m^2}{1 - x^2} \bigg) \Theta = 0
\end{equation*}
$\iff$
\begin{equation*}
(1 - x^2) \frac{d^2 \Theta}{d x^2} - 2x \frac{d \Theta}{d x} + \bigg(\lambda - \frac{m^2}{1 - x^2} \bigg) \Theta = 0
\end{equation*}
We will solve this equation first for the case $m = 0$ (the solutions will be Legendre polynomials) and use these results to construct solutions for the case $m \neq 0$ (the solutions here will be the associated Legendre functions). Setting $m = 0$ we get
\begin{equation*}
(1 - x^2) \frac{d^2 \Theta}{d x^2} - 2x \frac{d \Theta}{d x} + \lambda \Theta = 0
\end{equation*}
which has the form of a well known differential equation known as Legendre's equation. It can be solved by assuming a series solution of the form
\begin{equation*}
\Theta = a_0 + a_1 x + a_2 x^2 + a_3 x^3 + a_4 x^4 + \cdots + a_n x^n + \cdots
\end{equation*}
and then differentiating it term by term twice to get
\begin{equation*}
\Theta^{\prime} = a_1  + 2 a_2 x + 3 a_3 x^2 + 4 a_4 x^3 + \cdots + n a_n x^{n-1} + \cdots
\end{equation*}
and
\begin{equation*}
\Theta^{\prime \prime} = 2 a_2 + 6 a_3 x + 12 a_4 x^2 + 20 a_5 x^3 + \cdots + n (n - 1) a_n x^{n-2} + \cdots
\end{equation*}
We now substitute these into Legendre's equation and set the coefficient of each power of $x$ equal to zero (because $\Theta$ must satisfy Legendre's equation identically). We find that the coefficient of the $x^n$ term satisfies
\begin{equation*}
(n + 2)(n + 1) a_{n + 2} + (\lambda - n(n + 1)) a_n = 0
\end{equation*}
which implies
\begin{equation*}
a_{n + 2} = - \frac{(\lambda - n(n + 1))}{(n + 2)(n + 1)} a_n
\end{equation*}
This formula makes it possible to find any even coefficient as a multiple of $a_0$ and any odd coefficient as a multiple of $a_1$. The general solution of our Legendre equation is then a sum of two series involving two arbitrary constants $a_0$ and $ a_1$:
\begin{equation*}
\Theta = a_0 \bigg \{1 - \frac{\lambda}{2!} x^2 + \frac{\lambda (\lambda - 6)}{4!} x^4 - \frac{\lambda(\lambda - 6)(\lambda - 20)}{6!} x^6 + \cdots \bigg \}
\end{equation*}
\begin{equation*}
+ a_1 \bigg \{x - \frac{(\lambda - 2)}{3!} x^3 + \frac{(\lambda - 2)(\lambda - 12)}{5!} x^5 - \frac{(\lambda - 2)(\lambda - 12)(\lambda - 30)}{7!} x^7 + \cdots \bigg \}
\end{equation*}
Both of the series in this sum converge for $ x^2 < 1$ but in general they do not converge for $x^2 = 1$. This is a problem for us because in our change of variables we set $x = \cos \theta$ and we want solutions that converge for all possible values of $\theta$ including those that result in $x^2 = 1$. It turns out that the only way to get such solutions is to choose integer values of $\lambda$ that make either the $a_0$ or the $a_1$ series in the above sum terminate (the other series will generally be divergent so we remove it by setting the corresponding arbitrary constant equal to zero). This requires $\lambda$ to take values in the quadratic sequence $0$, $2$, $6$, $12$, $20$, $30$, $42$, $56 \ldots$ The $l$-th term of this sequence is $l(l + 1)$, so the separation constant $\lambda$ must be of this form, i.e., $\lambda = l(l + 1)$ for some $l = 0, 1, 2, 3, \dots$. When $l$ takes an even value the $a_0$ series will terminate and we can set $a_1 = 0$ to make the other series vanish. Conversely, when $l$ takes an odd value the $a_1$ series will terminate and we can set $a_0 = 0$ to make the other series vanish.

From the eigenvalue equation for $L^2$ given earlier ($L^2 Y = \lambda \hbar^2 Y$) it is clear that the magnitude of the orbital angular momentum is $L = \sqrt{l(l + 1)} \hbar$. It is interesting to see how the form of this arises mathematically from considering series solutions to Legendre's equation above. The parameter $l$ is called the orbital angular momentum quantum number.

Note that negative integral values of $l$ are allowed but they simply give solutions already obtained for positive values. For example, $l = -2$ gives $\lambda = 2$ and this makes the $a_1$ series terminate, yielding the polynomial solution
\begin{equation*}
\Theta = a_1 x
\end{equation*}
This is exactly the same solution as the one that would be obtained if $l = 1$. It is therefore customary to restrict $l$ to nonnegative values. Each possible value of $l$ gives a polynomial solution to Legendre's equation. For $l = 0$ we get $\Theta = a_0$, for $l = 1$ we get $\Theta = a_1 x$, for $l = 2$ we get $\Theta = a_0 - 3 a_0 x^2$, and so on. If the value of $a_0$ or $a_1$ in each polynomial equation is selected so that $\Theta = 1$ when $x = 1$ the resulting polynomials are called Legendre polynomials, denoted by $P_l(x)$. Given that for each $l$ we have $P_l(1) = 1$ the first few Legendre polynomials are
\begin{equation*}
P_0(x) = 1
\end{equation*}
\begin{equation*}
P_1(x) = x
\end{equation*}
\begin{equation*}
P_2(x) = \frac{1}{2}(3 x^2 - 1)
\end{equation*}
\begin{equation*}
P_3(x) = \frac{1}{2}\big(5 x^3 - 3x \big)
\end{equation*}
These are the physically acceptable solutions to Legendre's equation for $\Theta$ above.

We now consider the solutions for $m \neq 0$ of the equation
\begin{equation*}
(1 - x^2) \frac{d^2 \Theta}{d x^2} - 2x \frac{d \Theta}{d x} + \bigg(\lambda - \frac{m^2}{1 - x^2} \bigg) \Theta = 0
\end{equation*}
We now know that $\lambda = l(l + 1)$ so we can write this in and we can also add the subscript $l$ to $m$ as the solutions to this equation will involve a link the between the values of the orbital angular momentum and magnetic quantum numbers. The equation we need to solve becomes
\begin{equation*}
(1 - x^2) \frac{d^2 \Theta}{d x^2} - 2x \frac{d \Theta}{d x} + \bigg[l(l + 1) - \frac{m_l^2}{1 - x^2} \bigg] \Theta = 0
\end{equation*}
The link between $l$ and $m_l$ arises from the fact that we are constrained in trying to solve this equation: it encompasses the case $m_l = 0$ for which the physically acceptable solutions are the Legendre polynomials $P_l(x)$. Therefore the physically allowable solutions for the above equation must include the Legendre polynomials as a special case. We can find these by using the series approach again and it turns out that the physically acceptable solutions are the so-called associated Legendre functions which take the form
\begin{equation*}
P_l^{m_l}(x) = (1 - x^2)^{m_l/2}\frac{d^{m_l}}{d x^{m_l}}P_l(x)
\end{equation*}
Now, each Legendre polynomial $P_l(x)$ is a polynomial of degree $l$. Therefore the $m_l$-th order derivative in $P_l^{m_l}$ will equal zero if $|m_l| > l$, so for physically acceptable solutions we must impose the constraint $|m_l| \leq l$ in the differential equation for $\Theta$. This is where the link between the quantum numbers $l$ and $m_l$ comes from in the quantum theory of the hydrogen atom: given a value of $l$ the acceptable values of $m_l$ are integers in the range $-l \leq m_l \leq l$.

Finally, note two things: (1) $P_l^{m_l}(x)$ reduces to the Legendre polynomial $P_l(x)$ when $m_l = 0$, which is what we needed. (2) A negative value for $m_l$ does not change $m_l^2$ in the original differential equation so a solution for positive $m_l$ is also a solution for the corresponding negative $m_l$. Thus many references define the associated Legendre function $P_l^{m_l}(x)$ for $-l \leq m_l \leq l$ as $P_l^{|m_l|}(x)$.

To conclude, given values for the quantum numbers $l$ and $m_l$, the general solution of the equation for $\Theta$ can be written as
\begin{equation*}
\Theta(\theta) \propto P_l^{m_l}(\cos \theta)
\end{equation*}

\section{The radial equation for R}
To clarify where the principal quantum number comes from, the final equation we need to deal with is the radial equation
\begin{equation*}
\frac{1}{r^2} \frac{d}{d r}\bigg( r^2 \frac{d R}{d r}\bigg) + \bigg[ \frac{2m_e}{\hbar^2}(E - U) - \frac{\lambda}{r^2} \bigg] R = 0
\end{equation*}
Writing $\lambda = l(l + 1)$ and replacing $U$ with the formula for the potential energy we get
\begin{equation*}
\frac{1}{r^2} \frac{d}{d r}\bigg( r^2 \frac{d R}{d r}\bigg) + \bigg[ \frac{2m_e}{\hbar^2}\bigg(\frac{e^2}{4 \pi \epsilon_0 r} + E\bigg) - \frac{l(l + 1)}{r^2} \bigg] R = 0
\end{equation*}
$\iff$
\begin{equation*}
\frac{d^2 R}{d r^2} + \frac{2}{r} \frac{d R}{d r} + \frac{2m_e}{\hbar^2} \bigg[ E + \frac{e^2}{4 \pi \epsilon_0 r} - \frac{l(l + 1) \hbar^2}{2 m_e r^2} \bigg] R = 0
\end{equation*}
We are only interested in solutions for which the electron is bound within the atom, so we take $E < 0$ (the negative energy of the electron is the amount of energy that must be supplied to it to free it from the atom). In order to solve the above equation it is customary to make the change of variable
\begin{equation*}
\rho = \bigg(-\frac{8 m_e E}{\hbar^2}\bigg)^{1/2} r
\end{equation*}
and define the dimensionless constant
\begin{equation*}
\tau = \frac{e^2}{4 \pi \epsilon_0 \hbar} \bigg(-\frac{m_e}{2 E} \bigg)^{1/2}
\end{equation*}
If we then specify $ R = R(\rho)$ we have
\begin{equation*}
\frac{d R}{d r} = \frac{d R}{d \rho} \frac{d \rho}{d r} = \bigg(-\frac{8 m_e E}{\hbar^2}\bigg)^{1/2} \frac{d R}{d \rho}
\end{equation*}
\begin{equation*}
\frac{2}{r} \frac{d R}{d r} = \bigg(-\frac{8 m_e E}{\hbar^2}\bigg) \frac{2}{\rho} \frac{d R}{d \rho}
\end{equation*}
\begin{equation*}
\frac{d^2 R}{d r^2} = \bigg(-\frac{8 m_e E}{\hbar^2}\bigg)^{1/2} \frac{d^2 R}{d \rho^2} \frac{d \rho}{d r} = \bigg(-\frac{8 m_e E}{\hbar^2}\bigg) \frac{d^2 R}{d \rho^2}
\end{equation*}
\begin{equation*}
\frac{2 m_e}{\hbar^2} E + \frac{2 m_e}{\hbar^2} \frac{e^2}{4 \pi \epsilon_0 r} = \bigg(-\frac{8 m_e E}{\hbar^2}\bigg) \bigg \{\frac{1}{4}\frac{e^2}{4 \pi \epsilon_0 r} \bigg(-\frac{1}{E}\bigg) - \frac{1}{4}\bigg \} = \bigg(-\frac{8 m_e E}{\hbar^2}\bigg) \bigg(\frac{\tau}{\rho} - \frac{1}{4}\bigg)
\end{equation*}
\begin{equation*}
\frac{l(l + 1)}{r^2} = \bigg(-\frac{8 m_e E}{\hbar^2}\bigg) \frac{l(l + 1)}{\rho^2}
\end{equation*}
Using these results we can rewrite the differential equation as
\begin{equation*}
\bigg(-\frac{8 m_e E}{\hbar^2}\bigg) \bigg \{\frac{d^2 R}{d \rho^2} + \frac{2}{\rho} \frac{d R}{d \rho} + \bigg[\frac{\tau}{\rho} - \frac{1}{4} - \frac{l(l + 1)}{\rho^2}\bigg] R(\rho) \bigg \} = 0
\end{equation*}
$\iff$
\begin{equation*}
\frac{d^2 R}{d \rho^2} + \frac{2}{\rho} \frac{d R}{d \rho} + \bigg[\frac{\tau}{\rho} - \frac{1}{4} - \frac{l(l + 1)}{\rho^2}\bigg] R(\rho) = 0
\end{equation*}
To make further progress we consider the behaviour of this differential equation as $\rho \rightarrow \infty$. It reduces to
\begin{equation*}
\frac{d^2 R}{d \rho^2} - \frac{1}{4} R = 0
\end{equation*}
which is a straightforward second-order differential equation with auxiliary equation
\begin{equation*}
\zeta^2 - \frac{1}{4} = 0
\end{equation*}
\begin{equation*}
\implies \zeta = \pm \frac{1}{2}
\end{equation*}
The positive solution to the auxiliary equation implies a term in the general solution of the form $e^{\rho/2}$ which is unacceptable since it explodes as $\rho \rightarrow \infty$. Therefore we only accept the negative solution to the auxiliary equation and the general solution for $R$ as $\rho \rightarrow \infty$ must be of the form
\begin{equation*}
R \propto e^{-\rho/2}
\end{equation*}
This suggests we can try an exact solution of the full differential equation of the form
\begin{equation*}
R = e^{-\rho/2} F(\rho)
\end{equation*}
Differentiating this twice we get
\begin{equation*}
\frac{d R}{d \rho} = -\frac{1}{2} e^{-\rho/2} F(\rho) + e^{-\rho/2} F^{\prime} (\rho)
\end{equation*}
\begin{equation*}
\frac{d^2 R}{d \rho^2} = \frac{1}{4} e^{-\rho/2} F(\rho) - \frac{1}{2} e^{-\rho/2} F^{\prime}(\rho) - \frac{1}{2} e^{-\rho/2} F^{\prime}(\rho) + e^{-\rho/2} F^{\prime \prime} (\rho)
\end{equation*}
Substituting these into the differential equation
\begin{equation*}
\frac{d^2 R}{d \rho^2} + \frac{2}{\rho} \frac{d R}{d \rho} + \bigg[\frac{\tau}{\rho} - \frac{1}{4} - \frac{l(l + 1)}{\rho^2}\bigg] R(\rho) = 0
\end{equation*}
gives
\begin{equation*}
F^{\prime \prime}(\rho) + \frac{(2 - \rho)}{\rho} F^{\prime}(\rho) + \bigg[\frac{(\tau - 1)}{\rho} - \frac{l(l + 1)}{\rho^2}\bigg] F(\rho) = 0
\end{equation*}
$\iff$
\begin{equation*}
\rho^2 F^{\prime \prime}(\rho) + \rho (2 - \rho) F^{\prime}(\rho) + \big[\rho (\tau - 1) - l(l + 1)\big] F(\rho) = 0
\end{equation*}
We can now try to solve this latest version of the differential equation by the method of Frobenius, which involves assuming a generalised power series solution of the form
\begin{equation*}
F(\rho) = a_0 \rho^s + a_1 \rho^{s + 1} + a_2 \rho^{s + 2} + \cdots
\end{equation*}
Differentiating twice we get
\begin{equation*}
F^{\prime}(\rho) = s a_0 \rho^{s-1} + (s + 1) a_1 \rho^s + (s + 2) a_2 \rho^{s + 1} + \cdots
\end{equation*}
\begin{equation*}
F^{\prime \prime}(\rho) = (s - 1) s a_0 \rho^{s-2} + s (s + 1) a_1 \rho^{s-1} + (s + 1) (s + 2) a_2 \rho^s + \cdots
\end{equation*}
Then the terms appearing in the differential equation have the generalised power series forms
\begin{equation*}
\rho^2 F^{\prime \prime}(\rho) = (s - 1) s a_0 \rho^s + s (s + 1) a_1 \rho^{s+1} + (s + 1) (s + 2) a_2 \rho^{s+2} + \cdots
\end{equation*}
\begin{equation*}
2 \rho F^{\prime}(\rho) = 2 s a_0 \rho^s + 2 (s + 1) a_1 \rho^{s+1} + 2 (s + 2) a_2 \rho^{s + 2} + \cdots
\end{equation*}
\begin{equation*}
-\rho^2 F^{\prime}(\rho) = - s a_0 \rho^{s+1} - (s + 1) a_1 \rho^{s+2} - (s + 2) a_2 \rho^{s + 3} - \cdots
\end{equation*}
\begin{equation*}
(\tau - 1) \rho F(\rho) = (\tau - 1) a_0 \rho^{s+1} + (\tau - 1) a_1 \rho^{s + 2} + (\tau - 1) a_2 \rho^{s + 3} + \cdots
\end{equation*}
\begin{equation*}
-l(l + 1) F(\rho) = -l(l + 1) a_0 \rho^s - l(l + 1) a_1 \rho^{s + 1} - l(l + 1) a_2 \rho^{s + 2} - \cdots
\end{equation*}
Summing these terms (remembering that the sum must be identically equal to zero) we find the coefficient of $\rho^s$ to be
\begin{equation*}
[s(s - 1) + 2s - l(l+1)]a_0 = 0
\end{equation*}
$\implies$
\begin{equation*}
s(s + 1) - l(l+1) = 0
\end{equation*}
$\implies s = l$ or $s = -l - 1$. Now, when $s = -l - 1$ the first term of the power series for $F(\rho)$ is $a_0/\rho^{l+1}$ which explodes as $\rho \rightarrow 0$. This is unacceptable so we discard this solution and set $s = l$.

For the coefficient of $\rho^{s + n}$ we get
\begin{equation*}
[(s+n)(s+n-1) + 2(s+n) - l(l+1)]a_n + [(\tau - 1) - (s+n-1)]a_{n-1} = 0
\end{equation*}
Setting $s = l$ and rearranging gives us the recurrence equation
\begin{equation*}
a_n = \frac{(l + n - \tau)}{(l+n+1)(l+n) - l(l+1)} a_{n-1}
\end{equation*}
From this recurrence equation we observe that
\begin{equation*}
a_n \rightarrow \frac{1}{n}a_{n-1} = \frac{1}{n!}a_0
\end{equation*}
as $n \rightarrow \infty$. We deduce from this that the series for $F(\rho)$ becomes like $a_0 \rho^l \sum \frac{\rho^n}{n!}$ as $n \rightarrow \infty$ and therefore $R = e^{-\rho/2} F(\rho)$ becomes like $a_0 \rho^l e^{\rho/2}$. However, this diverges as $\rho \rightarrow \infty$ which is unacceptable, so we conclude that the series for $F(p)$ must terminate at some value of $n$ which we will call $N$. In this case we have $a_{N+1} = 0$ which the recurrence equation tells us can only happen if
\begin{equation*}
\tau = l + N + 1 \equiv \tilde{n}
\end{equation*}
This is how the principal quantum number $\tilde{n}$ first appears. Now, we have
\begin{equation*}
\tau = \frac{e^2}{4 \pi \epsilon_0 \hbar} \bigg(-\frac{m_e}{2 E} \bigg)^{1/2} = \tilde{n}
\end{equation*}
$\iff$
\begin{equation*}
\bigg(\frac{e^2}{4 \pi \epsilon_0}\bigg)^2 \bigg(-\frac{m_e}{2 \hbar^2}\bigg) \frac{1}{E} = \tilde{n}^2
\end{equation*}
\begin{equation*}
\iff E_{\tilde{n}} = \bigg(-\frac{m_e}{2 \hbar^2}\bigg) \bigg(\frac{e^2}{4 \pi \epsilon_0}\bigg)^2 \frac{1}{\tilde{n}^2}
\end{equation*}
These are the famous bound-state energy eigenvalues for $\tilde{n} = 1, 2, \ldots$. This is the same formula for the energy levels of the hydrogen atom that Niels Bohr obtained by intuitive means in his 1913 \emph{solar system} model of atomic structure.

As stated above, the integer $\tilde{n}$ is called the principal quantum number. Recall that $\tilde{n} = l + N + 1$ and $N$ cannot be smaller than zero. It follows that
\begin{equation*}
\tilde{n} - l - 1 \geq 0
\end{equation*}
$\iff$
\begin{equation*}
l \leq \tilde{n} – 1
\end{equation*}
This explains why for a given value of $\tilde{n}$ the allowable values of $l$ are $l = 0, 1, 2, \dots, (\tilde{n}-1)$.

Returning to the solution of
\begin{equation*}
\rho^2 F^{\prime \prime}(\rho) + \rho (2 - \rho) F^{\prime}(\rho) + \big[\rho (\tau - 1) - l(l + 1)\big] F(\rho) = 0
\end{equation*}
the above discussion suggests that we should look for a solution of the form
\begin{equation*}
F(\rho) = a_0 \rho^l L(\rho)
\end{equation*}
where $L(\rho)$ is a polynomial (rather than an infinite series). Differentiating this twice gives
\begin{equation*}
F^{\prime}(\rho) = a_0 l \rho^{l-1}L(\rho) + a_0 \rho^l L^{\prime}(\rho)
\end{equation*}
\begin{equation*}
F^{\prime \prime} (\rho) = a_0 (l-1) l \rho^{l-2}L(\rho) + 2 a_0 l \rho^{l-1} L^{\prime}(\rho) + a_0 \rho^l L^{\prime \prime}(\rho)
\end{equation*}
Substituting these into the differential equation and setting $\tau = \tilde{n}$ we get
\begin{equation*}
\rho^{l+2} L^{\prime \prime}(\rho) + (2l + 2 - \rho) \rho^{l+1}L^{\prime}(\rho) + (\tilde{n} - 1 - l)\rho^{l+1}L(\rho) = 0
\end{equation*}
$\iff$
\begin{equation*}
\rho L^{\prime \prime}(\rho) + (2l + 2 - \rho) L^{\prime}(\rho) + (\tilde{n} - 1 - l) L(\rho) = 0
\end{equation*}
$\iff$
\begin{equation*}
\rho L^{\prime \prime}(\rho) + (\alpha + 1 - \rho) L^{\prime}(\rho) + \tilde{n}^{*} L(\rho) = 0
\end{equation*}
where $\alpha \equiv 2l + 1$ and $\tilde{n}^{*} \equiv \tilde{n} - 1 - l$. This last form is a well known differential equation whose physically acceptable solutions in the present context are associated Laguerre polynomials given by the formula
\begin{equation*}
L_{\tilde{n}^{*}}^{(\alpha)} = \sum_{j = 0}^{\tilde{n}^{*}} (-1)^j \frac{(\tilde{n}^{*} + \alpha)!}{(\tilde{n}^{*}-j)!(\alpha + j)!}\frac{\rho^j}{j!}
\end{equation*}
For given quantum numbers $\tilde{n}$ and $l$, the solution of the radial equation for $R$ is then
\begin{equation*}
R_{\tilde{n}l}(r) \propto e^{-\rho/2}\rho^l L_{\tilde{n} - l - 1}^{(2l + 1)}
\end{equation*}

\section{Final form of the electronic wave function $\psi$}
Putting everything together, for given principal quantum number $\tilde{n}$, orbital quantum number $l$ and magnetic quantum number $m_l$, the wave function of the electron in the hydrogen atom is
\begin{equation*}
\psi_{\tilde{n} l m_l}(r, \theta, \phi) \propto e^{-\rho/2}\rho^l L_{\tilde{n} - l - 1}^{(2l + 1)} P_l^{m_l}(\cos \theta) e^{i m_l \phi}
\end{equation*}
where
\begin{equation*}
\rho = \bigg(-\frac{8 m_e E_{\tilde{n}}}{\hbar^2}\bigg)^{1/2} r
\end{equation*}
and
\begin{equation*}
E_{\tilde{n}} = \bigg(-\frac{m_e}{2 \hbar^2}\bigg) \bigg(\frac{e^2}{4 \pi \epsilon_0}\bigg)^2 \frac{1}{\tilde{n}^2}
\end{equation*}


\chapter{Derivation of exact solutions of the radial equation}\label{second-appendix}

We saw in Appendix~\ref{first-appendix} that for given quantum numbers $\tilde{n}$ and $l$, the solution for the radial equation is 
\begin{equation*}
R_{\tilde{n}l}(r) \propto e^{-\rho/2}\rho^l L_{\tilde{n} - l - 1}^{(2l + 1)}
\end{equation*}
where
\begin{equation*}
\rho = \bigg(-\frac{8 m_e E_{\tilde{n}}}{\hbar^2}\bigg)^{1/2} r
\end{equation*}
\begin{equation*}
E_{\tilde{n}} = \bigg(-\frac{m_e}{2 \hbar^2}\bigg) \bigg(\frac{e^2}{4 \pi \epsilon_0}\bigg)^2 \frac{1}{\tilde{n}^2}
\end{equation*}
and where
\begin{equation*}
L_{\tilde{n}^{*}}^{(\alpha)} = \sum_{j = 0}^{\tilde{n}^{*}} (-1)^j \frac{(\tilde{n}^{*} + \alpha)!}{(\tilde{n}^{*}-j)!(\alpha + j)!}\frac{\rho^j}{j!}
\end{equation*}
are the associated Laguerre polynomials, with $\alpha \equiv 2l + 1$ and $\tilde{n}^{*} \equiv \tilde{n} - 1 - l$. In this Appendix we will use these formulas to derive explicit forms for the first few radial functions ($R_{10}$, $R_{20}$, $R_{30}$, $R_{21}$, $R_{31}$, $R_{32}$), normalised over $r$ so that
\begin{equation*}
\int_0^{\infty} (R_{\tilde{n}l})^2 r^2 dr = 1
\end{equation*} 
(Here, the $r^2$ component comes from the fact that we are using spherical polar coordinates). These exact solutions will then be used to assess the accuracy of the computer approximations in this study. 

For the calculations below it will be convenient to re-express $\rho$ in terms of the Bohr radius
\begin{equation*}
a = \frac{4\pi\epsilon_0 \hbar^2}{e^2 m_e}
\end{equation*}
which is the radius of the innermost Bohr orbit, equal to $5.292 \times 10^{-11}$m. Putting this expression for $a$ into the expression for $E_{\tilde{n}}$ we get
\begin{equation*}
E_{\tilde{n}} = \bigg(-\frac{\hbar^2}{2 m_e}\bigg) \bigg(\frac{1}{\tilde{n}a}\bigg)^2
\end{equation*}
and putting this in turn into the above expression for $\rho$ we get
\begin{equation*}
\rho = \bigg(\frac{2}{\tilde{n}a}\bigg)r
\end{equation*}

\section{$R_{10}$}
With $\tilde{n} = 1$ and $l = 0$ we have
\begin{equation*}
R_{10} \propto e^{-\rho/2} L_{0}^{(1)} = e^{-\rho/2}
\end{equation*} 
and 
\begin{equation*}
\rho = \bigg(\frac{2}{a}\bigg)r
\end{equation*}
The constant of proportionality $A_{10}$ for $R_{10}$ is calculated so that 
\begin{equation*}
\int_0^{\infty} (R_{10})^2 r^2 dr = 1
\end{equation*} 
Therefore we need $A_{10}$ such that 
\begin{equation*}
\int_0^{\infty} A_{10}^2 (e^{-\rho/2})^2 r^2 dr =  A_{10}^2  \bigg(\frac{a}{2}\bigg)^3 \int_0^{\infty} e^{-\rho} \rho^2  d\rho = 1
\end{equation*} 
where in the second integral I have made the change of variable $r = \big(\frac{a}{2}\big)\rho$. Since 
\begin{equation*}
\int_0^{\infty} e^{-\rho} \rho^2  d\rho = 2!
\end{equation*}
we get
\begin{equation*}
A_{10} = \frac{2}{\sqrt{a^3}}
\end{equation*}
Therefore the exact normalised solution for the radial function $R_{10}$ is 
\begin{equation}
R_{10} = \frac{2}{\sqrt{a^3}} e^{-\rho/2} = \frac{2}{\sqrt{a^3}} e^{-r/a}
\end{equation}

\section{$R_{20}$}
With $\tilde{n} = 2$ and $l = 0$ we have
\begin{equation*}
R_{20} \propto e^{-\rho/2} L_{1}^{(1)} = e^{-\rho/2} (2 - \rho)
\end{equation*} 
and 
\begin{equation*}
\rho = \bigg(\frac{1}{a}\bigg)r
\end{equation*}
The constant of proportionality $A_{20}$ for $R_{20}$ is calculated so that 
\begin{equation*}
\int_0^{\infty} (R_{20})^2 r^2 dr = 1
\end{equation*} 
Therefore we need $A_{20}$ such that 
\begin{equation*}
\int_0^{\infty} A_{20}^2 ((2 - \rho) e^{-\rho/2})^2 r^2 dr =  A_{20}^2 a^3 \int_0^{\infty} (4\rho^2 - 4\rho^3 + \rho^4) e^{-\rho} d\rho = 1
\end{equation*} 
where in the second integral I have made the change of variable $r = a \rho$. Since 
\begin{equation*}
\int_0^{\infty} (4\rho^2 - 4\rho^3 + \rho^4) e^{-\rho} d\rho = 4 \cdot 2! - 4 \cdot 3! + 4! = 8
\end{equation*}
we get
\begin{equation*}
A_{20} = \frac{1}{2\sqrt{2a^3}}
\end{equation*}
Therefore the exact normalised solution for the radial function $R_{20}$ is 
\begin{equation}
R_{20} =  \frac{1}{2\sqrt{2a^3}} (2 - \rho) e^{-\rho/2} =  \frac{1}{2\sqrt{2a^3}} \bigg(2 - \frac{r}{a}\bigg) e^{-r/2a}
\end{equation}
\newpage

\section{$R_{30}$}
With $\tilde{n} = 3$ and $l = 0$ we have
\begin{equation*}
R_{30} \propto e^{-\rho/2} L_{2}^{(1)} = \bigg(3 - 3\rho + \frac{1}{2}\rho^2\bigg)  e^{-\rho/2}
\end{equation*} 
and 
\begin{equation*}
\rho = \bigg(\frac{2}{3a}\bigg)r
\end{equation*}
The constant of proportionality $A_{30}$ for $R_{30}$ is calculated so that 
\begin{equation*}
\int_0^{\infty} (R_{30})^2 r^2 dr = 1
\end{equation*} 
Therefore we need $A_{30}$ such that 
\begin{equation*}
\int_0^{\infty} A_{30}^2 \bigg(\bigg(3 - 3\rho + \frac{1}{2}\rho^2\bigg) e^{-\rho/2}\bigg)^2 r^2 dr 
\end{equation*}
\begin{equation*}
=  A_{30}^2 \bigg(\frac{3a}{2}\bigg)^3 \int_0^{\infty} \bigg(9\rho^2 - 18\rho^3 + 12\rho^4 - 3\rho^5 + \frac{1}{4}\rho^6 \bigg) e^{-\rho} d\rho = 1
\end{equation*} 
where in the second integral I have made the change of variable $r = \frac{3a}{2} \rho$. Since 
\begin{equation*}
\int_0^{\infty} \bigg(9\rho^2 - 18\rho^3 + 12\rho^4 - 3\rho^5 + \frac{1}{4}\rho^6 \bigg) e^{-\rho} d\rho = 9 \cdot 2! - 18 \cdot 3! + 12 \cdot 4! - 3 \cdot 5! + \frac{1}{4} \cdot 6! = 9 
\end{equation*}
we get
\begin{equation*}
A_{30} = \frac{2}{9\sqrt{3a^3}}
\end{equation*}
Therefore the exact normalised solution for the radial function $R_{30}$ is 
\begin{equation}
R_{30} =  \frac{2}{9\sqrt{3a^3}} \bigg(3 - 3\rho + \frac{1}{2}\rho^2\bigg) e^{-\rho/2} =  \frac{2}{81\sqrt{3a^3}} \bigg(27 - 18\frac{r}{a} + 2\frac{r^2}{a^2}\bigg) e^{-r/3a}
\end{equation}

\section{$R_{21}$}
With $\tilde{n} = 2$ and $l = 1$ we have
\begin{equation*}
R_{21} \propto e^{-\rho/2}\rho L_{0}^{(3)} = e^{-\rho/2} \rho 
\end{equation*} 
and 
\begin{equation*}
\rho = \bigg(\frac{1}{a}\bigg)r
\end{equation*}
The constant of proportionality $A_{21}$ for $R_{21}$ is calculated so that 
\begin{equation*}
\int_0^{\infty} (R_{21})^2 r^2 dr = 1
\end{equation*} 
Therefore we need $A_{21}$ such that 
\begin{equation*}
\int_0^{\infty} A_{21}^2 (e^{-\rho/2}\rho)^2 r^2 dr =  A_{21}^2 a^3 \int_0^{\infty} e^{-\rho} \rho^4 d\rho = 1
\end{equation*} 
where in the second integral I have made the change of variable $r = a \rho$. Since 
\begin{equation*}
\int_0^{\infty} e^{-\rho} \rho^4 d\rho = 4! = 24
\end{equation*}
we get
\begin{equation*}
A_{21} = \frac{1}{2\sqrt{6a^3}}
\end{equation*}
Therefore the exact normalised solution for the radial function $R_{21}$ is 
\begin{equation}
R_{21} =  \frac{1}{2\sqrt{6a^3}} \rho e^{-\rho/2} =  \frac{1}{2\sqrt{6a^3}} \frac{r}{a} e^{-r/2a}
\end{equation}

\section{$R_{31}$}
With $\tilde{n} = 3$ and $l = 1$ we have
\begin{equation*}
R_{31} \propto e^{-\rho/2}\rho  L_{1}^{(3)} = e^{-\rho/2} (4 - \rho) \rho
\end{equation*} 
and 
\begin{equation*}
\rho = \bigg(\frac{2}{3a}\bigg)r
\end{equation*}
The constant of proportionality $A_{31}$ for $R_{31}$ is calculated so that 
\begin{equation*}
\int_0^{\infty} (R_{31})^2 r^2 dr = 1
\end{equation*} 
Therefore we need $A_{31}$ such that 
\begin{equation*}
\int_0^{\infty} A_{31}^2 ((4 - \rho) \rho e^{-\rho/2})^2 r^2 dr =  A_{31}^2 \bigg(\frac{3a}{2}\bigg)^3 \int_0^{\infty} (16\rho^4 - 8\rho^5 + \rho^6) e^{-\rho} d\rho = 1
\end{equation*} 
where in the second integral I have made the change of variable $r = \frac{3a}{2} \rho$. Since 
\begin{equation*}
\int_0^{\infty} (16\rho^4 - 8\rho^5 + \rho^6) e^{-\rho} d\rho = 16 \cdot 4! - 8 \cdot 5! + 6! = 18
\end{equation*}
we get
\begin{equation*}
A_{31} = \frac{1}{9\sqrt{6a^3}}
\end{equation*}
Therefore the exact normalised solution for the radial function $R_{31}$ is 
\begin{equation}
R_{31} =  \frac{1}{9\sqrt{6a^3}} (4 - \rho) \rho e^{-\rho/2} =  \frac{4}{81\sqrt{6a^3}} \bigg(6 - \frac{r}{a}\bigg)\frac{r}{a} e^{-r/3a}
\end{equation}

\section{$R_{32}$}
Finally, with $\tilde{n} = 3$ and $l = 2$ we have
\begin{equation*}
R_{32} \propto e^{-\rho/2}\rho^2  L_{0}^{(5)} = e^{-\rho/2} \rho^2
\end{equation*} 
and 
\begin{equation*}
\rho = \bigg(\frac{2}{3a}\bigg)r
\end{equation*}
The constant of proportionality $A_{32}$ for $R_{32}$ is calculated so that 
\begin{equation*}
\int_0^{\infty} (R_{32})^2 r^2 dr = 1
\end{equation*} 
Therefore we need $A_{32}$ such that 
\begin{equation*}
\int_0^{\infty} A_{32}^2 (\rho^2 e^{-\rho/2})^2 r^2 dr =  A_{32}^2 \bigg(\frac{3a}{2}\bigg)^3 \int_0^{\infty} e^{-\rho}\rho^6 d\rho = 1
\end{equation*} 
where in the second integral I have made the change of variable $r = \frac{3a}{2} \rho$. Since 
\begin{equation*}
\int_0^{\infty} e^{-\rho} \rho^6 d\rho = 6! = 720
\end{equation*}
we get
\begin{equation*}
A_{32} = \frac{1}{9\sqrt{30a^3}}
\end{equation*}
Therefore the exact normalised solution for the radial function $R_{32}$ is 
\begin{equation}
R_{32} = \frac{1}{9\sqrt{30a^3}} \rho^2 e^{-\rho/2} =  \frac{4}{81\sqrt{30a^3}} \bigg(\frac{r}{a}\bigg)^2 e^{-r/3a}
\end{equation}


\chapter{Maple code for plotting B-spline sets}\label{third-appendix}
The sample code here is for the top left plot in Figure~\ref{fig:bsplinesets}, with $n=7$ and $k=3$.

\includegraphics{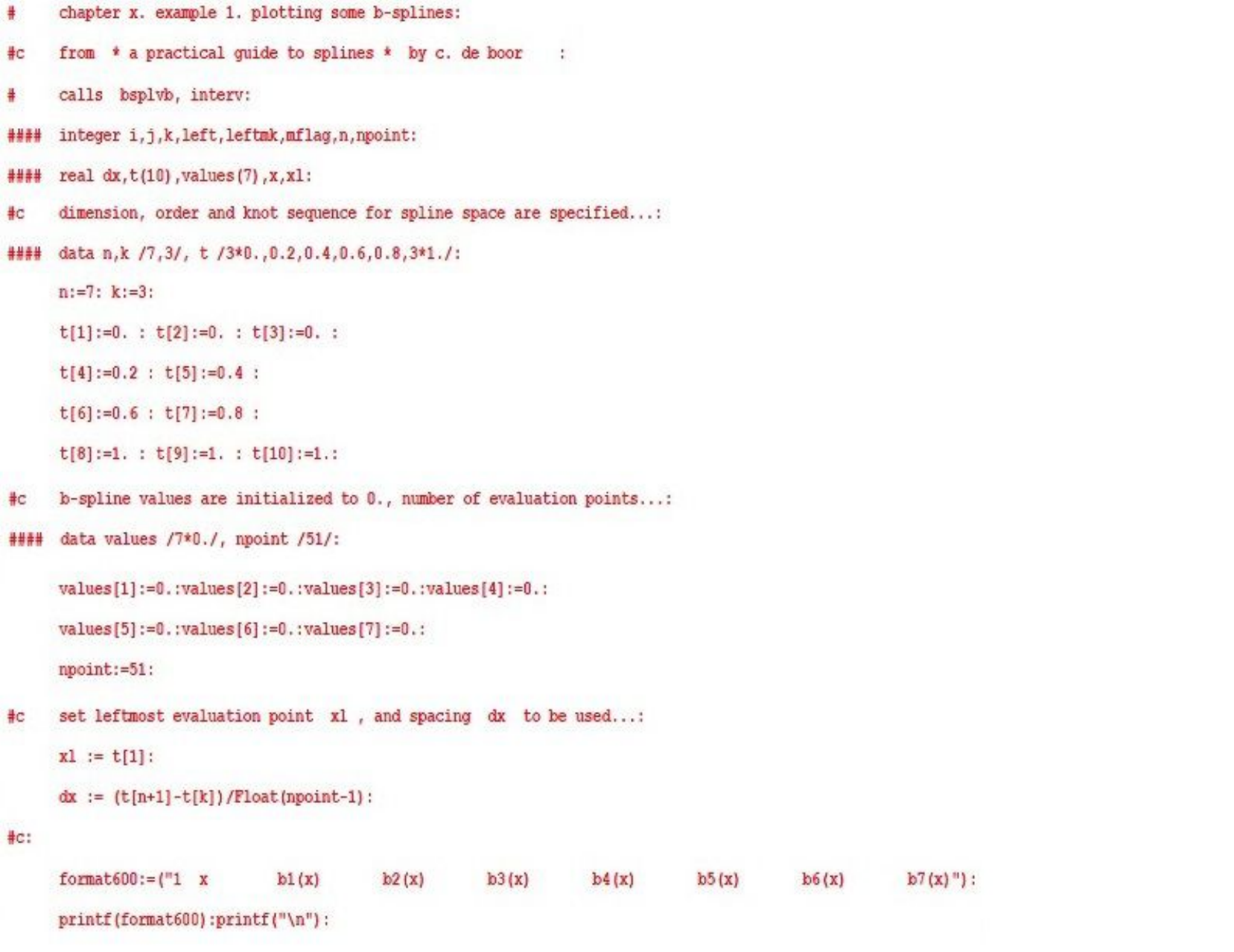} 
\newpage
\includegraphics{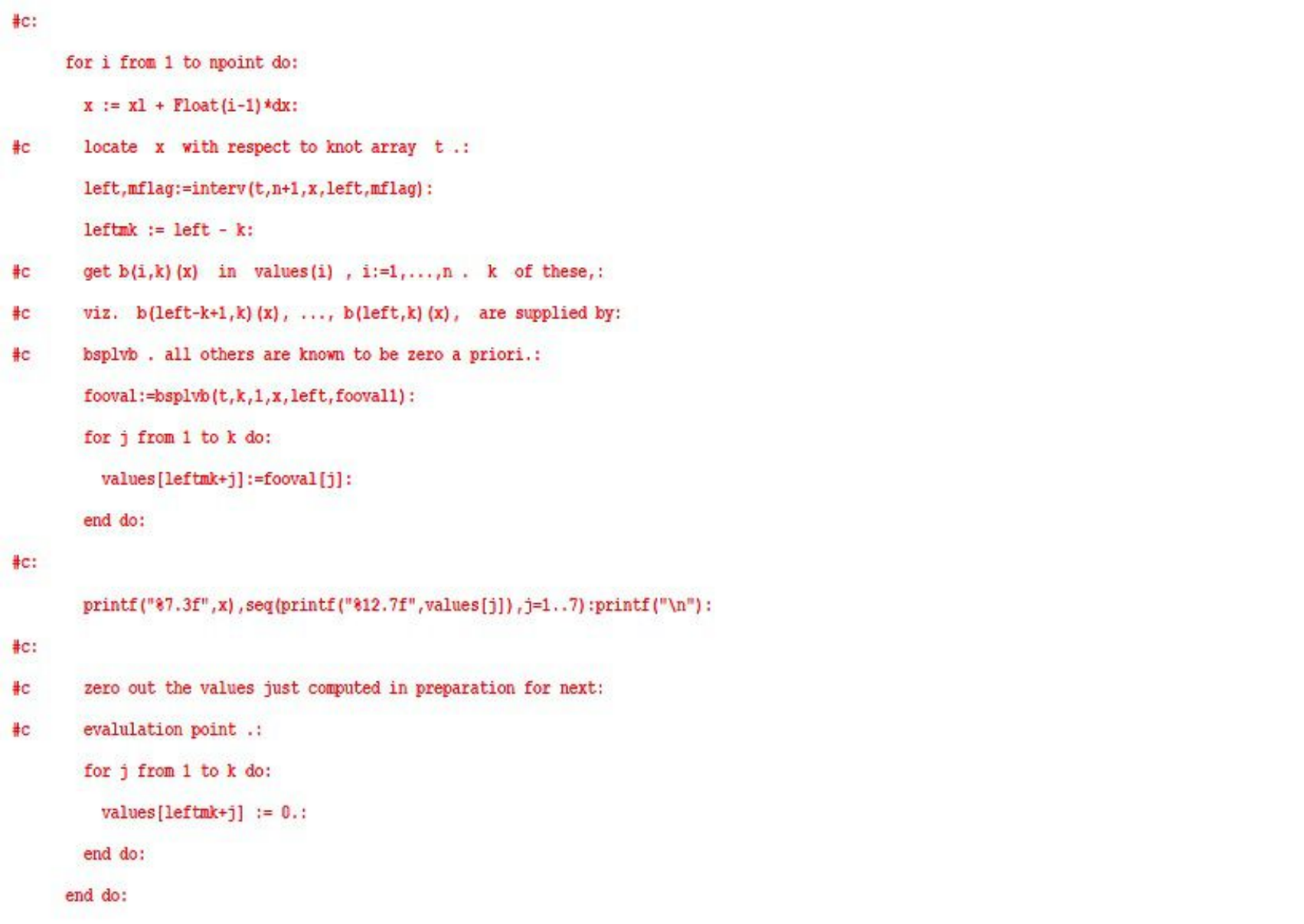} 
\newpage
\includegraphics{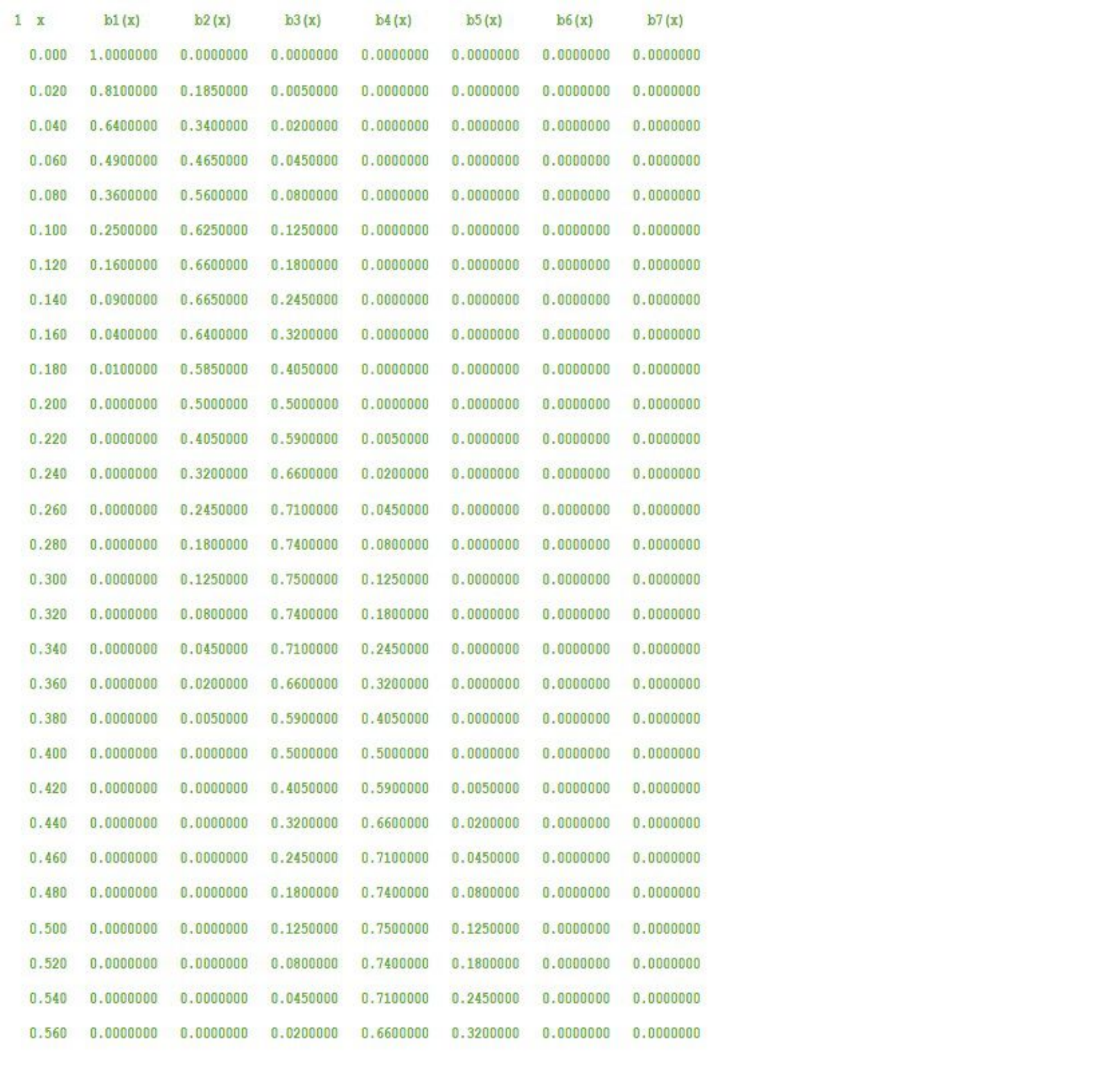}
\newpage 
\includegraphics{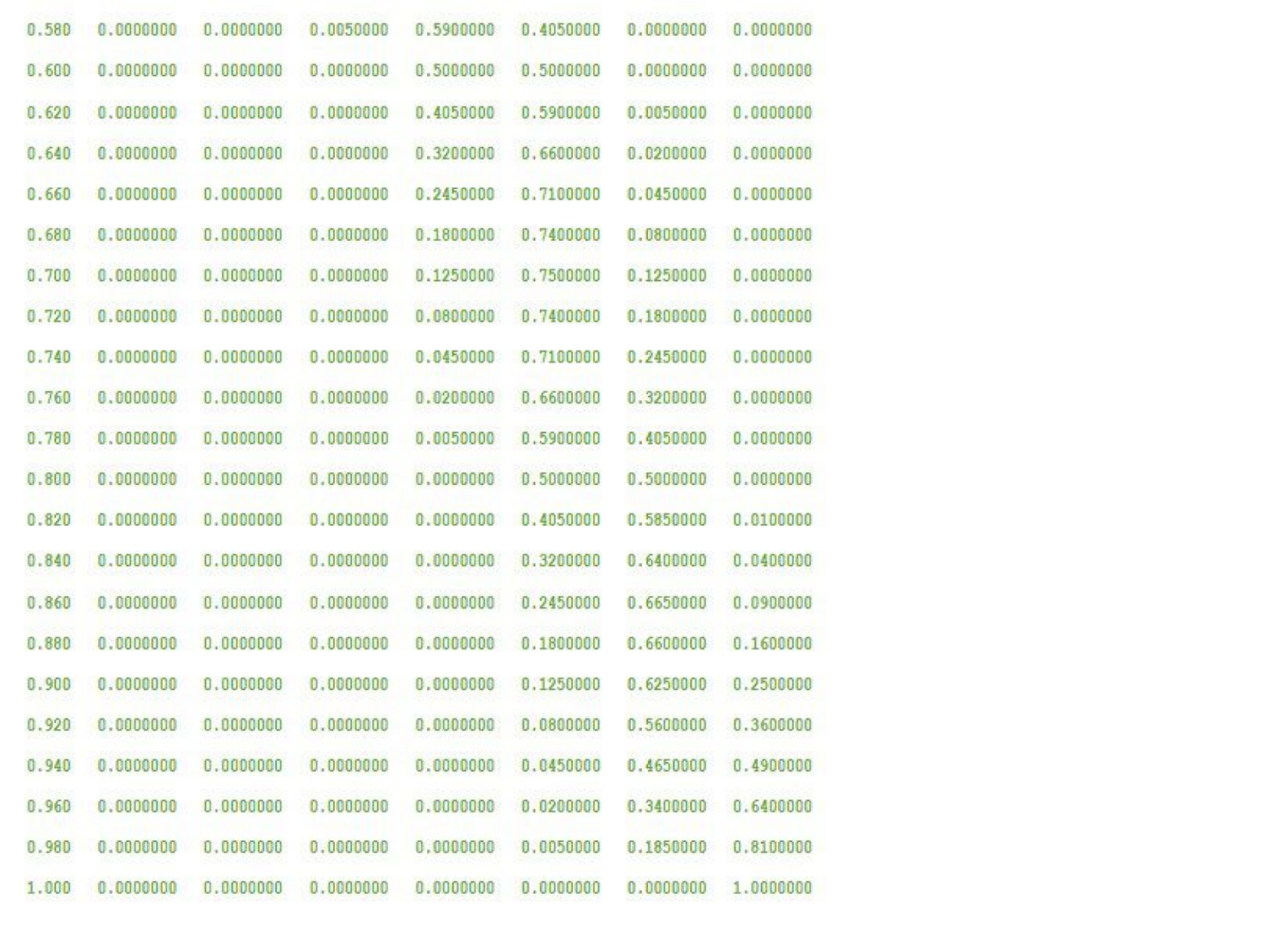}
\newpage 
\includegraphics{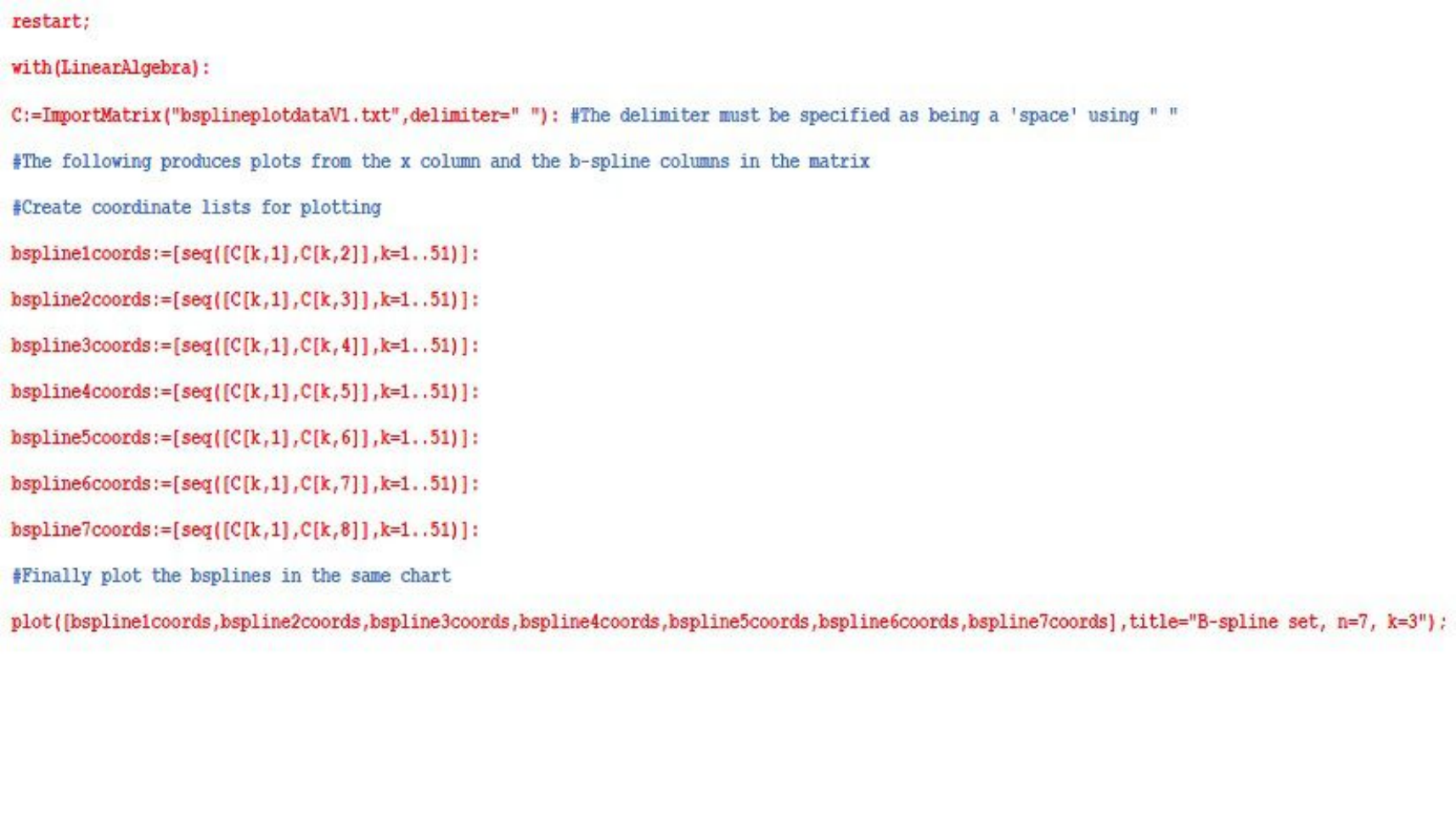} 


\chapter{Maple code for spline interpolation}\label{fourth-appendix}
The code presented here is for the cubic spline interpolation in Figure~\ref{fig:hypergeometric}.

\includegraphics{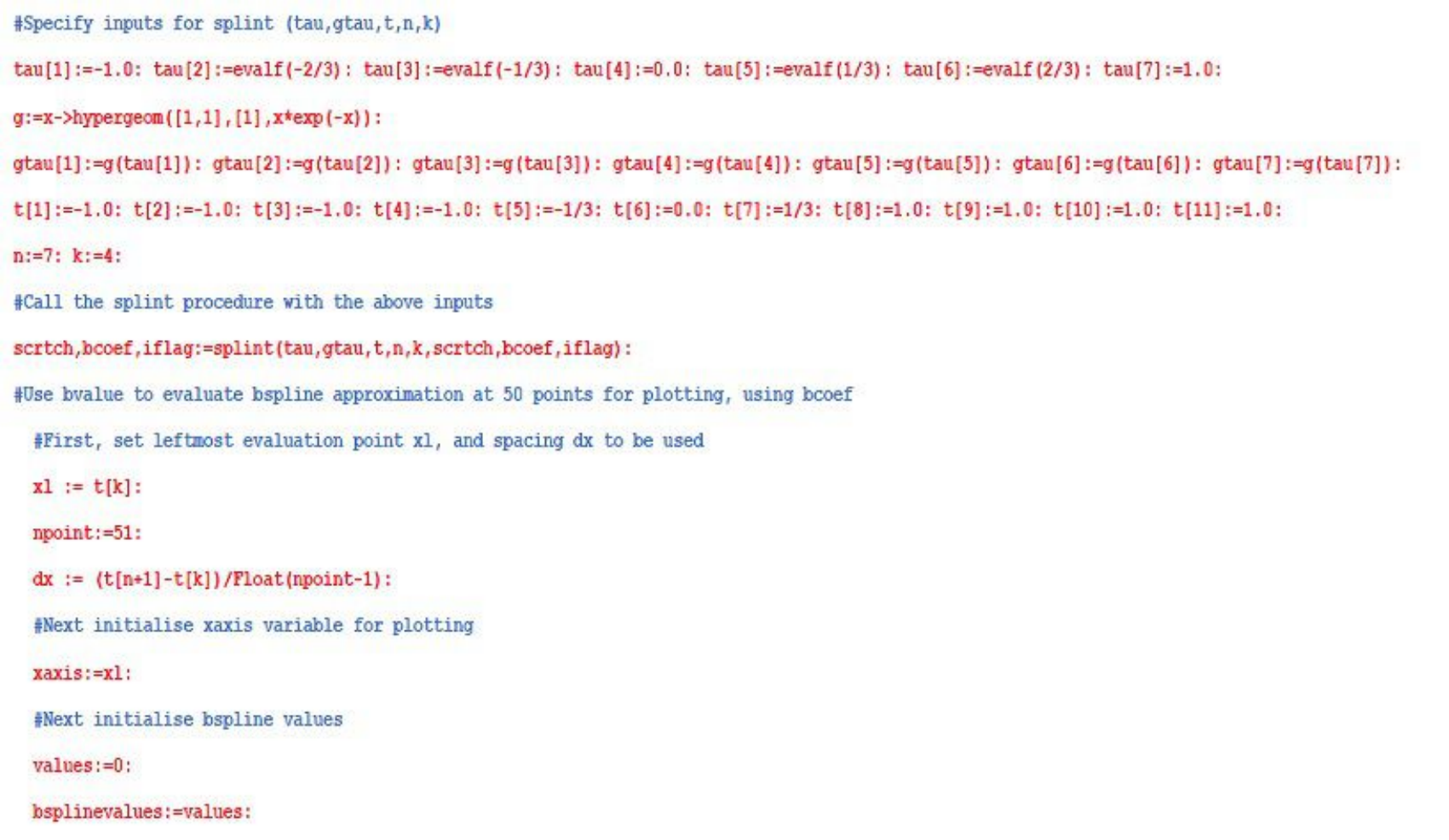}
\newpage 
\includegraphics{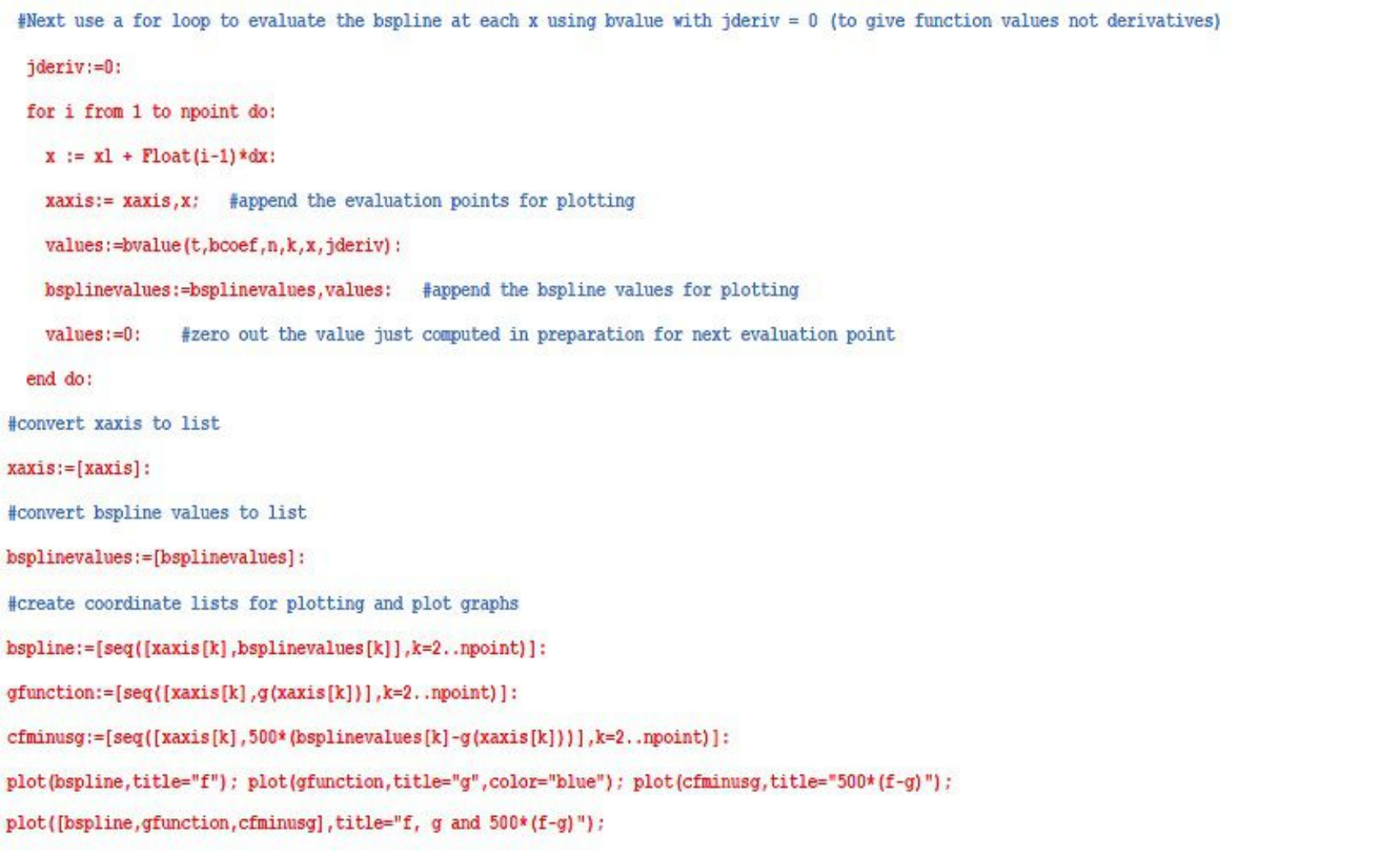}


\chapter{Maple code for wave functions with zero angular momentum}\label{fifth-appendix}
The sample code presented here is for the experiment with box [0, 10], 10 intervals and two collocation points within each interval in Section 4.1. Modifications were made to COLPNT and DIFEQU, and also shown is the Maple code used for post-output processing after calling COLLOC. Using the notation of Chapter XV of de Boor, we thus have m = 2 (two side conditions), k = 2 collocation points per interval, and k + m = 4 (the order of the approximation is 4 so we are using cubic approximations). We will be using knot sequence (0, 1, 2, ..., 10) within the interval [0, 10] and there will therefore be 10 polynomial pieces. 

In the procedure COLPNT, in order to use equally spaced collocation points instead of Gaussian points within each interval, we will need to replace the Gaussian collocation points in the section for k = 2 with the two equally spaced points -1/3  and 1/3 in the interval [-1, 1] (see my annotations in the COLPNT procedure below). These will then transform into two equally spaced collocation points within each interval. A number of alterations are needed in the DIFEQU procedure.

\includegraphics{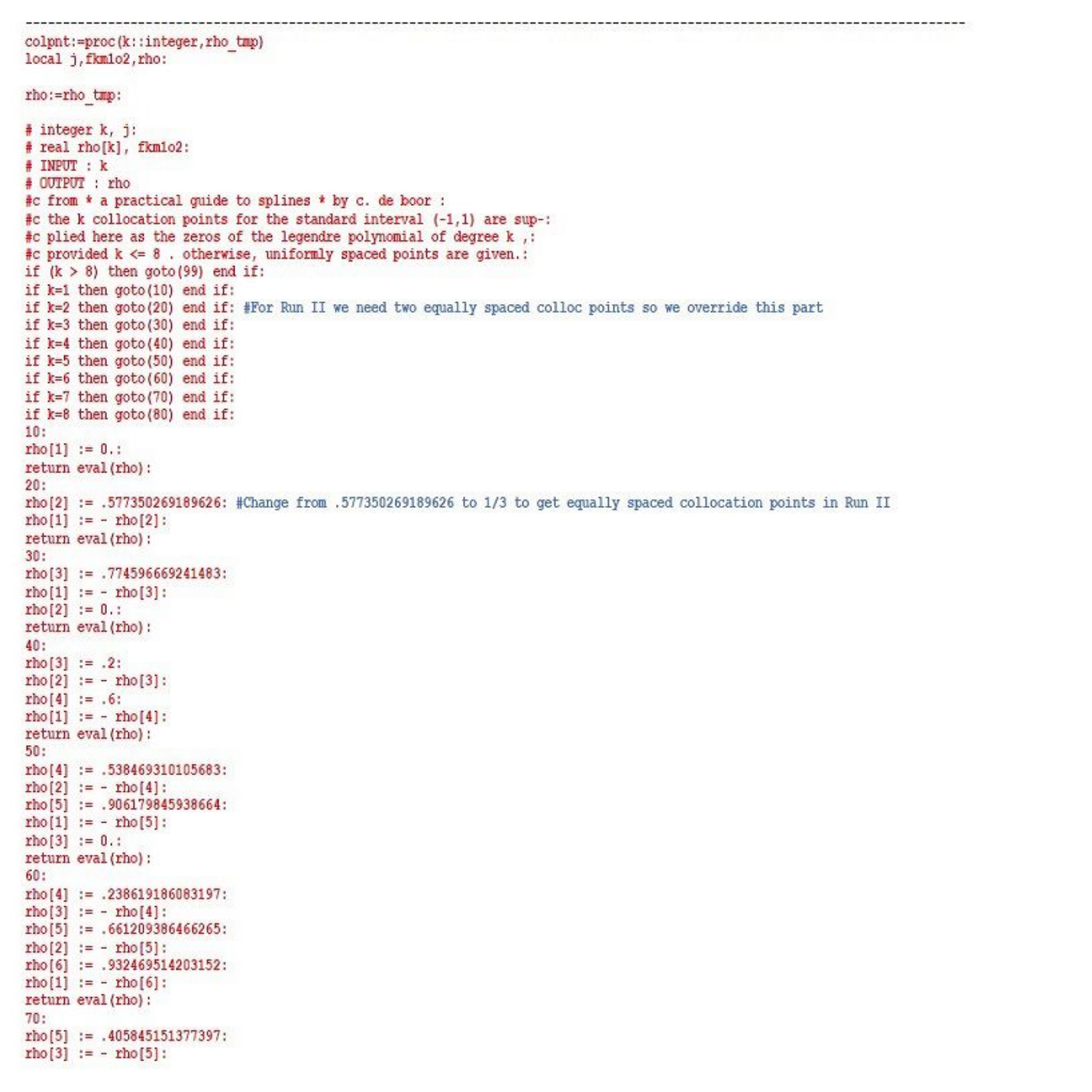}
\newpage
\includegraphics{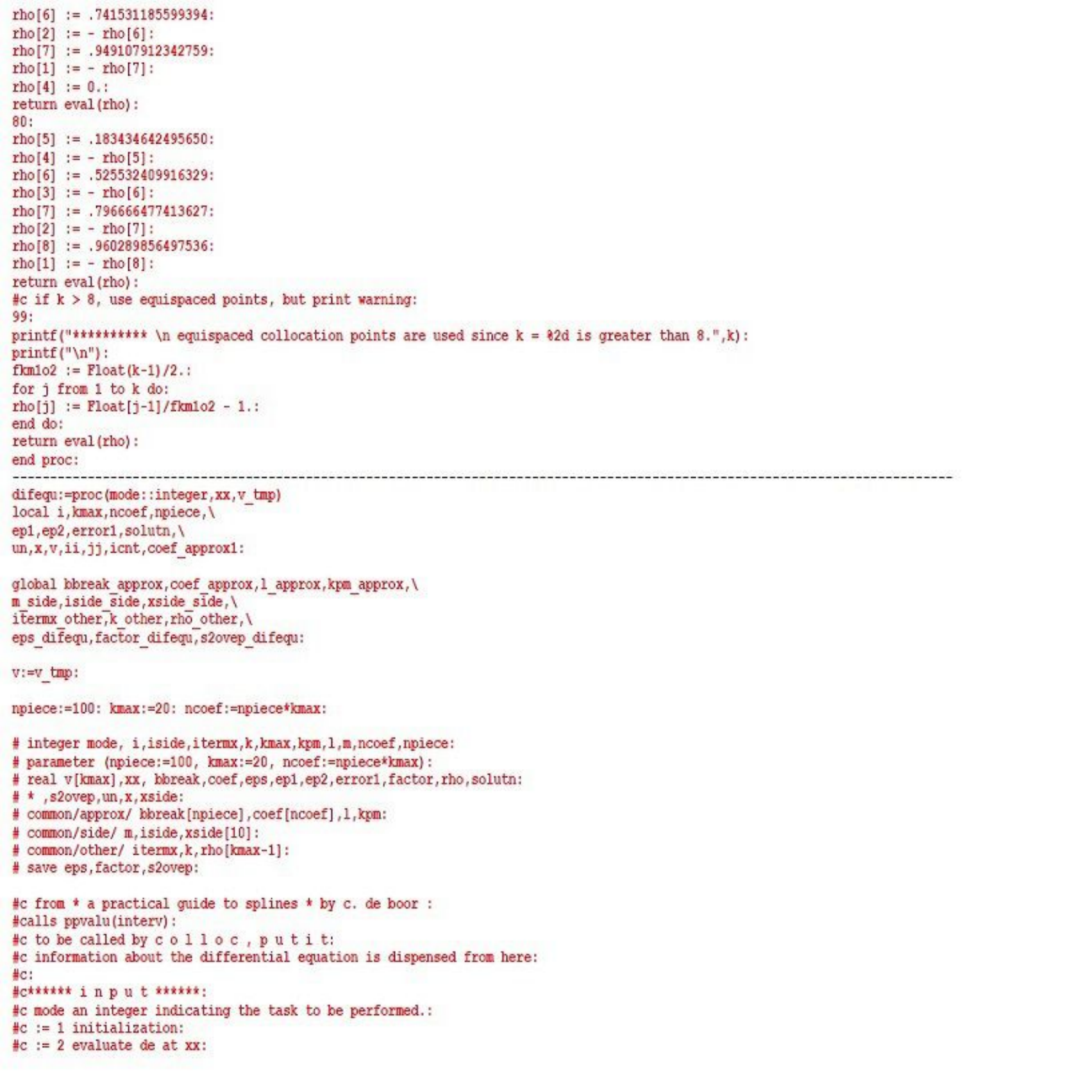}
\newpage
\includegraphics{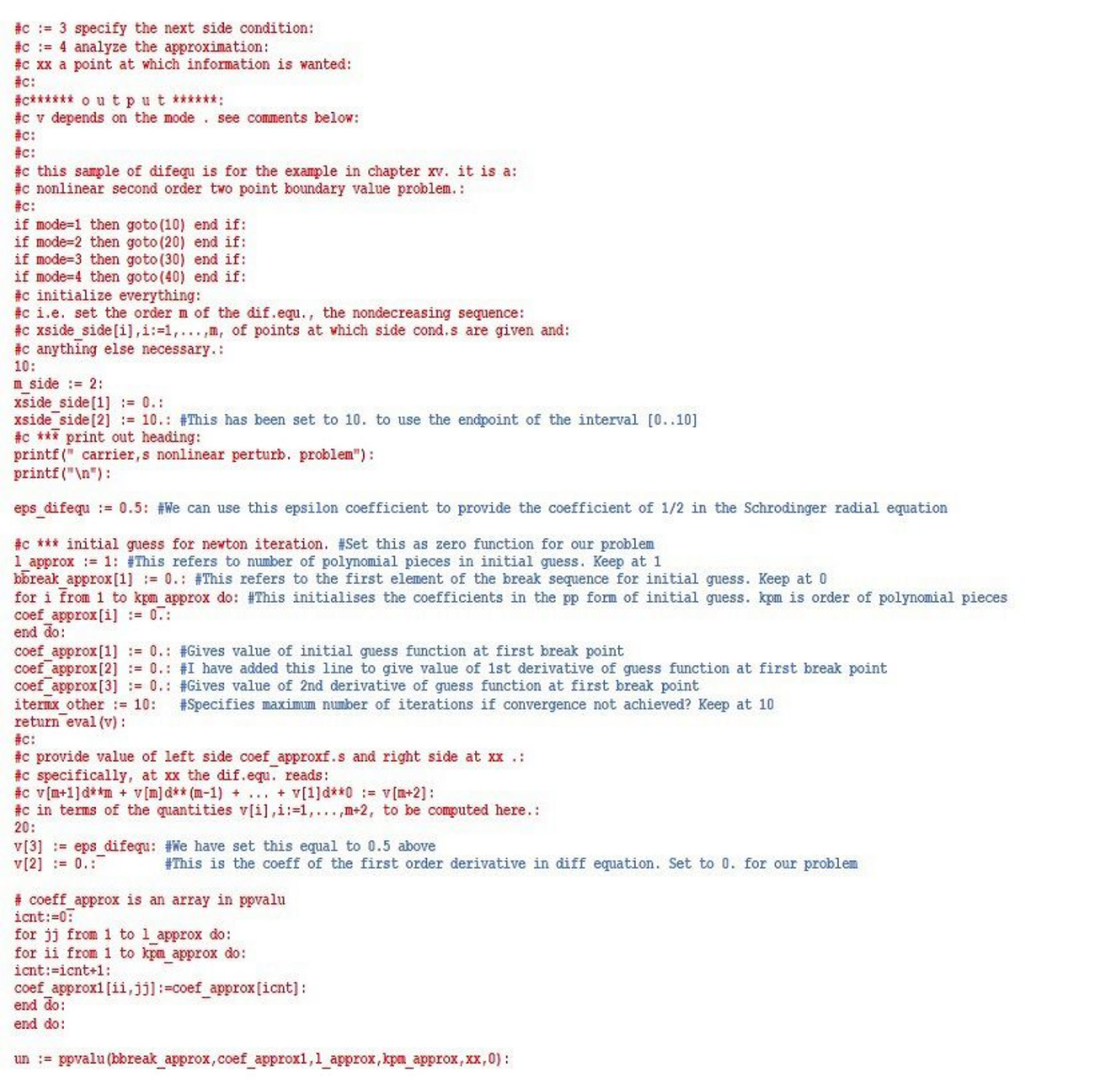}
\newpage
\includegraphics{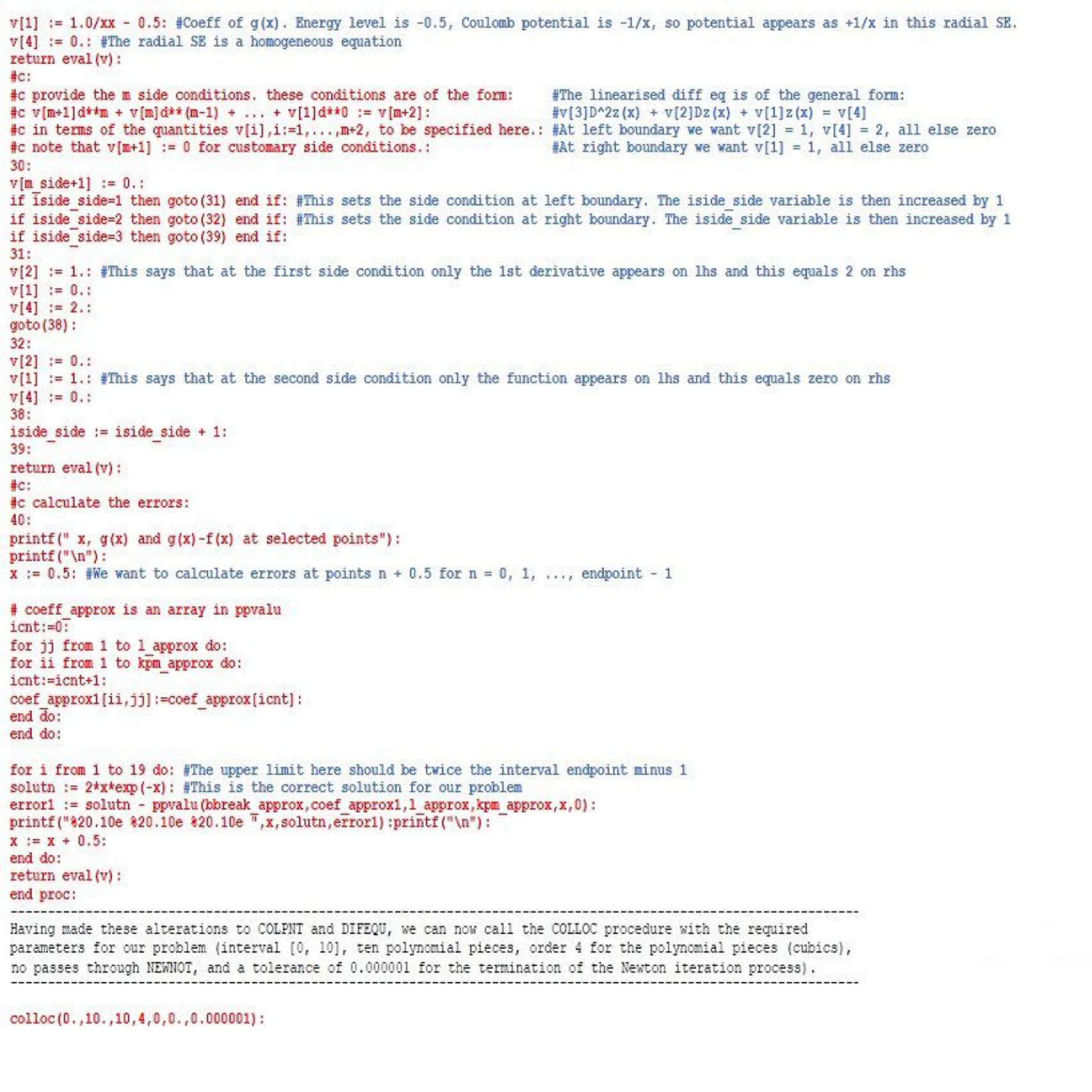}
\newpage
\includegraphics{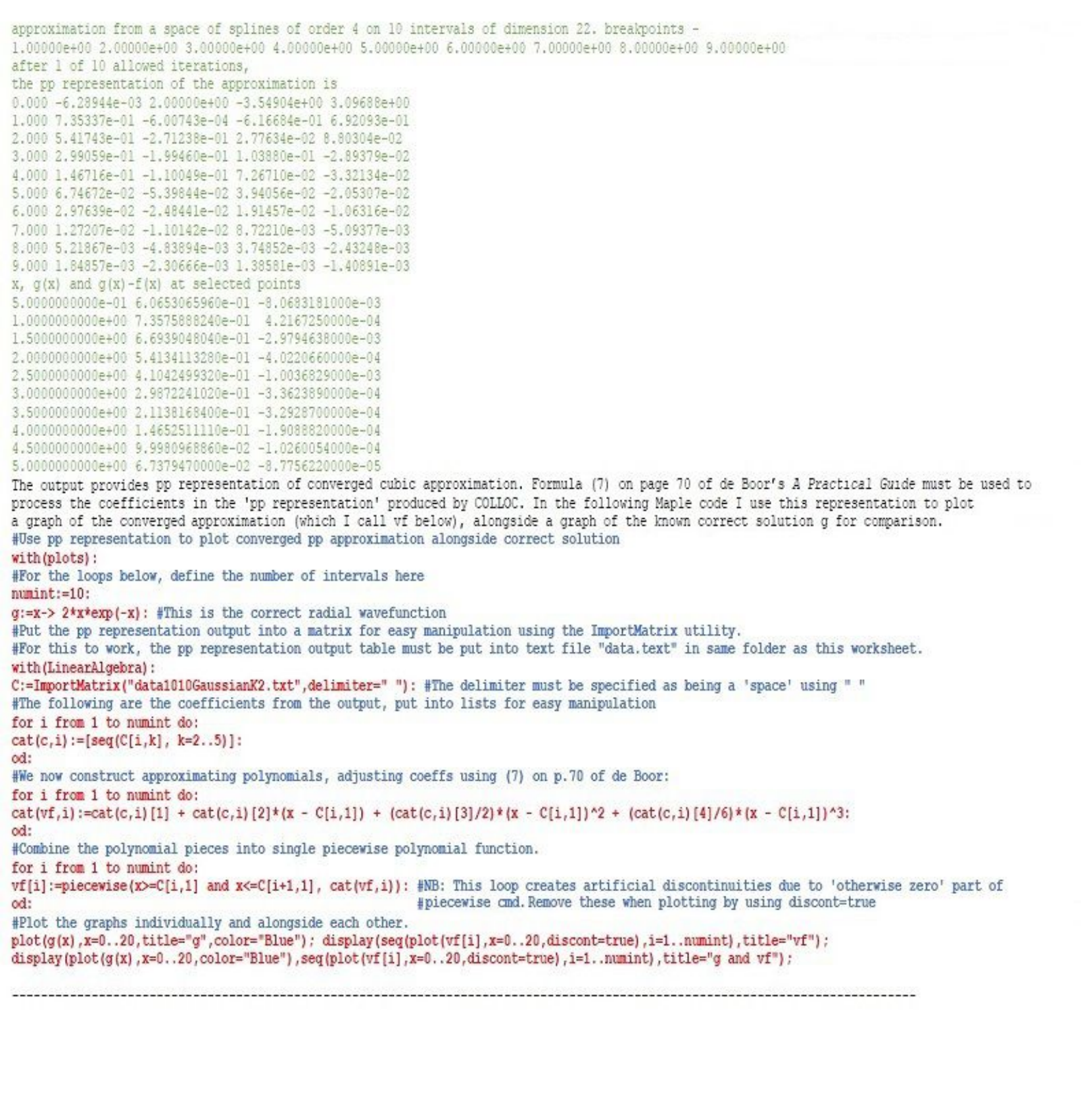}


\chapter{Maple code for wave functions with angular momentum}\label{sixth-appendix}

The following shows the modifications made to DIFEQU and the post-output processing code for the experiment with box [0, 50], 30 intervals and six collocation points within each interval in Section 4.2. Using the notation of Chapter XV of de Boor, we have m = 2 (two side conditions), k = 6 collocation points per interval, and k + m = 8 (the order of the approximation is 8 so we will be using heptic approximations). We will be using 30 equally spaced knots within the interval [0, 50] so there will be 30 polynomial pieces. In the procedure COLPNT (not shown here), in order to use equally spaced collocation points instead of Gaussian points within each interval, we will need to replace the Gaussian collocation points in the section for k = 6 with the six equally spaced points -5/7, -3/7, -1/7, 1/7, 3/7 and 5/7 in the interval [-1, 1]. These will then transform into six equally spaced collocation points within each interval.

\includegraphics{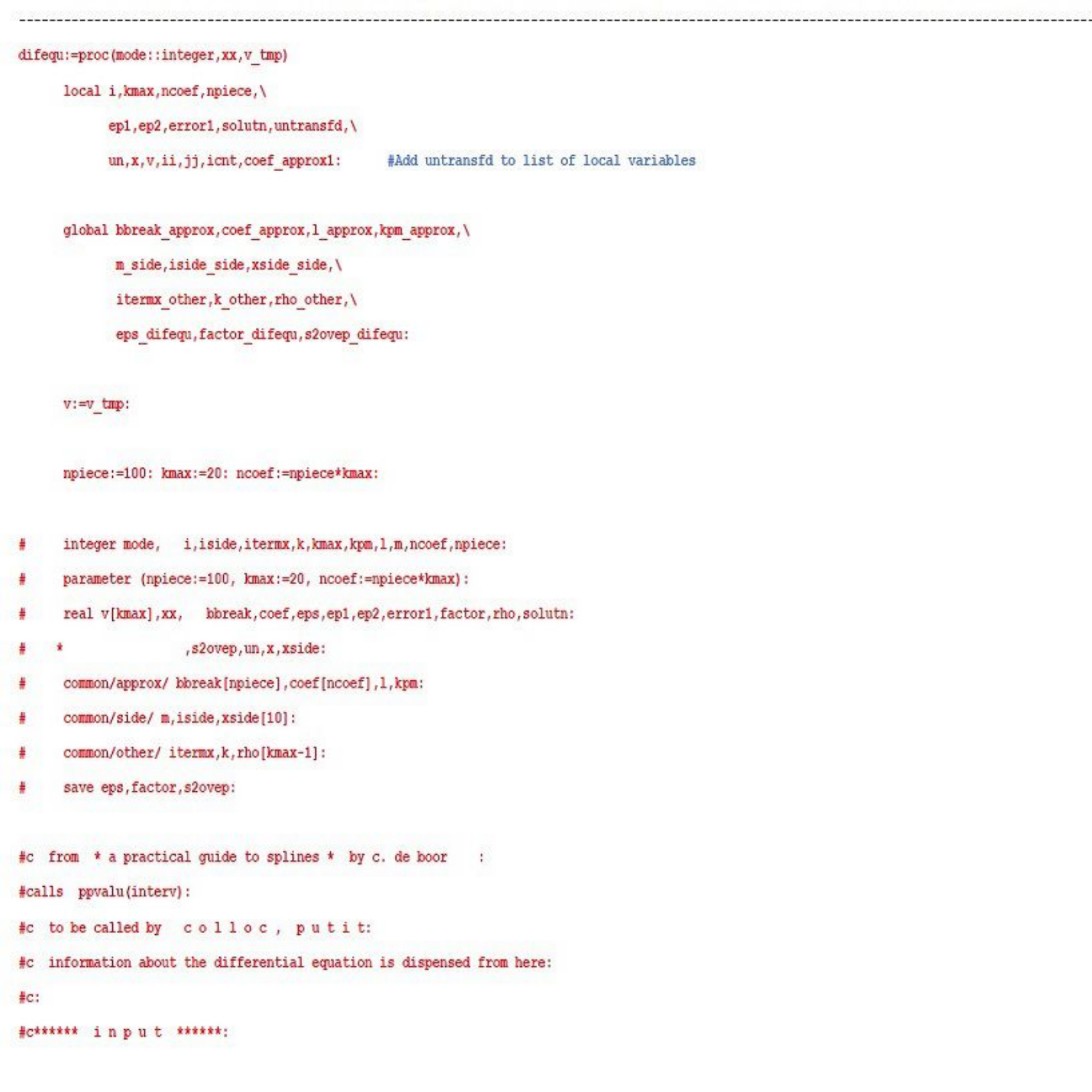}
\newpage
\includegraphics{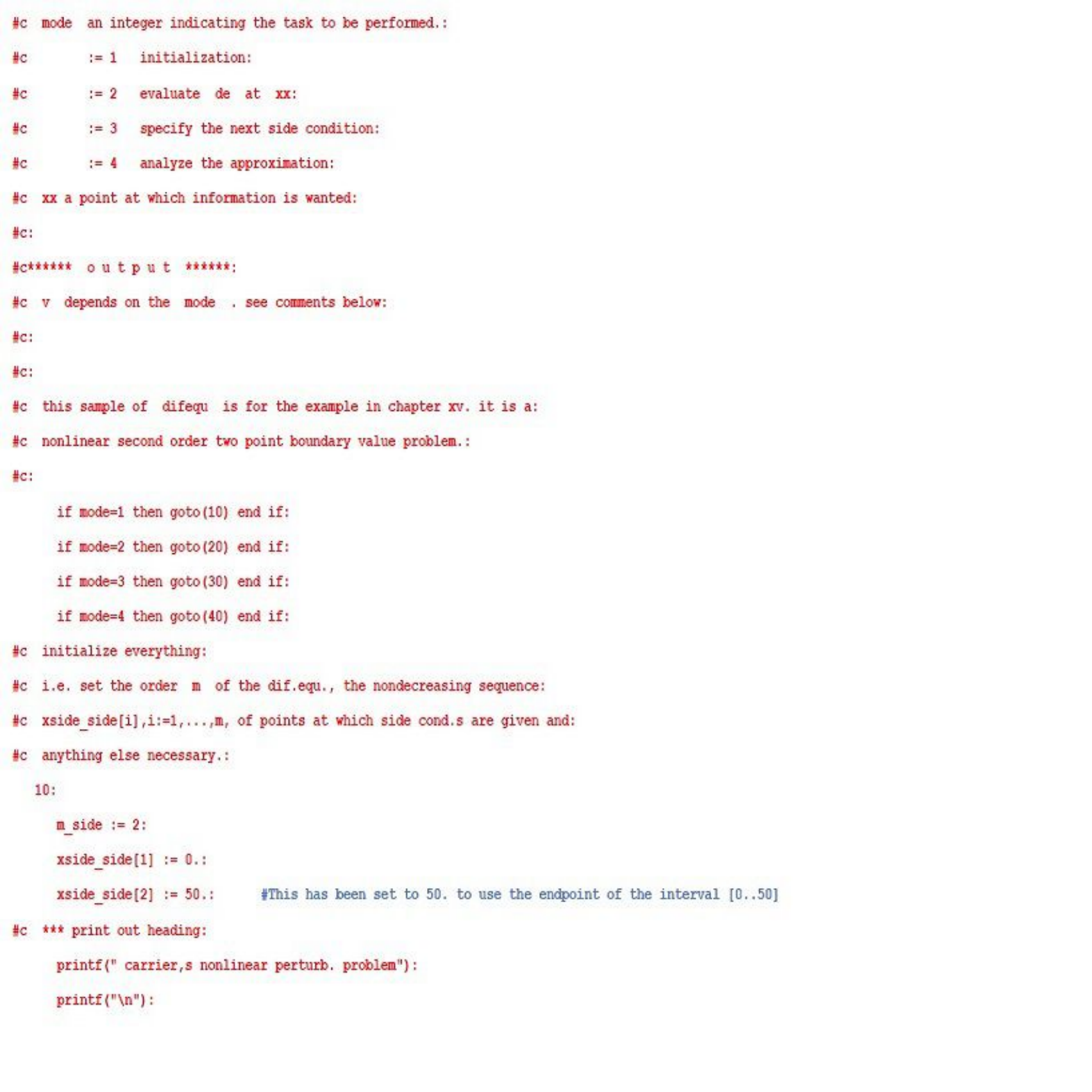}
\newpage
\includegraphics{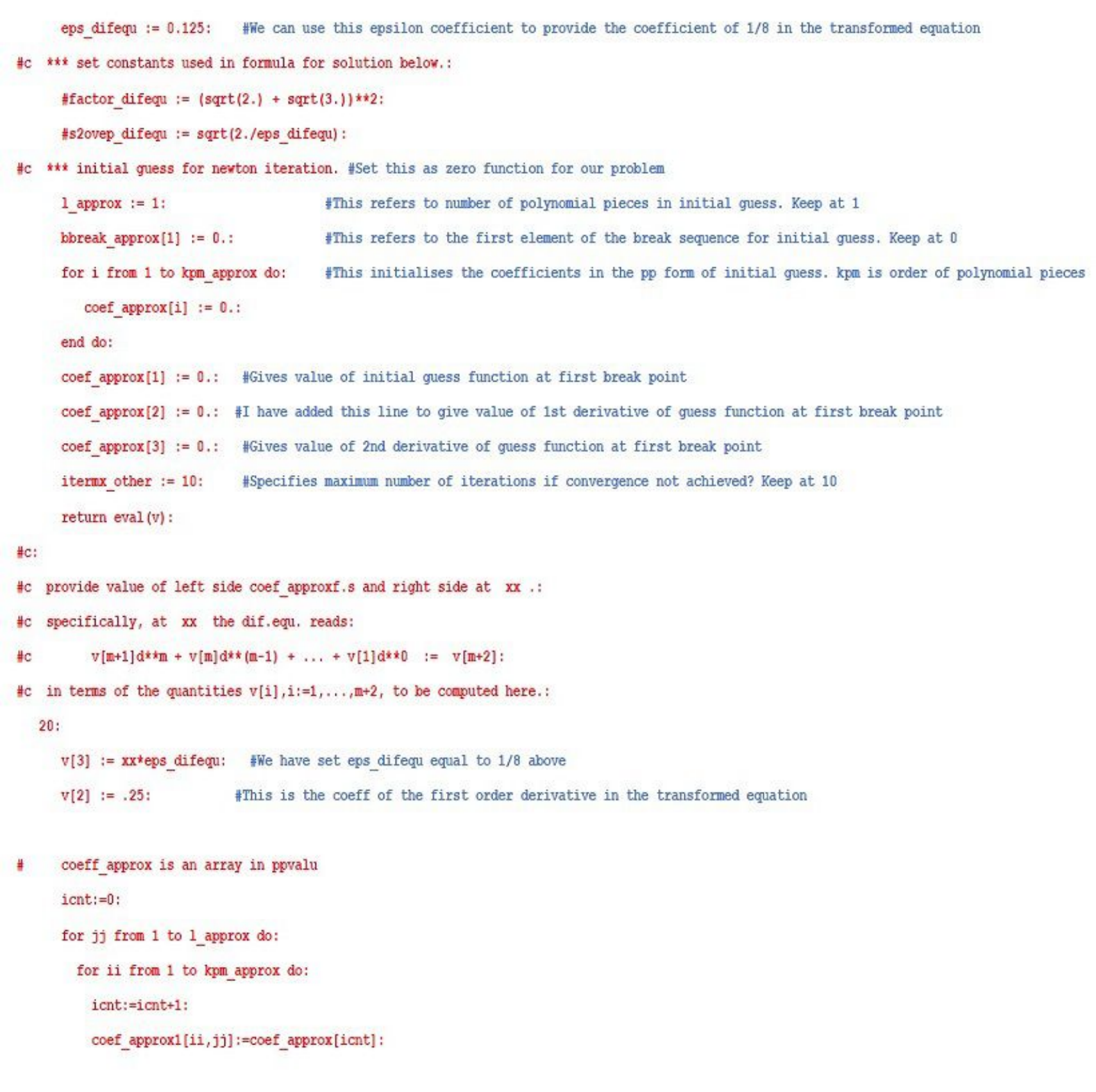}
\newpage
\includegraphics{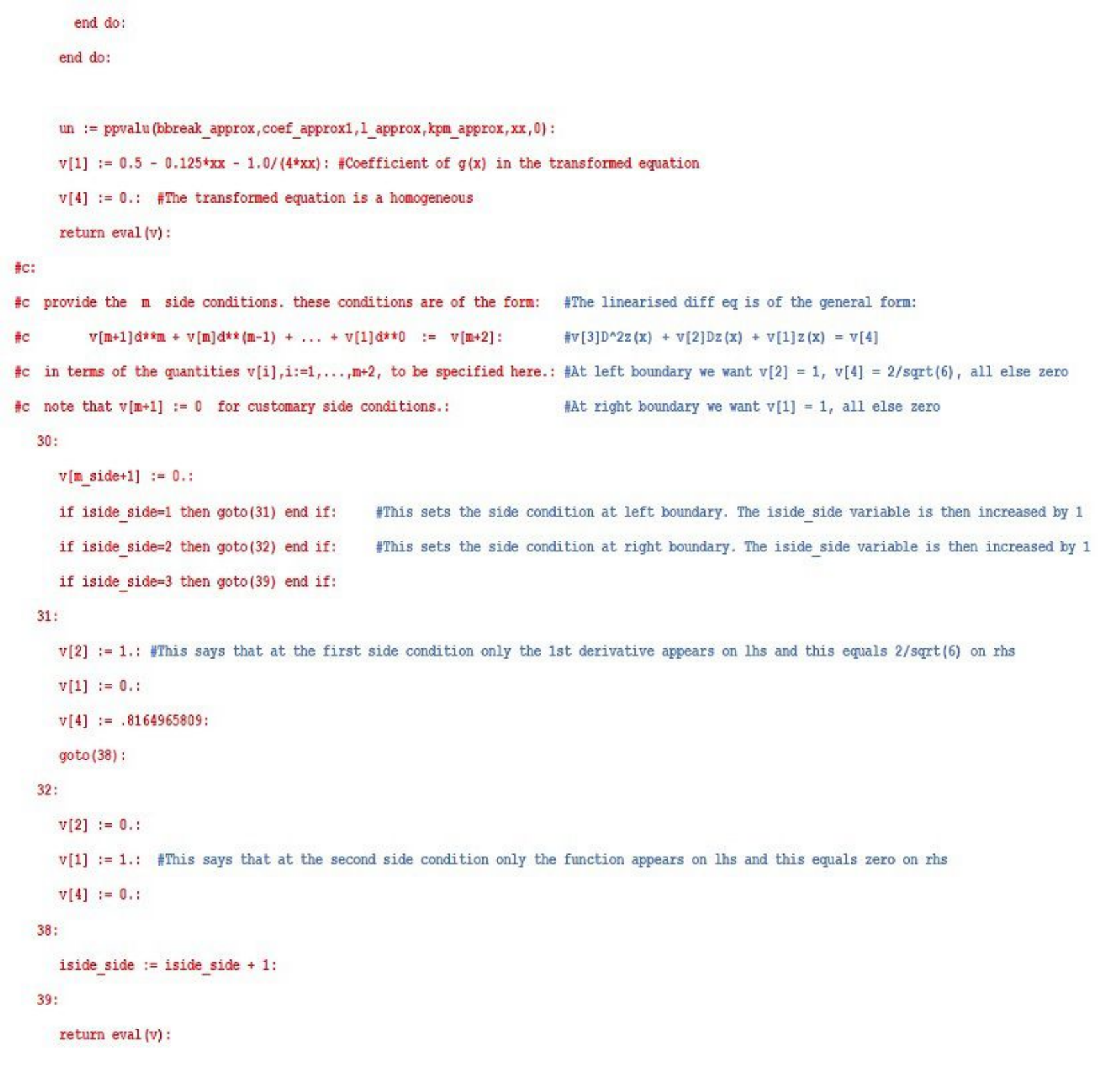}
\newpage
\includegraphics{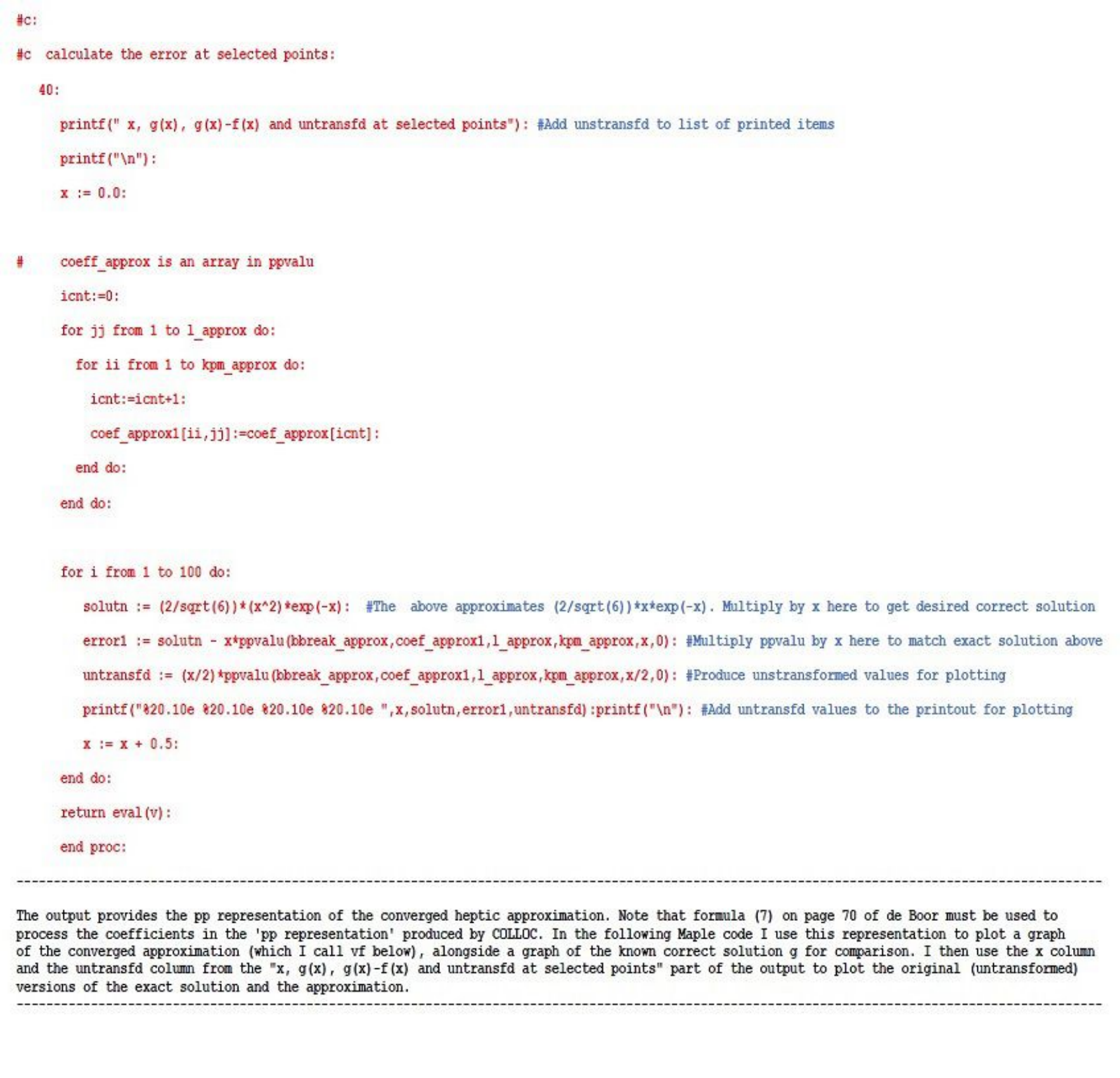}
\newpage
\includegraphics{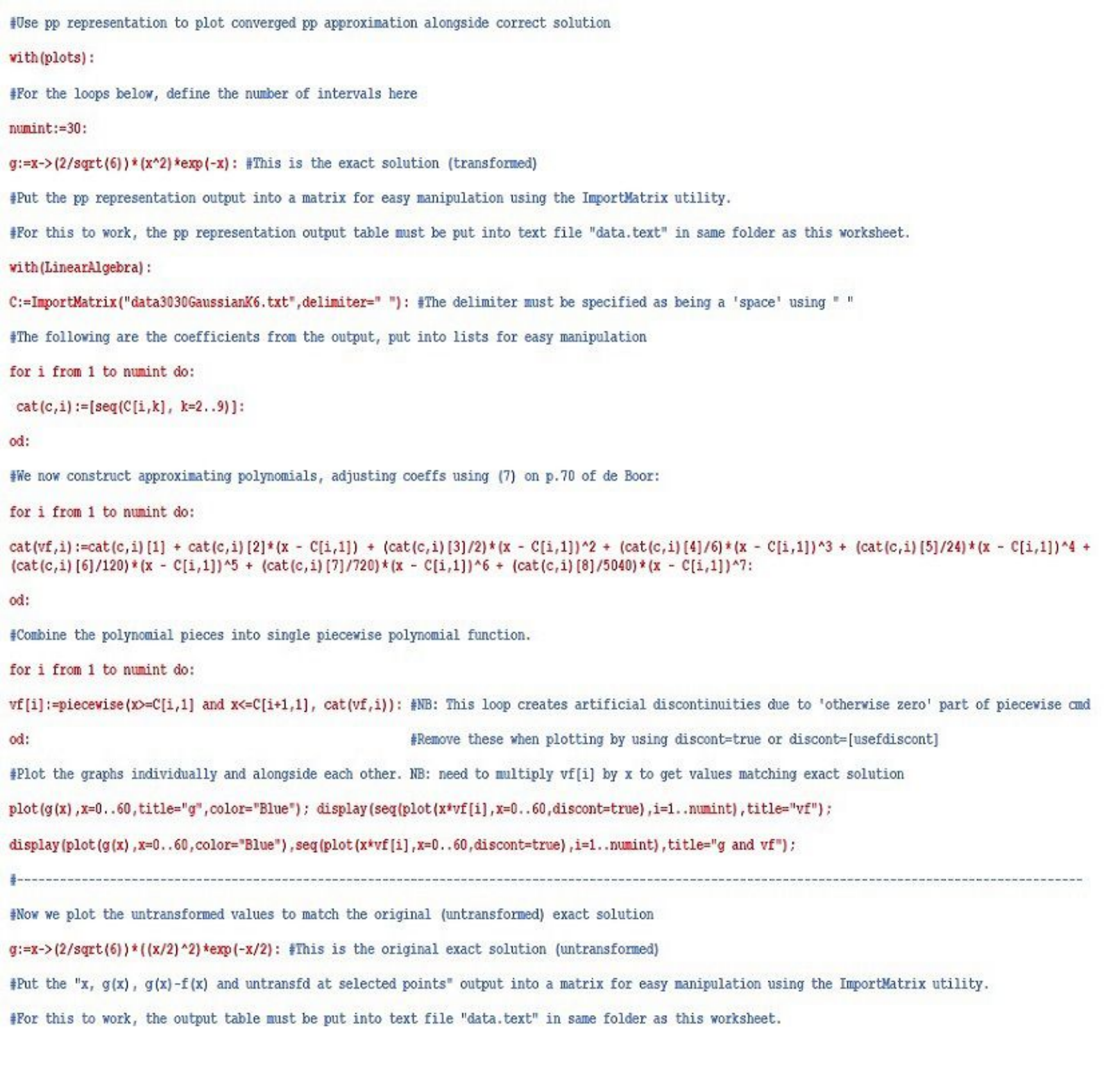}
\newpage
\includegraphics{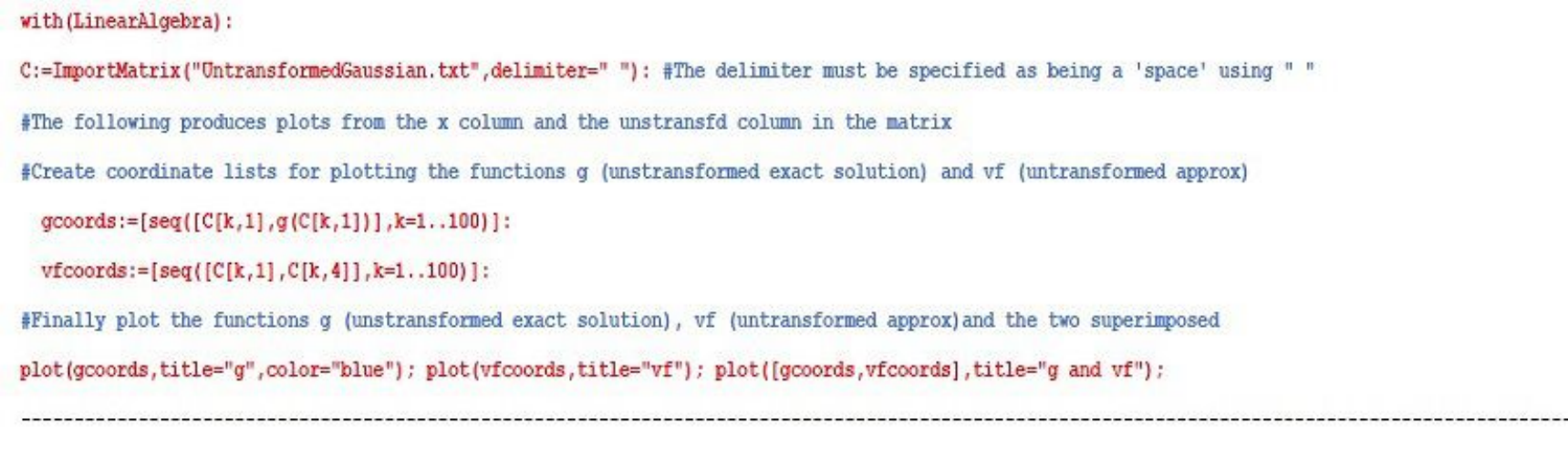}
  

\chapter{Maple code for nonlinear Schr\"{o}dinger equation}\label{seventh-appendix}

The following shows the modifications made to DIFEQU and the post-output processing code for the experiment with $\epsilon = 0.1$, box [0, 1], 20 intervals and six collocation points within each interval in Section 5.1. Using the notation of Chapter XV of de Boor, we have m = 2 (two side conditions), k = 6 collocation points per interval, and k + m = 8 (the order of the approximation is 8 so we will be using heptic approximations). We will be using 20 equally spaced knots within the interval [0, 1] so there will be 20 polynomial pieces. In the procedure COLPNT (not shown here), in order to use equally spaced collocation points instead of Gaussian points within each interval, we will need to replace the Gaussian collocation points in the section for k = 6 with the six equally spaced points -5/7, -3/7, -1/7, 1/7, 3/7 and 5/7 in the interval [-1, 1]. These will then transform into six equally spaced collocation points within each interval.

\includegraphics{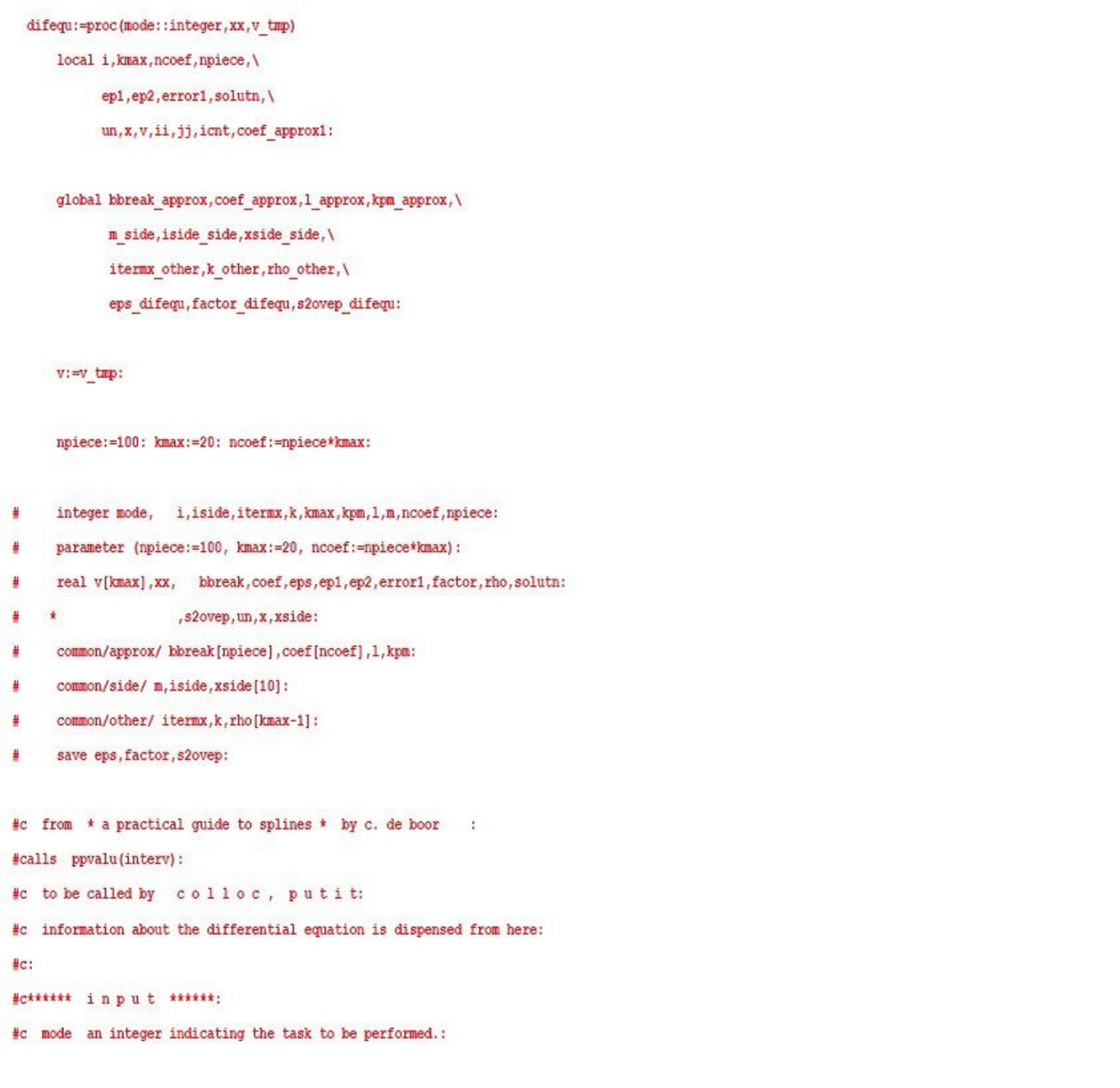}
\newpage
\includegraphics{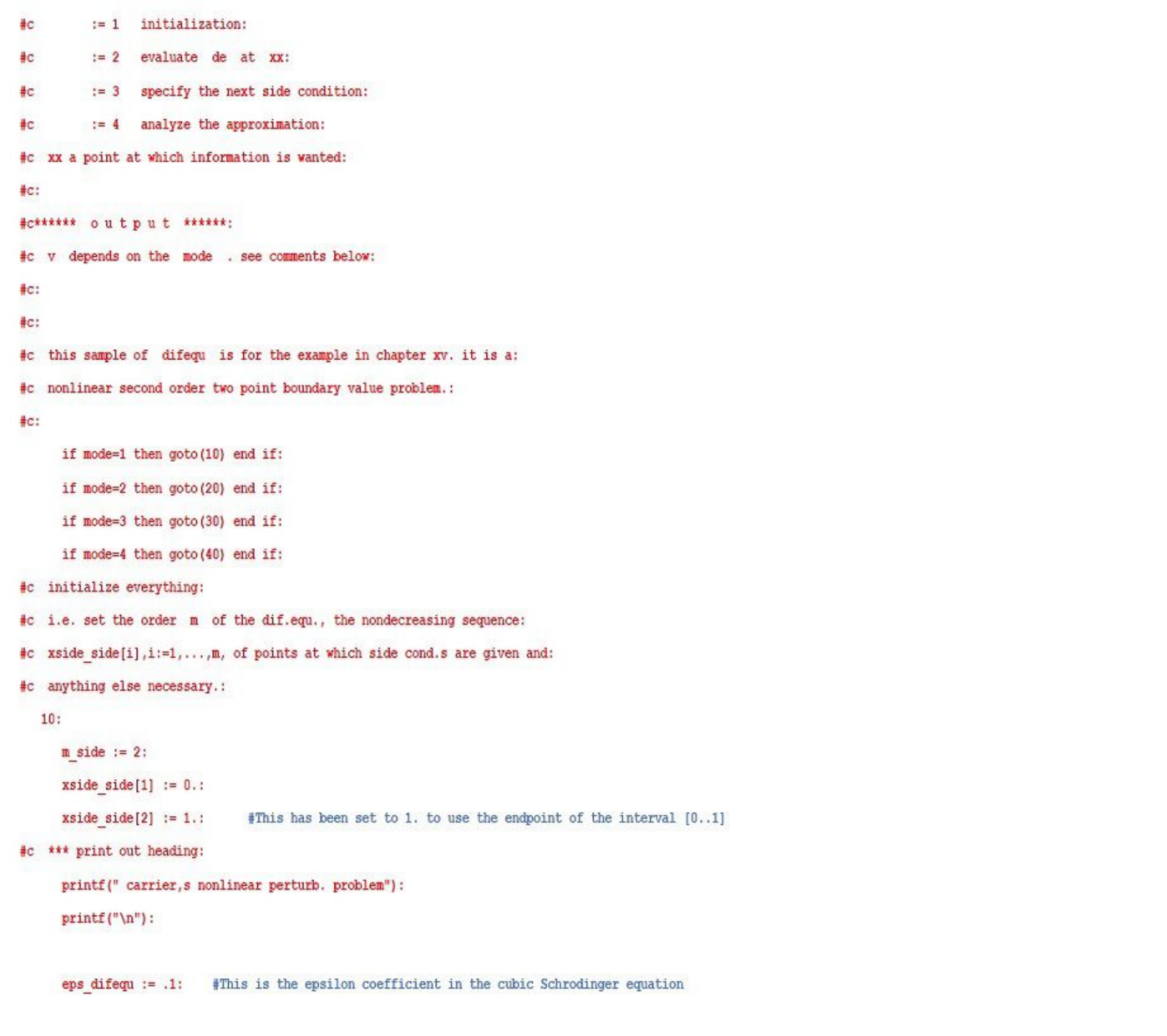}
\newpage
\includegraphics{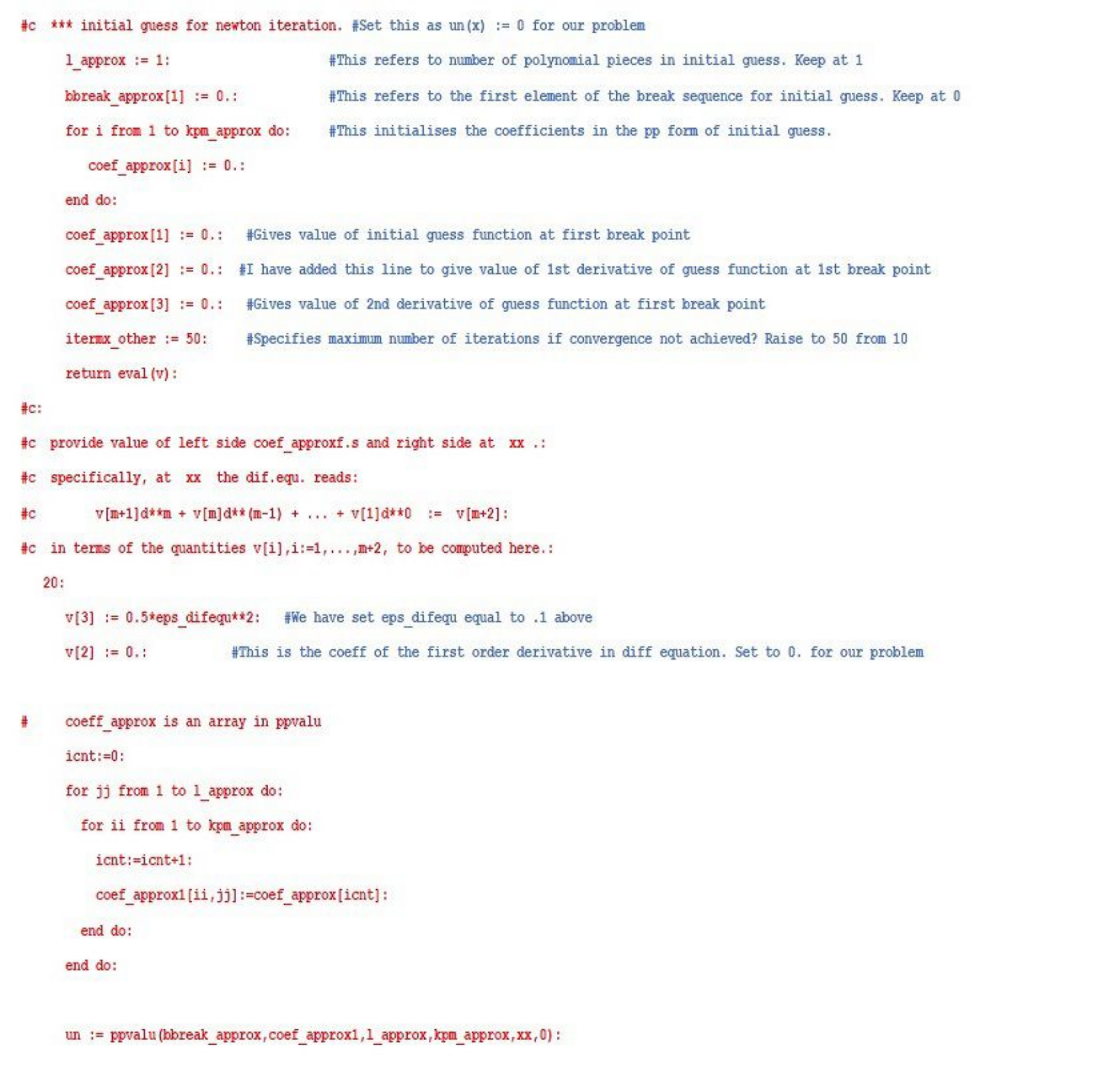}
\newpage
\includegraphics{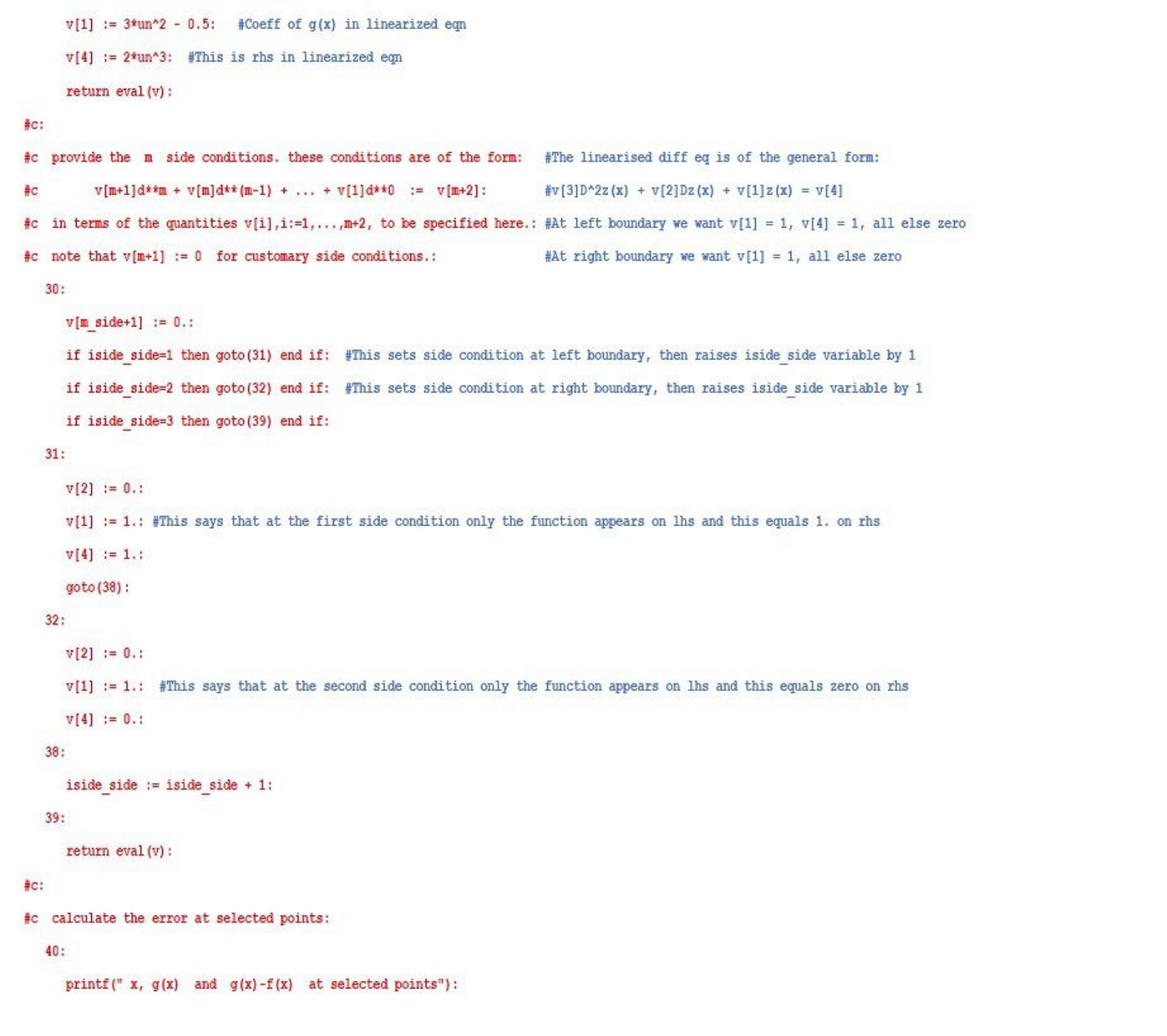}
\newpage
\includegraphics{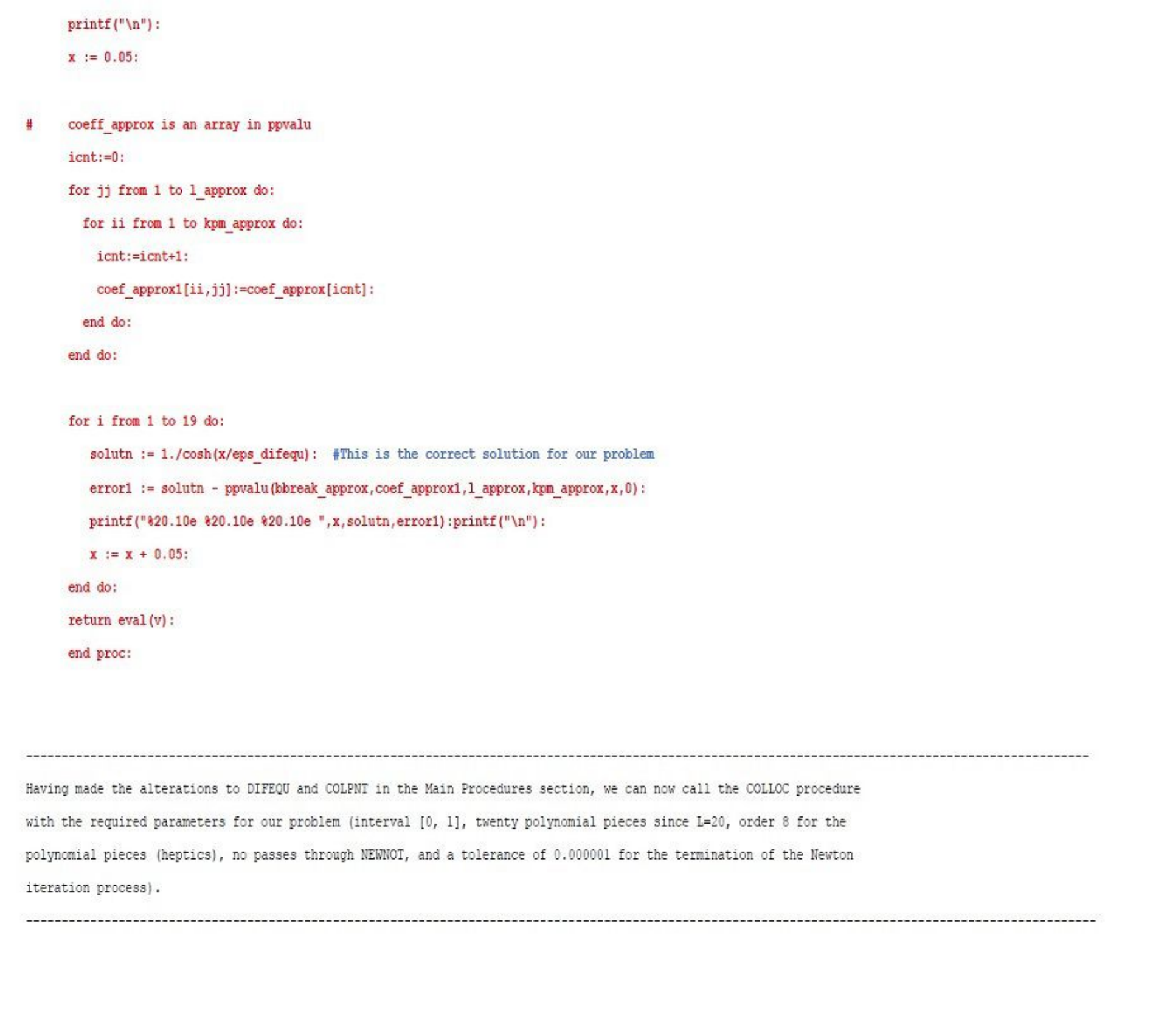}
\newpage
\includegraphics{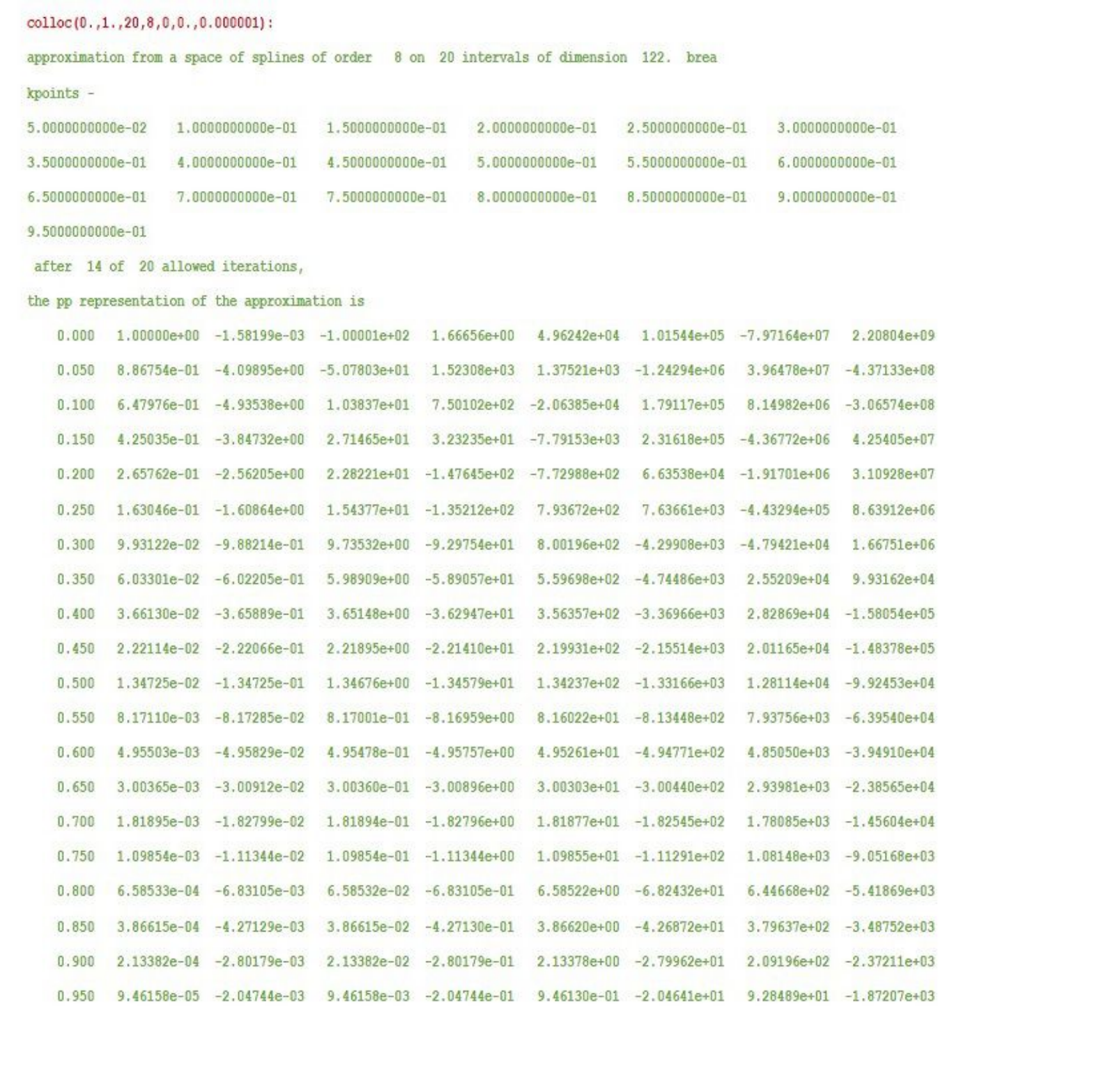}
\newpage
\includegraphics{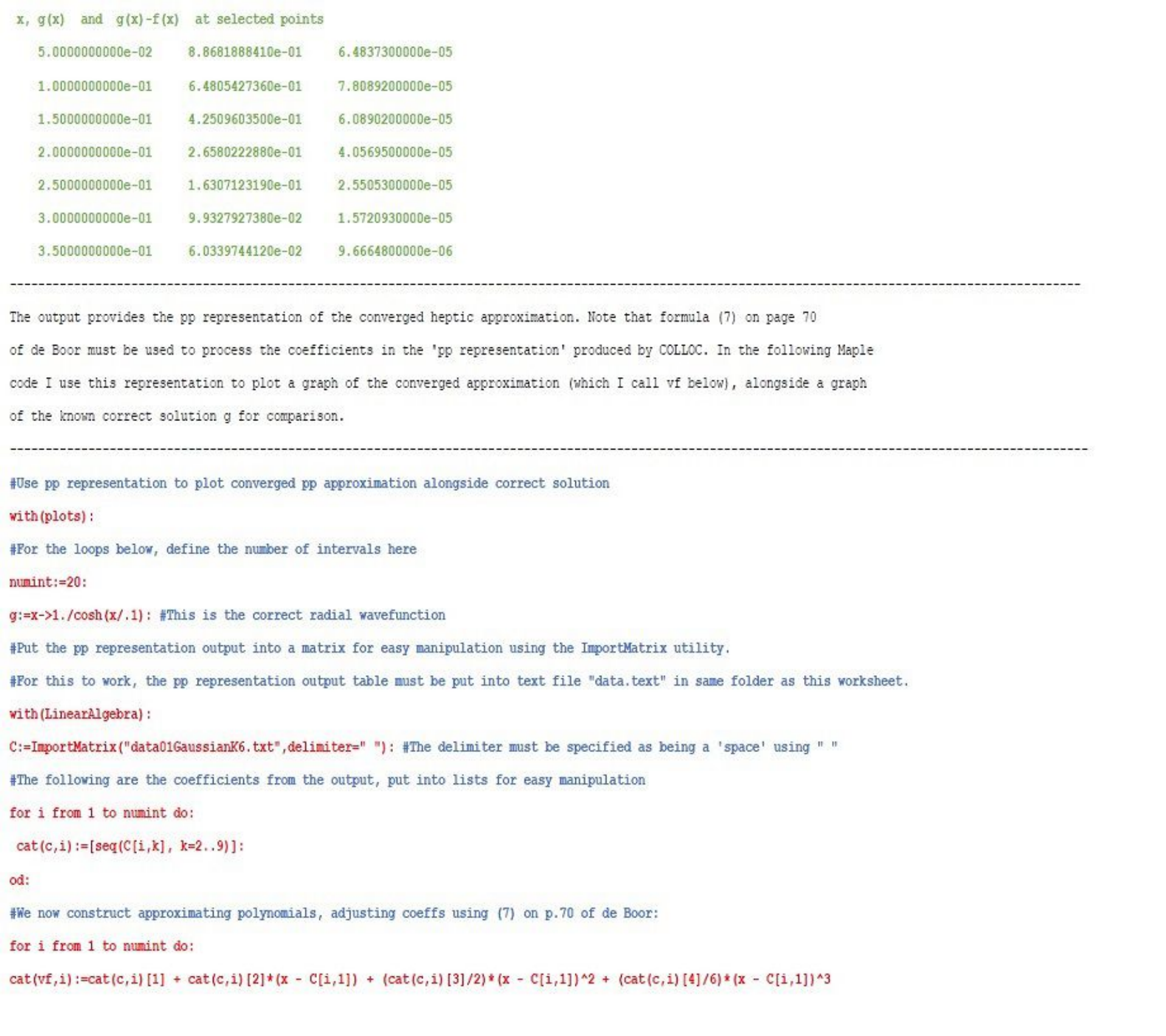}
\newpage
\includegraphics{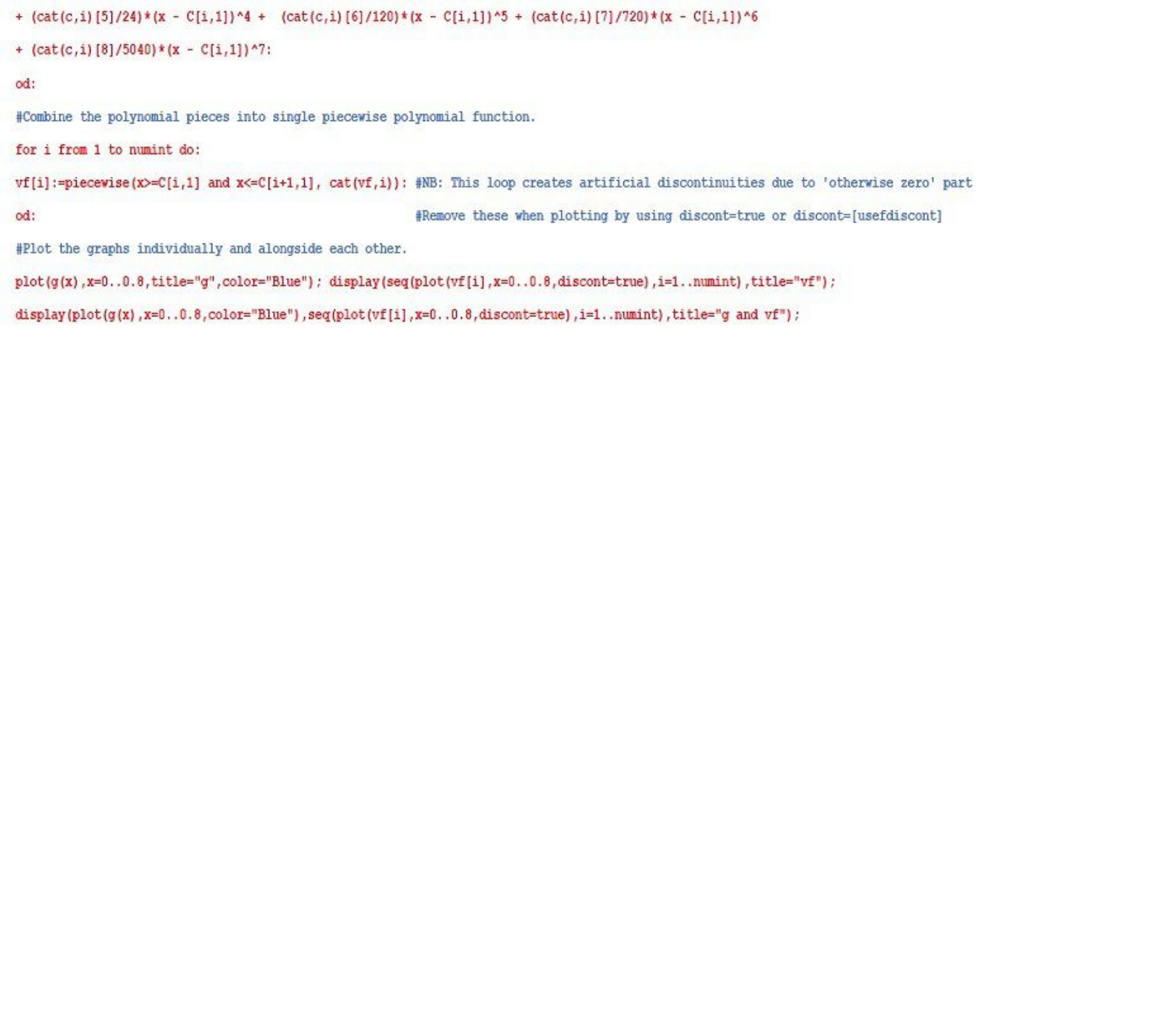}



\addcontentsline{toc}{chapter}{Bibliography}
\begin{thebibliography}{99}
\bibitem{shore} Shore, B, 1973, Solving the radial Schr\"{o}dinger equation by using cubic-spline basis functions, The Journal of Chemical Physics, Vol 58, No 9, pp. 3855-3866.
\bibitem{bottcher} Morrison, J, Bottcher, C, 1993, Spline collocation methods for calculating orbital energies,  Journal of Physics B: Atomic, Molecular and Optical Physics, 26 p. 3999.
\bibitem{sapirstein} Sapirstein, J, Johnson, W, 1996, The use of basis splines in theoretical atomic physics, Journal of Physics B: Atomic, Molecular and Optical Physics, 29 pp. 5213-5225.
\bibitem{morrison} Morrison, J, et al, 1996, Spline collocation calculation for $H_2^{+}$,  Journal of Physics B: Atomic, Molecular and Optical Physics, 29 p. 2375.
\bibitem{sanchez} Sanchez, I, Martin, F, 1997, Representation of the electronic continuum of $\text{H}_2$ with B-spline basis, Journal of Physics B: Atomic, Molecular and Optical Physics, Vol 30, No 3. 
\bibitem{martin} Martin, F, 1999, Ionization and dissociation using B-splines: photoionization of the hydrogen molecule,  Journal of Physics B: Atomic, Molecular and Optical Physics, 32 pp. R197-R231.
\bibitem{odero} Odero, D, Peacher, J, Madison, D, 2001, Numerical solutions of quantum mechanical problems using the basis-spline collocation method, International Journal of Modern Physics C, Vol 12, No 7, pp. 1093-1108.  
\bibitem{bachau} Bachau, H, Cormier, E, Decleva, P, Hansen, J, Martin, F, 2001, Applications of B-splines in atomic and molecular physics, Reports on Progress in Physics, Vol 64, No 12. 
\bibitem{ting} Ting-Yun, S, Cheng-Guang, B, Bai-Wen, L, 2001, Energy spectra of the confined atoms obtained by using B-splines, Communications in Theoretical Physics, Vol 35, No 2. 
\bibitem{chang} Chang, T, Fang, T, 2004, Multiple excitation in photoionization using B-splines, Radiation Physics and Chemistry, Vol 70, pp. 173-190.
\bibitem{deboor} de Boor, C, Swartz, B, 1973, Collocation at Gaussian Points, SIAM Journal on Numerical Analysis, Vol 10, No 4, pp. 582-606. 
\bibitem{deboor2} de Boor, C, 1978, A Practical Guide to Splines (New York: Springer).
\bibitem{beiser} Beiser, A, 2003, Concepts of modern physics, 6th Edition (McGraw-Hill Higher Education, New York). 
\bibitem{brehm} Brehm, J, 1989, Introduction to the structure of matter (John Wiley and sons, Inc.). 
\bibitem{griffiths} Griffiths, D, 1995, Introduction to Quantum Mechanics, 3rd Edition (Prentice Hall). 
\bibitem{jain} Jain, M, 2007, Quantum mechanics: a textbook for undergraduates (PHI Learning Private Limited, New Delhi). 
\bibitem{boas} Boas, M, 2006, Mathematical methods in the physical sciences, 3rd Edition (John Wiley and Sons, Inc.). 
\bibitem{floer} Floer, A, Weinstein, A, 1986, Nonspreading wavepackets for the cubic Schr\"{o}dinger equation with a bounded potential, Journal of Functional Analysis 69, pp. 397-408. 
\bibitem{ebaid} Ebaid, A, Khaled, S, 2011, New types of exact solutions for nonlinear Schr\"{o}dinger equation with cubic nonlinearity, Journal of Computational and Applied Mathematics 235, pp. 1984-1992. 
\bibitem{naumkin} Naumkin, I, 2016, Sharp asymptotic behaviour of solutions for cubic nonlinear Schr\"{o}dinger equations with a potential, Journal of Mathematical Physics 57, 051501.    
\bibitem{omalley} O'Malley, R, 1991, Singular perturbation methods for ordinary differential equations, Applied Mathematical Sciences, 89 (Springer-Verlag, New York).
\bibitem{powell} Powell, M, 1981, Approximation theory and methods (Cambridge University Press).  
\bibitem{rivlin} Rivlin, T, 1981, An introduction to the approximation of functions (Dover).   
\end{thebibliography}
\end{document}